\documentclass[11pt,english,a4paper]{article}

\usepackage{babel}
\usepackage{array}
\usepackage{delarray}
\usepackage{indentfirst}
\usepackage{amsmath}
\usepackage{amsfonts}
\usepackage{latexsym}
\usepackage[dvips]{graphicx}
\usepackage{verbatim}
\usepackage{fontenc}
\usepackage{amssymb}
\usepackage{enumerate}

\input xy
\xyoption{all}

\def\N{\mathbb N}
\def\Z{\mathbb Z}
\def\Q{\mathbb Q}
\def\R{\mathbb R}
\def\C{\mathbb C}

\def\Pun{\mathbb{P}^{1}}

\def\dem{\textbf{Proof. }}
\def\fin{\hspace{0.3cm} $\square$}

\def\rmk{\textbf{Remark. }}

\def\pgcd{\textrm{Gcd}}

\def\CbeA{\mathcal{A}}
\def\CbeC{\mathcal{C}}
\def\CbeD{\mathcal{D}}
\def\CbeE{\mathcal{E}}

\def\CbeH{\mathcal{H}}

\newtheorem{theo}{Theorem}[section]
\newtheorem{stheo}{Theorem}[subsection]
\newtheorem{sstheo}{Theorem}[subsubsection]
 
\newtheorem{slemm}[stheo]{Lemma} 
\newtheorem{sslemm}[sstheo]{Lemma}
\newtheorem{prop}[theo]{Proposition}
\newtheorem{sprop}[stheo]{Proposition}
\newtheorem{ssprop}[sstheo]{Proposition}

\newtheorem{sdefi}[stheo]{Definition}

\newtheorem{coro}[theo]{Corollary}
\newtheorem{scoro}[stheo]{Corollary}
\newtheorem{sscoro}[sstheo]{Corollary}
\newtheorem{nota}[theo]{Notation}
\newtheorem{snota}[stheo]{Notation}
\newtheorem{ssnota}[sstheo]{Notation}
\newtheorem{asum}[theo]{Assumption}

\begin{document}
\bibliographystyle{alpha}
{\LARGE\textbf{Using hyperelliptic curves to find positive polynomials that are not sum of three squares in $\R (x,y)$.}}
\vspace{0.5cm}
\begin{center}
{\large\textbf{Val\'ery Mah\'e}}
\vspace{0.5cm}
$$\begin{array}{rl}
\textbf{adress: }&\textrm{School of mathematics}\\
&\textrm{University of East Anglia}\\
&\textrm{NORWICH NR4 7TJ, UK}\\
&\\
\textbf{email :}&\textrm{V.mahe@uea.ac.uk}\\
\end{array}$$
\end{center}
\vspace{0.5cm}
\noindent\textbf{Abstract: }This article deals with a quantitative aspect of Hilbert's seventeenth problem : producing a collection of real polynomials in two variables of degree 8 in one variable which are positive but are not a sum of three squares of rational fractions.
        As explained by Huisman and Mah\'e, a given monic squarefree positive polynomial in two variables x and y of degree in y divisible by 4 is a sum of three squares of rational fractions if and only if the jabobian variety of some hyperelliptic curve (associated to P) has an "antineutral" point.
        Using this criterium, we follow a method developped by Cassels, Ellison and Pfister to solve our problem : at first we show the Mordell-Weil rank of the jacobian variety J associated to some polynomial is zero (this step is done by doing a 2-descent), and then we check that the jacobian variety J has no antineutral torsion point.

\vspace{0.5cm}
\noindent\textbf{Acknowledgements: } This work is part of my PhD thesis at the universit\'e de Rennes 1 and was financed by the french governement. The editing of this article was supported by EPSRC grant EP/E012590/1. I want to thank warmly my thesis advisors D. Lubicz and L. Mah\'e. I also thank G. Everest and P. Satg\'e for their suggestions and comments.

\tableofcontents

\section*{Introduction}
\addcontentsline{toc}{section}{Introduction}
\markboth{INTRODUCTION}{INTRODUCTION}

Let $P\in\R [X_{1},\cdots ,X_{n}]$ be a polynomial. If $P$ is a sum of squares in $\R (X_{1},\cdots ,X_{n}),$ 
then $P$ is a positive polynomial (i.e. such that $P(x_{1},\cdots ,x_{n})$ is greater than or equal to $0$  for every $(x_{1},\cdots ,x_{n})\in\R^{n}$). 
Conversely, when $P$ is a positive polynomial, 
one can ask whether $P$ is a sum of squares in $\R (X_{1},\cdots ,X_{n}).$ 
That question is called Hilbert's 17th problem and is answered positevely by Artin in the year 1927 (see \cite{Artin}).

A related question is to compute the minimal number $r$ such that every positive polynomial can be written as a sum of $r$ squares in $\R (X_{1},\cdots ,X_{n}).$ The answer is not completely known. Hilbert proves that every positive polynomial $P\in\R [X,Y]$ is a sum of four squares in $\R (X,Y)$ (see \cite{Hilbert}). He proves a little more: every positive polynomial $P\in\R [X,Y]$ of total degree less than $4$ is a sum of three squares of polynomials. The first of these two results is generalized by Pfister in the following way: every positive polynomial $P\in\R [X_{1},\cdots ,X_{n}]$ is a sum of $2^{n}$ squares in $\R (X_{1},\cdots ,X_{n})$ (see \cite{Pfister1967}).%

There is no known effective characterization of sums of three squares in $\R (X,Y).$ 
However, in 1971, Cassels, Ellison and Pfister show that 
%
%
Motzkin's polynomial 
$$M(X,Y)=1+X^{2}Y^{4}+X^{4}Y^{2}-3X^{2}Y^{2}$$ 
is a positive polynomial (thus a sum of four squares in $\R (X,Y)$), but is not a sum of three squares in $\R (X,Y)$ (see \cite{Cassels-Ellison-Pfister}).

To prove this theorem, Cassels, Ellison and Pfister consider for each positive polynomial $F(X,Y)=1+A(X)Y^{2}+B(X)Y^{4}\in\R (X,Y)$ with $B(A^{2}-4B)\neq 0$ the elliptic curve $\CbeE_{F}$ 
defined over $\R(x)$ by the equation 
$$-\beta^{2}=\alpha (\alpha^{2}-2A(x)\alpha +A(x)^{2}-4B(x))$$ 
and they show that $F$ is a sum of three squares in $\mathbb{R}(X,Y)$ if and only if 
$\CbeE_{F}$ has a $\R (x)$-point $(\alpha ,\beta )$ such that $\alpha$ and $-\left(\alpha^{2}-2A\left( x\right)\alpha +A\left( x\right)^{2}-4B\left( x\right)\right)$ are sums of two squares in $\R (x)$ (i.e. take only nonnegative values on $\R ).$ 
A similar method allows Christie (in $1976$, see \cite{Christie}), then Mac\'e (in $2000$, see \cite{Mace}), and then Mac\'e and Mah\'e (in 2005, see \cite{Mahe-Mace2005}) to construct other families of positive polynomials in two variables that are not a sum of three squares of rational fractions. 




Using a totally different strategy (based on Noether-Lefschetz's theorem), Colliot-Th\'el\`ene proves in 1992 the existence of positive polynomials in two variables of even degree greater than or equal to $6$ in one variable that are not a sum of three squares of rational fractions. (see \cite{Colliot1993}).

Using the method of Cassels, Ellison and Pfister one can only study  polynomials of the form  $Y^{4}+A(X)Y^{2}+B(X)$: we need the elliptic curve $\CbeE_{F}$ given by the equation $-\beta^{2}=\alpha (\alpha^{2}-2A(x)\alpha +A(x)^{2}-4B(x)).$ During the year 2001, Huisman and Mah\'e generalized the construction of Cassels, Ellison and Pfister by introducing the concept of antineutral point (see Definition \ref{21-def-pt-antineutre}). %
In \cite{Huisman-Mahe2001}, Huisman and Mah\'e show that a nonconstant monic squarefree polynomial $P(X,Y)$ of degree in $Y$ divisible by $4$ is a sum of three squares in the field $\mathbb{R}(X,Y)$ if and only if a $\R (x)$-point of the jacobian of the hyperelliptic curve $\CbeC$ 
defined over $\R (x)$ by the equation $z^{2}+P(x,y)=0$ is antineutral.%

In this article we generalize the method of Cassels, Ellison and Pfister using the results of Huisman and Mah\'e in order to construct families of positive polynomials in two variables of degree $8$ in one variable that are not a sum of three squares. As a corollary we get a positive polynomial with coefficients in $\Q$ of degree $8$ in one variable that is not a sum of three squares (such an example was not known before).

\begin{section}{Notation.}

To simplify our statements we consider an elliptic curves as a genus $1$ hyperelliptic curve (i.e. we do not assume a hyperelliptic curve to have genus at least $2$). 


For the background on Mumford's representation, semi-reduced divisors, reduced divisors ans Cantor's algorithms we refer to \cite{Mumford} and \cite{Cantor} (see also \cite{Cohen-Frey-al} and \cite{Gaudry}). A semi-reduced divisor with Mumford's representation $(u,v)$ is denoted by $\textrm{div}(u,v)$. A linear equivalence class with Mumford's representation $(u,v)$ is denoted by $<u,v>$.

When $D$ is a divisor on a curve $\CbeC$ defined over a field $k$ and $K$ is an extension of $k$ we denote by $\textrm{Supp}_{K}(D)$ the support of $D$ as a divisor on $\CbeC\times_{k}K$.

For every abelian group $A$ and for every $n\in\N^{*}$ we denote 
\begin{itemize}
\item[* ] by $[n]_{A}$ (or $[n]$) the multiplication-by-$n$ automorphism of $A$,
\item[* ] by $A[n]$ the kernel of $[n]_{A}$ and
\item[* ] by $A_{tors}$ the torsion subgroup of $A$.
\end{itemize}
For the background on places of function fields we refer to \cite{Stichtenoth} (in this article we use its notation; in particular by a function field over a field $k$ we mean a transcendance degree $1$ extension of $k$). If $F_{1}$ is a function field and $F_{2}/F_{1}$ is a finite extension and $\mathcal{P}$ is a place of $F_{2}$ above a place $\mathfrak{p}$ of $F_{1}$ we denote by $e(\mathcal{P}|\mathfrak{p} )$ the ramification index of $\mathcal{P}$ above $\mathfrak{p}$ and by $f(\mathcal{P}|\mathfrak{p} )$ the relative degree. 

\noindent\rmk%
In this article when we refer to a hypothesis by giving its number, we mean the corresponding hypothesis in either Assumptions \ref{Asummption23} or Assumptions \ref{Asummption24}

\end{section}

\begin{section}{Statement of the results.}


\begin{nota}\label{Notations_Principales_21} Let $\eta,$ $\omega$ and $\rho$ be real numbers. We assume that $|\omega |$ and $|\eta |$ are distinct. Denote by $b_{1}$ the element 
$$b_{1}:=\frac{\rho^{2}-\eta^{2}}{\omega^{2}-\eta^{2}}+\frac{\eta^{2}-\omega^{2}}{4}.$$
In this article we look at the polynomial 
$$P(x^{2},y^{2}):=\left( y^{2}+1\right)\left( y^{2}+C\left( x^{2}\right)\right)\left( y^{4}+\left( 1+C\left( x^{2}\right)\right) y^{2}+B\left( x^{2}\right)\right)$$
with $B(x):=(x+b_{1})^{2}-\eta^{2}$ and $C(x):=2(x+b_{1})+\omega^{2}-\eta^{2}-1.$ This polynomial is defined on the field $k:=\Q (\eta ,\omega ,\rho )$. We denote by $\CbeC$ the hyperelliptic curve defined over $k(x)$ by the affine equation $z^{2}+P(x^{2},y^{2})=0$.
\end{nota}

\begin{asum}\label{Asummption22} We assume the following three inequalities 
$$\omega >1+|\eta |,\hspace{0.5cm}\omega^{2}-\eta^{2}>2\omega\hspace{0.2cm}\textrm{and}\hspace{0.2cm}b_{1}>1+\frac{\omega^{2}-\eta^{2}}{2}.$$
\end{asum}

\begin{asum}\label{Asummption23} We assume that all the following elements are different from $0$
\begin{enumerate}
\item\label{hypothesis1} $\eta$ and $\rho$, 
\item\label{hypothesis2} $\omega^{2}-\eta^{2}-2-2\eta$ and $\omega^{2}-\eta^{2}-2+2\eta$,
\item\label{hypothesis3} $\left(\omega^{2}-\eta^{2}-2\right)^{2}-4\eta^{2}-4,$
\item\label{hypothesis4} $\omega^{2}-\eta^{2}-1+2\eta$ and  $\omega^{2}-\eta^{2}-1-2\eta$,
\item\label{hypothesis5} $\left(\omega^{2}-\eta^{2}-1\right)^{2}-4\eta^{2}-1$,
\item\label{hypothesis6} $\omega^{2}-\eta^{2}-2\eta$ and $\omega^{2}-\eta^{2}+2\eta$.
\end{enumerate}
\end{asum}

\begin{asum}\label{Asummption24} We assume that none of the following elements is a square in $k$:
\begin{enumerate}[a.]
\item\label{hypothesisA} $\left(\left(\omega^{2}-\eta^{2}\right)^{2}-4\omega^{2}\right)=\left(\omega^{2}-\eta^{2}-2\omega\right)\left(\omega^{2}-\eta^{2}+2\omega\right)$,
\item\label{hypothesisB} $(2b_{1}-2+\omega^{2}-\eta^{2})(\omega^{2}-\eta^{2}-2\omega )$,
\item\label{hypothesisC} $(2b_{1}-2+\omega^{2}-\eta^{2})(\omega^{2}-\eta^{2}+2\omega )$,
\item\label{hypothesisDsupplementaire} $2(\omega^{2}-\eta^{2}-2\omega )(b_{1}-1-\omega )$,
\item\label{hypothesisEsupplementaire} $2(\omega^{2}-\eta^{2}+2\omega )(b_{1}-1+\omega )$
\item\label{hypothesisD} $2\left( 2b_{1}-2+\omega^{2}-\eta^{2}\right)\left( b_{1}-1+\omega\right)$, 
\item\label{hypothesisE} $2\left( 2b_{1}-2+\omega^{2}-\eta^{2}\right)\left( b_{1}-1-\omega\right)$, 
\item\label{hypothesisF} $\left(\left( b_{1}-1\right)^{2}-\omega^{2}\right)\left(\left(\omega^{2}-\eta^{2}\right)^{2}-4\omega^{2}\right)$ 
\item\label{hypothesisG} $2\left(\omega^{2}-\eta^{2}\right)\left(\omega^{2}-\eta^{2}-2\omega\right)\left(\left(\omega +1\right)^{2}-\eta^{2}\right)^{n}$ (for each $n\in\{ 0,1\}$),
\item\label{hypothesisH} $2\left(\omega^{2}-\eta^{2}\right)\left(\omega^{2}-\eta^{2}+2\omega\right)\left(\left(\omega -1\right)^{2}-\eta^{2}\right)^{n}$ (for each $n\in\{ 0,1\}$) 
\item\label{hypothesisI} $\left(\left( b_{1}-1\right)^{2}-\omega^{2}\right)^{n_{1}}\left(\left(\omega -1\right)^{2}-\eta^{2}\right)^{n_{2}}\left(\left(\omega +1\right)^{2}-\eta^{2}\right)^{n_{3}}$ (with $(n_{1},n_{2},n_{3})$ a nontrivial triplet of elements of $\{ 0,1\}$),
\item\label{hypothesisJ} $2\left(\omega^{2}-\eta^{2}\right)\left(2b_{1}-2+\omega^{2}-\eta^{2}\right)$, 
\item\label{hypothesisK} $\begin{array}[t]{l}
2^{n_{1}}\left(\omega^{2}-\eta^{2}+2\omega\right)\left( b_{1}-1+\omega\right)\left(\omega^{2}-\eta^{2}\right)^{n_{1}}\left( 2b_{1}-2+\omega^{2}-\eta^{2}\right)^{1-n_{1}}\\
\times\left(\omega -1-\eta\right)^{1-n_{2}}\left(\omega -1+\eta\right)^{n_{2}}\\
\end{array}$\\
(with $n_{1}$, $n_{2}\in\N$),
\item\label{hypothesisL} $\begin{array}[t]{l}
2^{1-n_{1}}\left(\omega^{2}-\eta^{2}+2\omega\right)\left(\omega^{2}-\eta^{2}\right)^{n_{1}}\left( 2b_{1}-2+\omega^{2}-\eta^{2}\right)^{n_{1}}\\
\times\left(\omega -1-\eta\right)^{1-n_{2}}\left(\omega -1+\eta\right)^{n_{2}},\\
\end{array}$ \\
(with $n_{1}$, $n_{2}\in\N$),
\item\label{hypothesisM} $b_{1}^{2}-\eta^{2}$, and
\item\label{hypothesisN} $2b_{1}+\omega^{2}-\eta^{2}-1$.
\end{enumerate}
\end{asum}

\begin{theo}\label{theo-final-pas-somme-3-carres-46}  We use the notation in \ref{Notations_Principales_21}. Then, under assumptions \ref{Asummption22}, \ref{Asummption23} and \ref{Asummption24} the polynomial 
$$P(x^{2},y^{2})=\left( y^{2}+1\right)\left( y^{2}+C\left( x^{2}\right)\right)\left( y^{4}+\left( 1+C\left( x^{2}\right)\right) y^{2}+B\left( x^{2}\right)\right)$$
is positive but is not a sum of three squares in $\R (x,y)$.
\end{theo}

\begin{coro} We use the notation in \ref{Notations_Principales_21}. Under Assumptions \ref{Asummption22}, %
%
if $\eta,$ $\omega ,$ and $\rho$ are algebraically independant over $\Q$, then the polynomial 
$$P(x^{2},y^{2})=\left( y^{2}+1\right)\left( y^{2}+C\left( x^{2}\right)\right)\left( y^{4}+\left( 1+C\left( x^{2}\right)\right) y^{2}+B\left( x^{2}\right)\right)$$
is positive but is not a sum of three squares in $\R (x,y)$.
\end{coro}

\begin{coro}  
Consider the two polynomials 
$$B(x):=
x^{2}+\frac{14063}{22}x+\frac{196743825}{1936}\textrm{ and }
C(x):=
2x+\frac{27835}{22}.$$ 
Then the positive polynomial 
$$P(x^{2},y^{2}):=\left( y^{2}+1\right)\left( y^{2}+C\left( x^{2}\right)\right)\left( y^{4}+\left( 1+C\left( x^{2}\right)\right) y^{2}+B\left( x^{2}\right)\right)\in\Q (x,y)$$
is not a sum of three squares in $\R (x,y)$.
\end{coro}
\dem Apply Theorem \ref{theo-final-pas-somme-3-carres-46} with $\eta :=23$ $\omega :=34$ and $\rho :=547$. 
\fin%
\\


In \cite{Pfister1967}, Pfister showed the product of two sums of $2^{n}$ squares is a sum of $2^{n}$ squares. In general a product of two sums of three squares is not a sum of three squares. Looking for antineutral torsion points  we can give examples of products of four sums of three squares in $\R (x,y)$ that are a sum of three squares of polynomials. 
      \begin{prop}\label{formule-4-tors-non-triviale} Let $\alpha ,$ $\beta ,$ and $\gamma\in\R (x)$ be three rational fractions. Denote %
$a:=1+\alpha^{2}(1+\beta^{2})(1+\gamma^{2}),$ $b:=1+\alpha^{2}(1+\beta^{2})^{2}(1+\gamma^{2})$ and \linebreak$c:=1+\alpha^{2}(1+\beta^{2})(1+\gamma^{2})^{2}.$

        Then the polynomial $P(x,y):=(y^{2}+1)(y^{2}+a)(y^{2}+b)(y^{2}+c)$ is a sum of three squares in $\R (x,y)$:
        $$\begin{array}{rcl}
          P(x,y)&=&\left( \frac{(a-1)y(y^{2}+a)}{\alpha}+\alpha\beta\gamma ((1-\beta\gamma )y+\beta +\gamma )(y^{2}+1)\right)^{2}\\
          &&+\left( \frac{(a-1)(y^{2}+a)}{\alpha}+\alpha\beta\gamma (1-\beta\gamma -(\beta +\gamma )y)(y^{2}+1)\right)^{2}\\
          &&+\left( (y^{2}+1)(y^{2}+a-\beta\gamma (a-1))\right)^{2}.
        \end{array}$$
      \end{prop}
      For more details about Proposition \ref{formule-4-tors-non-triviale} we refer to subsection \ref{ref-heuristic-prod-sum3squares}.

\end{section}

\begin{section}{Sums of three squares and antineutral points.}

\begin{subsection}{The results of Huisman and Mah\'e}

\begin{snota}\label{HMnot-anti1}
\textrm{Let $\Sigma$ be the galois group $\textrm{Gal}(\C (x)/\R (x))=\textrm{Gal}(\C /\R )$ and $\sigma$ be its nontrivial element. Let $\fbox{2}_{\R (x)}$ be the group of nontrivial sums of two squares in $\R (x)$.} 

\textrm{Let $\CbeD$ be a geometrically integral, smooth, projective curve over $\R (x)$ with odd genus. Let $\CbeD ':=\CbeD\times_{\R (x)}\C (x)$ be its complexification and $p:\CbeD '\longrightarrow\CbeD$ be the projection. The galois group $\Sigma$ acts naturally on $\CbeD '$. This action induces an action of $\Sigma$ on the Picard group $\textrm{Pic}(\CbeD ')$.} 

\textrm{The projection $p$ induces a morphism $p^{*}$ from $\textrm{Pic}(\CbeD )$ into $\textrm{Pic}(\CbeD ')$. The image of $p^{*}$ is contained in the subgroup $\textrm{Pic}(\CbeD ')^{\Sigma}$ of $\Sigma$-invariants elements of $\textrm{Pic}(\CbeD ').$}
\end{snota}

\begin{snota}\label{HMnot-anti1bis}We define a group homomorphism 
$$\delta :\textrm{Pic}(\CbeD ')^{\Sigma}\longrightarrow H^{1}(\Sigma , \C (x)(\CbeD ')^{\times}/\C (x)^{\times}).$$ Let $cl(A)\in\textrm{Pic}(\CbeD ')^{\Sigma}$ be the class of a divisor $A$. 
Because of the $\Sigma$-invariance of $cl(A),$ the divisor $A-\sigma^{*}A$ is the principal divisor associated to a function $f\in \C (x)(\CbeD ')^{\times}.$ The principal divisor of $\R (x)(\CbeD)$ associated to $f\sigma (f)$ is $0.$ Thus $f\sigma (f)$ is an element of $\R (x)^{\times}.$ The element $\delta (cl(A))$ is chosen as the class of $f$ in $H^{1}(\Sigma , \C (x)(\CbeD ')^{\times}/\C (x)^{\times})$.
\end{snota}

\begin{slemm}\label{suite_exacte_antineutralite} Using the notation in \ref{HMnot-anti1} and \ref{HMnot-anti1bis}, the following is an exact sequence:
$$\xymatrix{ 0 \ar[r] & \textrm{Pic}(\CbeD ) \ar[r]^-{p^{*}} & \textrm{Pic}(\CbeD ')^{\Sigma} \ar[r]^-{\delta} & H^{1}(\Sigma , \C (x)(\CbeD ')^{\times}/\C (x)^{\times}) \ar[r] & 0. \\}$$
\end{slemm}
\rmk%
The map $\delta$ is a coboundary map; it can be defined by looking at the long exact sequence associated to the short exact sequence :
$$\xymatrix{ 0 \ar[r] & \C (x)(\CbeD ')^{\times}/\C (x)^{\times} \ar[r]^-{\textrm{div}} & \textrm{Div}(\CbeD ') \ar[r]^-{\textrm{cl}} & \textrm{Pic}(\CbeD ') \ar[r] & 0. \\}$$

\begin{snota}\label{not-anti2}
We use the notation of Lemma \ref{suite_exacte_antineutralite}. The map\linebreak $1+\sigma : \C (x)(\CbeD ')\longrightarrow \R (x)(\CbeD )$ induces a monomorphism 
$$\eta: H^{1}(\Sigma , \C (x)(\CbeD ')^{\times}/\C (x)^{\times})\longrightarrow \R (x)^{\times}/\fbox{2}_{\R (x)}.$$ 
We denote by $\varpi :\textrm{Pic}^{0}(\CbeD ')^{\Sigma}\longrightarrow\R (x)^{\times}/\fbox{2}_{\R (x)}$ the restriction of the map $\eta\circ\delta$ to $\textrm{Pic}^{0}(\CbeD ')^{\Sigma}$.
\end{snota}
\begin{sdefi}\label{21-def-pt-antineutre}
Let us use the notation in \ref{not-anti2}. An element $\beta\in\textrm{Pic}^{0}(\CbeD ')^{\Sigma}$ is said to be antineutral when $\varpi (\beta)=-1.$%
%
\end{sdefi}
\begin{sprop}\label{prop-finale-HM}
Let $P(Y)\in \R (x)[Y]$ be a squarefree nonconstant monic totally positive polynomial of degree divisible by $4$. Let $\CbeD$ be 
the hyperelliptic curve defined over $\R (x)$ by the affine equation $z^{2}+P(y)=0$. We use the notation in \ref{not-anti2} (relative to $\CbeD$). 
Then $P(Y)$ is a sum of three squares in $\R (x,Y)$ if and only if 
$\textrm{Pic}^{0}(\CbeD ')^{\Sigma}$ has an antineutral element.
\end{sprop}
\dem%
See \cite{Huisman-Mahe2001} Theorem 6.5.
\fin%

\end{subsection}

\begin{subsection}{An effective version.}

\begin{snota}\label{nota-22-caract-antin1}
Let $k$ be a subfield of $\R .$ Let $k'$ be the field $k(i)$. Let $\sigma$ be the nontrivial element of the galois group $\textrm{Gal}(k'(x)/k(x))=\textrm{Gal}(k'/k).$ Let $Q\in k(x)[y]$ be a monic polynomial such that $(y^{2}+1)Q(y^{2})$ is squarefree. Let $\CbeC$ be the 
hyperelliptic curve defined over $k(x)$ by the affine equation 
$$\CbeC :z^{2}+(y^{2}+1)Q(y^{2})=0$$
and let $\CbeC ':=\CbeC\times_{k(x)}k'(x)$ be its complexification. 

Let $g$ be the degree of the polynomial $Q$. We assume that $g$ is odd and that the rational fraction $d:= -Q(-1)\in k(x)$ is nontrivial. 
Let $\widetilde{\CbeC}$ be the 
$k(x)$-hyperelliptic curve given in coordinates $(s,t)$ by the affine equation
$$\widetilde{\CbeC}:t^{2}=-\frac{s}{d}(s-d)^{2g}Q\left(-\left(\frac{s+d}{s-d}\right)^{2}\right)$$ 
and let $\widetilde{\CbeC}':=\widetilde{\CbeC}\times_{k(x)}k'(x)$ be its complexification. The two curves $\CbeC '$ and $\widetilde{\CbeC}'$ have a $\C (x)$-rational point.
\end{snota}

\begin{snota}\label{nota-22-caract-antin2}
The map $\gamma : \begin{array}[t]{ccc}
k'(x)(\CbeC ') & \longrightarrow &  k'(x)(\widetilde{\CbeC}')\\
A(y,z) & \longmapsto & A\left( i\frac{s+d}{s-d},\frac{2idt}{(s-d)^{g+1}}\right)\\
\end{array}$ is an isomorphism. 
Denote by $\omega$ the $k'(x)$-automorphism $\sigma^{-1}\circ\gamma\circ\sigma\circ\gamma^{-1}$ of $k'(x)(\widetilde{\CbeC}')$. It sends $s$ to $\frac{d^{2}}{s}$ and $t$ to $(-1)^{g}\frac{d^{g+1}t}{s^{g+1}}$. The $k'(x)$-automorphism $\omega$ induces a $k'(x)$-automorphism $\widetilde{\omega}$ of $\widetilde{\CbeC}'$ and a $k'(x)$-automorphism $\Omega$ of $\textrm{Jac}(\widetilde{\CbeC}')$.

\end{snota}

\noindent\rmk%
The curve $\CbeC '$ has two $\C (x)$-rational Weierstrass points: \linebreak$(s,t)=(i,0)$ and $(s,t)=(-i,0)$. The map $\gamma$ is obtained by considering a map from $\CbeC '$ to $\widetilde{\CbeC}'$ that sends $(i,0)$ to infinity and $(-i,0)$ to $(0,0)$.\\

\noindent\rmk%
The degree of the polynomial $\frac{s}{d}(s-d)^{2g}Q\left(-\left(\frac{s+d}{s-d}\right)^{2}\right)$ is odd. Thus Mumford representation can be used to compute in the group $\textrm{Jac}(\widetilde{\CbeC}')(\C (x))$ (which can be identified to the group $\textrm{Pic}^{0}(\widetilde{\CbeC}')$). \\

\noindent\rmk%
The notion of antineutral point involves an action of $\Sigma$ on \linebreak$\textrm{Pic}^{0}(\CbeC ')$. 
This action gives an action of $\Sigma$ on $\textrm{Jac}(\widetilde{\CbeC}')\simeq\textrm{Jac}(\CbeC ')$ which is not the action induced by the natural action of $\Sigma$ on $\widetilde{\CbeC}\times_{k(x)}k'(x)$ but 
is the action obtained by consi\-dering the natural action of $\Sigma$ on $\CbeC '=\CbeC\times_{k(x)}k'(x)$ and by looking at $\CbeC$ as a $k'(x)/k(x)$-form of $\widetilde{\CbeC}$. In other words we study the action of $\Sigma$ on $\textrm{Jac}(\widetilde{\CbeC}')$ associated to the the $1$-cocycle $\Sigma\longrightarrow\textrm{Aut}_{k'(x)}(\textrm{Jac}(\widetilde{\CbeC}'))$ whose value at $\sigma$ is $\Omega$ (see \cite{Raynaud}). 



\begin{slemm}\label{lemm-image-diviseur-conj-complexe-22-antin}We use notation \ref{nota-22-caract-antin1} and \ref{nota-22-caract-antin2}. By linearity we extend $\widetilde{\omega}$ to the group $\textrm{Div}^{0}(\widetilde{\CbeC}')$.

Let $D=\textrm{div}(u,v)\in\textrm{Div}^{0}(\widetilde{\CbeC}')$ be a reduced divisor. Let $m$ be the degree of $u$. We assume that $u(0)$ is nontrivial. 
Let $e$ be the quotient of the Euclidean division of $m+1=\deg_{s}(u)+1$ by $2$ and $\epsilon$ be its remainder.

Then the divisor $\widetilde{\omega}(D)+\textrm{div}(s^{e})$ is semi-reduced and $\widetilde{\omega}(D)+\textrm{div}(s^{e})=\textrm{div}\left(\frac{1}{u(0)}s^{2e}u\left(\frac{d^{2}}{s}\right) ,\hat{v}\right)$ where $\hat{v}$ denotes the remainder of the Euclidean division of $(-1)^{g}\left(\frac{s}{d}\right)^{g+1}v(\frac{d^{2}}{s})$ by $s^{2e}u(\frac{d^{2}}{s}).$
\end{slemm}
\dem%
Denote by $\infty$ the unique point at infinity of $\widetilde{\CbeC}'$. Denote by $\mathcal{P}_{0}$ the point $(0,0)$. The point $\mathcal{P}_{0}$ does not belong to $\textrm{Supp}(D).$ Write $D$ as $D=\left(\displaystyle\sum_{i=1}^{r}n_{i}\mathcal{P}_{i}\right) -m\infty$ where $\mathcal{P}_{i}$ is a point of $\widetilde{\CbeC}'$ (notice that $m=\displaystyle\sum_{i=1}^{r}n_{i}\deg (\mathcal{P}_{i})$). The divisor $\widetilde{\omega}(D)+\textrm{div}(s^{e})$ is semi-reduced because $\widetilde{\omega}(D)$ is equal to 
$$\left(\displaystyle\sum_{i=1}^{r}n_{i}\widetilde{\omega}(\mathcal{P}_{i})\right) -m\mathcal{P}_{0} =\left(\displaystyle\sum_{i=1}^{r}n_{i}\widetilde{\omega}(\mathcal{P}_{i})\right) -m\infty -\textrm{div}(s^{e})+(1-\epsilon )\left(\mathcal{P}_{0}-\infty\right) .$$
Let $p_{s}:\widetilde{\CbeC}'\longrightarrow\mathbb{P}^{1}_{k'(x)}$ be the projection associated to the $s$-coordinate. By definition of Mumford's representation, $p_{s}(D)$ is the principal divisor associated to the function $u$. 

Let $\check{\omega}:\mathbb{P}^{1}_{k'(x)}\longrightarrow\mathbb{P}^{1}_{k'(x)}$ be the morphism associated to the $k'(x)$-automorphism of $k'(x)(s)$ that sends $s$ on $\frac{d^{2}}{s}$. We know that $p_{s}\circ\widetilde{\omega}=\check{\omega}\circ p_{s}$ (look at the associated equality on the associated function fields). By linearity this equality is still valid for divisors: we have $p_{s}(\widetilde{\omega}(D)=\check{\omega}( p_{s}(D)).$ In particular Mumford's representation of $\widetilde{\omega}(D)+\textrm{div}(s^{e})$ is of the form $\left(\frac{1}{u(0)}s^{2e}u\left(\frac{d^{2}}{s}\right) ,w\right)$ for some polynomial $w\in k'(x)[y]$.

The function $t-v$ vanishes at all elements of the support of $D$ that are distinct from $\infty$. Since $\widetilde{\omega}\circ\textrm{div}=\textrm{div}\circ\omega$ and since $2e=m+1-\epsilon\le g+1$, this implies that  
$(-1)^{g}\frac{s^{g+1}}{d^{g+1}}\omega(t-v)=\left(t-(-1)^{g}\left(\frac{s}{d}\right)^{g+1}v(\frac{d^{2}}{s})\right)$ vanishes at 
elements of the support of $\widetilde{\omega}(D)+\textrm{div}(s^{e})$. In particular the remainder $\widehat{v}$ of the Euclidean division of $(-1)^{g}\left(\frac{s}{d}\right)^{g+1}v(\frac{d^{2}}{s})$ by $s^{2e}u(\frac{d^{2}}{s})$ is such that $t-\widehat{v}$ vanishes all elements of the support of $\widetilde{\omega}(D)+\textrm{div}(s^{e})$. 
%

To check that $\left(\frac{1}{u(0)}s^{2e}u\left(\frac{d^{2}}{s}\right) ,\hat{v}\right)$ satisfies the definition of Mumford's representation for $\widetilde{\omega}(D)+\textrm{div}(s^{e})$, we only need to prove the divisibility of $\hat{v}^{2}+\frac{s}{d}(s-d)^{2g}Q\left(-\left(\frac{s+d}{s-d}\right)^{2}\right)$ by $\frac{1}{u(0)}s^{2e}u\left(\frac{d^{2}}{s}\right)$. This relation of divisibility is obtained by applying $\omega$ to the relation of divisibility of \linebreak$v^{2}+\frac{s}{d}(s-d)^{2g}Q\left(-\left(\frac{s+d}{s-d}\right)^{2}\right)$ by $u$.
\fin%

\begin{sprop}\label{theo-caracterisation-antineutre-final}
We use the notation \ref{not-anti2} and \ref{nota-22-caract-antin1} and \ref{nota-22-caract-antin2}. Denote by $\sigma_{\star\CbeC}$ (respectively $\sigma_{\star\widetilde{\CbeC}}$) the action of $\sigma$ on $\textrm{Jac}(\widetilde{\CbeC}')$ (or $\textrm{Div}^{0}(\widetilde{\CbeC}')$) induced by the natural action of $\sigma$ on $\CbeC '=\CbeC\times_{k(x)}k'(x)$ (respectively $\widetilde{\CbeC}'=\widetilde{\CbeC}\times_{k(x)}k'(x)$). Put $\tau :=\sigma\circ\omega$. 
Let $\beta=<u,v>$ be a $\C (x)$-point of $\textrm{Jac}(\widetilde{\CbeC})$ such that $u(0)\neq 0$. Denote by $\check{v}$ the unique polynomial of degree less than or equal to $\textrm{deg}(u)$ such that $\check{v}(0)=0$ and $\check{v}\equiv v\bmod u$.
\begin{enumerate}
\item The point $\beta$ is invariant under $\sigma_{\star\CbeC}$ if and only if one of the two following conditions holds
\begin{enumerate}
\item either $\deg_{s}(u)$ is even, $s^{\deg_{s}(u)}\tau (u)=\sigma (u(0))u(s)$ and the remainder of the Euclidean division of $-\left(\frac{s}{d}\right)^{g+1}\tau (v)$ by $u$ is $v,$ 
\item or the degree of $u$ is $g$ and $\left(\frac{s}{d}\right)^{g+1}\tau (\check{v})=\check{v}$ and $\sigma (u(0))(f-\check{v}^{2})=su(s)s^{g}\tau (u(s)).$
\end{enumerate}
\item If $\beta$ is invariant under $\sigma_{\star\CbeC}$ and the degree of $u$ is strictly less than $g$, then $\varpi (\beta )$ is trivial. 
\item If $\beta$ is invariant under $\sigma_{\star\CbeC}$ and the degree of $u$ is $g$, then $\beta$ is antineutral if and only if 
$u(0)$ is a sum of two squares in $\R (x)$. 

Moreover if $\beta$ is antineutral, then $-d^{g-1}=u(0)h\tau (h)$ where $h$ denotes the function $\frac{t+\check{v}}{su(s)}.$
\end{enumerate}
\end{sprop}
\dem%
Let $D$ be the semi-reduced divisor $\textrm{div}(u,v)$. The point $\beta$ is invariant under $\sigma_{\star\CbeC}$ if and only if the divisors $D$ and $\sigma_{\star\CbeC}(D)=\sigma_{\star\widetilde{\CbeC}}(\widetilde{\omega}(D))$ are linearly equivalent. To study the invariance of $\beta$ under $\sigma_{\star\CbeC}$, we use Cantor's algorithm (see \cite{Cantor}).

Let $e$ be the quotient of the Euclidean division of $\deg_{s}(u)+1$ by $2$ and $\epsilon$ be its remainder. Following Lemma \ref{lemm-image-diviseur-conj-complexe-22-antin} the Mumford's representation of $\sigma_{\star\widetilde{\CbeC}}(\widetilde{\omega}(D))+\textrm{div}(s^{e})$ is $\left(\frac{1}{\sigma (u(0))}s^{2e}\sigma (u)\left(\frac{d^{2}}{s}\right) ,\hat{v}\right)$ where $\hat{v}$ denotes the remainder of the Euclidean division of $-\left(\frac{s}{d}\right)^{g+1}\sigma (v)(\frac{d^{2}}{s})=-\left(\frac{s}{d}\right)^{g+1}\tau (v)$ by $s^{2e}\sigma (u)(\frac{d^{2}}{s})=s^{2e}\tau (u)$ (the number $g$ is odd: see notation \ref{nota-22-caract-antin1}).
\begin{description}
\item[Case $1$: if the degree of $u$ is strictly less than $g$. ] Then 
the divisor \linebreak$\sigma_{\star\widetilde{\CbeC}}(\widetilde{\omega}(D))+\textrm{div}(s^{e})$ has weight less than or equal to $g$. In particular it is reduced. This means that it can be linearly equivalent to $D$ only when it is equal to $D$. In that case the degree of $u$ is even (notice that the degree of $s^{2e}\sigma (u)(\frac{d^{2}}{s})$ is even). 

If $\beta$ is invariant under $\sigma_{\star\CbeC}$ then $D-\sigma_{\star\widetilde{\CbeC}}(\widetilde{\omega}(D)) =\textrm{div}(s^{e})$ and thus $\varpi (\beta )$ is trivial.
\item[Case $2$: if the degree of $u$ is equal to $g$. ]
Denote by $\widetilde{D}$ the reduced divisor $\textrm{div}\left(\frac{f-\check{v}^{2}}{su(s)},w\right)$ where $w$ denotes the remainder of the Euclidean division of $-\check{v}$ by $\frac{f-\check{v}^{2}}{su(s)}$ (notice that its degree is $g$). Applying Cantor's Algorithm we get that $\textrm{div}(su(s),\check{v}) = \widetilde{D}-\textrm{div}\left(\frac{t+\check{v}}{su(s)}\right)$. 

The divisor $\sigma_{\star\widetilde{\CbeC}}(\widetilde{\omega}(D))-D$ is equal to 
$$\begin{array}{rl}
&\textrm{div}\left(\frac{1}{\sigma (u(0))}s^{\deg (u)}\sigma (u)\left(\frac{d^{2}}{s}\right) ,\widetilde{v}\right)  -\textrm{div}(s^{e-1})-div(su(s),\check{v})\\
=&\textrm{div}\left(\frac{1}{\sigma (u(0))}s^{\deg (u)}\sigma (u)\left(\frac{d^{2}}{s}\right) ,\widetilde{v}\right) -\widetilde{D} +\textrm{div}\left(\frac{t+\check{v}}{s^{e}u(s)}\right)\\
\end{array}$$ 
where $\widetilde{v}$ is the remainder of the Euclidean division of $-\left(\frac{s}{d}\right)^{g+1}\sigma (v)(\frac{d^{2}}{s})$ by $s^{\deg (u)}\sigma (u)$. In particular $\sigma_{\star\widetilde{\CbeC}}(\widetilde{\omega}(D))-D$ is principal if and only if the two reduced divisors $\textrm{div}\left(\frac{1}{\sigma (u(0))}s^{\deg (u)}\sigma (u)\left(\frac{d^{2}}{s}\right) ,\widetilde{v}\right)$ and \linebreak$\widetilde{D}=\textrm{div}\left(\frac{f-\check{v}^{2}}{su(s)},w\right)$ are equal. 

When $\sigma (u(0))(f-\check{v}^{2})=su(s)s^{g}\tau (u(s)),$ applying $\tau$ we show that $\widetilde{v}=w$ holds if and only if $\left(\frac{s}{d}\right)^{g+1}\tau (\check{v})=\check{v}$. As a consequence  $\beta$ is invariant under $\sigma_{\star\CbeC}$ if and only if $\sigma (u(0))(f-\check{v}^{2})=su(s)s^{g}\tau (u(s))$ and $\left(\frac{s}{d}\right)^{g+1}\tau (\check{v})=\check{v}.$

If $\beta$ is invariant under $\sigma_{\star\CbeC}$ then $\sigma (u (0))=\frac{s^{g+1}u(s)\tau (u)(s)}{f-\check{v}^{2}}$ belongs to $\R (x)$ (it is $\tau$-invariant hence $\sigma$-invariant) and $\sigma_{\star\CbeC}(D)-D=\textrm{div}\left(\frac{t+\check{v}}{s^{e}u(s)}\right)$. In that case, since $\sigma (u(0))\left(\frac{t+\check{v}}{s^{e}u(s)}\right)\tau\left(\frac{t+\check{v}}{s^{e}u(s)}\right)=-1$, the image $\varpi (\beta )$ is the class of $-u(0)$. 
%
%
%
\fin%
\end{description}

\end{subsection}

\end{section}

\begin{section}{How to find antineutral torsion points.}

\begin{subsection}{The $2$-primary torsion.}

The following Proposition helps us to restrict our study of the existence of an antineutral torsion point to the search of an antineutral $2$-primary torsion point.
\begin{sprop}%
We use the notation of Proposition \ref{prop-finale-HM}.
%
Then \linebreak$\textrm{Pic}^{0}(\CbeD ')^{\Sigma}$ has an antineutral torsion element if and only if $\textrm{Pic}^{0}(\CbeD ')^{\Sigma}$ has an antineutral $2$-primary torsion element.
\end{sprop}
\dem%
Assume the existence of an antineutral torsion element $D$. The order of $D$ is $2^{n}m$ with $m\in\N$ odd and $n\in\N$. 
The morphism $\eta\circ\delta$ takes values in the exponent $2$ group $\R (x)^{\times}/\fbox{2}_{\R (x)}$. Thus the image of a double under $\eta\circ\delta$ is trivial. The integer $m$ being odd, $D$ and $mD$ have the same image under $\eta\circ\delta$. As a consequence, the $2^{n}$-torsion point $mD$ is an antineutral point. 
\fin%


\end{subsection}

\begin{subsection}{A morphism characterizing doubles in the rational points group of some jacobian varieties.}

\begin{snota}\label{notation-Schaefer}
\rm Let $K$ be a characteristic $0$ field and $\overline{K}$ be an algebraic closure of $K.$ 
\rm 
Let $\CbeH$ be a hyperelliptic curve defined over $K$ by an affine equation $\CbeH :y^{2}=f(x)$ where $f(x)$ is a monic separable polynomial of odd degree.

\rm Let $g$ be the genus of $\CbeH .$ The polynomial $f(x)$ has degree $2g+1.$ Let $\alpha_{1},\cdots ,\alpha_{2g+1}$ be the roots of $f$ in $\overline{K}.$ %
\rm Denote by $\infty$ the point at infinity of the curve $\CbeH$ 
and by $P_{i}$ the point $(\alpha_{i},0)$. 
\rm Let $W:=\{ P_{1},\cdots ,P_{2g+1},\infty\}$ be the set of Weierstrass points of $\CbeH $. Denote by 
$\textrm{Div}_{W}^{0}(\CbeH )$ the set $$\{ D\in\textrm{Div}(\CbeH )\vert\textrm{deg}(D)=0\textrm{ and }\textrm{Supp}_{\overline{K}}(D)\cap W=\emptyset\} .$$

\rm Denote by $L$ the algebra $K[T]/(f(T))$. 
The class of a polynomial \linebreak$u\in K[T]$ in $L^{\times }/L^{\times 2}$ is denoted by $[u].$ Denote by $\phi_{\CbeH}$ the morphism
$$\phi_{\CbeH}:\begin{array}[t]{ccc}
\textrm{Div}_{W}^{0}(\CbeH )&\longrightarrow &L^{\times }/L^{\times 2}\\
\displaystyle\sum_{i\in I}n_{i}(Q_{i}-\deg (Q_{i})\infty)&\longmapsto &\left[\displaystyle\prod_{i\in I}(x(Q_{i})-T)^{n_{i}}\right]
\end{array}.$$
In terms of Mumford's representation, the map $\phi_{\CbeH}$ sends a semi-reduced divisor $\textrm{div}(u,v)$ (with $u$ monic as usual) to the class $\left[(-1)^{\deg (u)}u(T)\right]$.
\end{snota}

\begin{sprop}\label{def-Cas-Sch+im-in-ker}
We use the notation \ref{notation-Schaefer}. Then $\phi_{\CbeH}$ induces a morphism $\pi_{\CbeH}:\textrm{Jac}(\CbeH )(K)\longrightarrow L^{\times }/L^{\times 2}$ with kernel $2\textrm{Jac}(\CbeH )(K)$ and with image contained in the kernel of the norm map $N_{L/K}:L^{\times}/L^{\times 2}\longrightarrow K^{\times}/K^{\times 2}$.
\end{sprop}
\dem%
See \cite{Schaefer1995} Theorems 1.1 and 1.2. 
\fin\\

\noindent To use Proposition \ref{def-Cas-Sch+im-in-ker} in next sections we recall the following characterization of the squares in a quadratic extension.

\begin{sprop}\label{caract_carres_ext_quad}
Let $k_{0}$ be a field of characteristic different from $2$. Let $\delta$ be an element of $k_{0}$ which is not a square in $k_{0}$. Denote by $k$ the quadratic extension $k:=k_{0}[U]/(U^{2}-\delta ).$%
\begin{enumerate}
\item 
Let $\alpha$ and $\beta$  be two elements of $k_{0}$. We assume that $\alpha$ is nontrivial. Then $\alpha U+\beta$ is a square in $k$ if and only if there exists $\gamma ,$ $\eta\in k_{0}$ such that $N_{k/k_{0}}(\alpha U+\beta )=\gamma^{2}$ and $\frac{\beta +\gamma}{2}=\eta^{2}.$
\item Let $\beta$ be an element of $k_{0}.$ Then $\beta$ is a square in $k$ if and only if $\beta$ or $\delta\beta$ is a square in $k_{0}.$ 
\end{enumerate}
\end{sprop}

\end{subsection}



\begin{subsection}{A family of jacobian varieties without antineutral torsion point.}

\begin{snota}\label{notations-sec-25-intro}
\rm Let $B$ and $C$ be two elements of $\R (x)$. We consider the polynomial $P(x,y^{2}):=(y^2+1)(y^{2}+C)(y^{4}+(1+C)y^{2}+B)$ which is assumed to be squarefree. Let $\CbeC$ be the hyperelliptic curve defined over $\R (x)$ by the affine equation $\CbeC :z^{2}+P(x,y^{2})=0.$ 
\end{snota}

\begin{snota}\label{notations-sec-25}
\rm We use the notation of Proposition \ref{theo-caracterisation-antineutre-final} relative to the curve $\CbeC$. In particular we introduce $d:=(1-C)(B-C)$ and the three polynomials 
$g_{1}(s):=\frac{-(s+d)^{2}+(s-d)^{2}}{-4d}=s$, 
$g_{2}(s):=\frac{-(s+d)^{2}+C(s-d)^{2}}{C-1}$ and 
\linebreak$g_{3}(s):=\frac{(s+d)^{4}-(1+C)(s+d)^{2}(s-d)^{2}+B(s-d)^{4}}{B-C}$. 

The curve $\CbeC ':=\CbeC\times_{\R (x)}\C (x)$ is birationally equivalent to the curve $\widetilde{\CbeC}':=\widetilde{\CbeC}\times_{\R (x)}\C (x)$ with $\widetilde{\CbeC}$ the hyperelliptic curve defined over $\R (x)$ by the affine equation $\widetilde{\CbeC}:t^{2}=g_{1}(s)g_{2}(s)g_{3}(s).$

We denote by $\sigma$ be the complex conjugation and by $\tau$ the involution
$$\begin{array}[t]{cccc}
\tau :&\C (x)(\widetilde{\CbeC}) & \longrightarrow &  \C (x)(\widetilde{\CbeC})\\
&a(s)t+b(s)& \longmapsto & -\sigma (a)\left(\frac{d^{2}}{s}\right)\frac{d^{4}t}{s^{4}}+\sigma (b)\left(\frac{d^{2}}{s}\right) .\\
\end{array}$$
The involution $\tau$ is the same as the involution $\tau$ used in Proposition \ref{theo-caracterisation-antineutre-final}.

\rm For every index $i=1,2,3,$ let $k_{i}$ be the algebra $\C(x)[T]/(g_{i}(T)).$ 
Denote by $\pi_{\widetilde{\CbeC}}:\textrm{Jac}(\widetilde{\CbeC})(\C (x))\longrightarrow k_{1}^{\times}/k_{1}^{\times 2}\times k_{2}/k_{2}^{\times 2}\times k_{3}/k_{3}^{\times 2}$ the morphism defined by applying Proposition \ref{def-Cas-Sch+im-in-ker} to the curve $\widetilde{\CbeC}$ and denote by $\pi_{\widetilde{\CbeC},i}:\textrm{Jac}(\widetilde{\CbeC})(\C (x))\longrightarrow\nolinebreak k_{i}/k_{i}^{\times 2}$ its i-th coordinate. If $\alpha$ corresponds to the linear equivalence class of a semi-reduced divisor $\textrm{div}(u,v)$ with $u$ coprime to $g_{i},$ then $\pi_{\widetilde{\CbeC},i}(\alpha )$ is the class in $k_{i}/k_{i}^{\times 2}$ of $(-1)^{\deg (u)}u(T)$.
\end{snota}
\begin{sprop}\label{prop-25-2-torsion} We use the notation \ref{notations-sec-25-intro} and \ref{notations-sec-25}. We assume $B,$ $C$ and $(1+C)^{2}-4B$ are not squares in $\C (x).$

Then the $2$-torsion elements of $\textrm{Jac}\left(\widetilde{\CbeC}\right)\left(\C\left( x\right)\right)$ are the neutral element and the three points 
$<g_{1},0>,$ $<g_{2},0>$ and $<g_{1}g_{2},0>.$
\end{sprop}
\dem%
The hypotheses imply that the three polynomials $g_{1}$, $g_{2}$ and $g_{3}$ are irreducible (for a detailed proof see \cite{These-V-Mahe}). This is sufficient since the $2$-torsion points are the points $<u,0>$ with $u$ a divisor of degree less than the genus $g$ of $f$ (see \cite{Mumford}). \fin%
\begin{sprop}\label{T_double?_part_1} We use the notation \ref{notations-sec-25-intro} and \ref{notations-sec-25}. We assume $B,$ $C$ and $(1+C)^{2}-4B$ are not squares in $\C (x)$. 

If $<g_{1},0>$ is a double in $\textrm{Jac}\left(\widetilde{\CbeC}\right)\left(\C\left( x\right)\right) ,$ then 
either $(B-C)\in\C (x)^{\times 2}$ or $C(B-C)\in\C (x)^{\times 2}$.

\end{sprop}
\dem%
The image of an element of $2\textrm{Jac}\left(\widetilde{\CbeC}\right)\left(\C\left( x\right)\right)$ by $\pi_{\widetilde{\CbeC},2}$ is trivial. In particular, if $<g_{1},0>$ is a double in $\textrm{Jac}\left(\widetilde{\CbeC}\right)\left(\C\left( x\right)\right)$ then the class of $-T$ in $k_{2}$ 
is a square. 

Denote by $L_{2}$ the field $\C (x)[U]/(U^{2}-4C(B-C)^{2})$ (this is a field since $C$ is not a square in $\C (x)$). 
The morphism %
$\varphi_{2}:\begin{array}[t]{rcl}
L_{2}&\longrightarrow &k_{2}\\
U&\longmapsto &T+(1+C)(B-C).\\
\end{array}$ %
is an isomorphism. As a consequence, if the class of $-T$ in $k_{2}$ is a square, then the class of $-U+(1+C)(B-C)$ in $L_{2}$ is a square.

The field $L_{2}$ is a quadratic extension of $\C (x)$. 
We apply Proposition \ref{caract_carres_ext_quad}:
the class of $-U+(1+C)(B-C)$ in $L_{2}$ is a square if and only if either $(B-C)$ or $C(B-C)$ is a square in $\C (x)$. %
\fin%

\begin{sprop}\label{section-25-8-torsion} We use the notation \ref{notations-sec-25-intro} and \ref{notations-sec-25}. Denote by $\sigma_{\star\CbeC}$ the action of $\sigma$ on $\textrm{Jac}(\widetilde{\CbeC}')$ induced by the natural action of $\sigma$ on \linebreak$\CbeC '=\CbeC\times_{\R (x)}\C (x)$.

Then the point $<s-d,8d^{3}>$ is a $8$-torsion element of $\textrm{Jac}\left(\widetilde{\CbeC}\right)\left(\C\left( x\right)\right)$. Its double is not $\sigma_{\star\CbeC}$-invariant but $[4]<s-d,8d^{3}>$ is equal to $<g_{1}g_{2},0>$.
\end{sprop}
\dem%
The equality $[4]<s-d,8d^{3}>=<g_{1}g_{2},0>$ can be checked by using Cantor's algorithm for addition in $\textrm{Jac}\left(\widetilde{\CbeC}\right)\left(\C\left( x\right)\right)$. 
From this equality we deduce that $<s-d,8d^{3}>$ is a $8$-torsion point.\\
The remainder of the Euclidean division of $-\frac{s^{4}}{d^{4}}\tau (16d^{2}s-8d^{3})=8\frac{s^{4}}{d}-16s^{3}$ by $(s-d)^{2}$ is $-\left( 16d^{2}s-8d^{3}\right)$. It is different from $16d^{2}s-8d^{3}$. 
Following Proposition \ref{theo-caracterisation-antineutre-final}, the point $[2]<s-d,8d^{3}>=<(s-d)^{2},16d^{2}s-8d^{3}>$ is not $\sigma_{\star\CbeC}$-invariant.
\fin%
\begin{sprop}\label{coro-torsion-finie} Let $B$ and $C$ be two elements of $\R (x).$ Let $\CbeC$ be the hyperelliptic curve defined on $\R (x)$ by the affine equation
$$z^{2}+(y^2+1)(y^{2}+C)(y^{4}+(1+C)y^{2}+B)=0.$$ 
We assume that $B$, $C$, $C(B-C)$, $B-C,$ $(B-C)(1-C)$ and $(1+C)^{2}-4B$ are not squares in $\C (x)$. 
Then the $2$-primary torsion subgroup of $\textrm{Jac}(\CbeC )(\C (x))$ is finite.
\end{sprop}
\dem%
Since $C$, $1-C$, $B-C$, $B$ ans $(1+C)^{2}-4B$ are different from $0$, the polynomial $P(x,y^{2}):=(y^2+1)(y^{2}+C)(y^{4}+(1+C)y^{2}+B)$ is squarefree. We use the notation \ref{notations-sec-25}.

Following Propositions \ref{prop-25-2-torsion} and \ref{section-25-8-torsion}, the group $\textrm{Jac}\left(\widetilde{\CbeC}\right)\left(\C\left( x\right)\right)[2]$ is generated by $<g_{1},0>$ and $<g_{1}g_{2},0>=[4]<s-d,8d^{3}>$. If $T$ is a $4$-torsion element of $\textrm{Jac}\left(\widetilde{\CbeC}\right)\left(\C\left( x\right)\right)$, then 
$$2T=n_{1}<g_{1},0>+n_{2}<g_{1}g_{2},0>=n_{1}<g_{1},0>+4n_{2}<s-d,8d^{3}>$$ 
with $n_{1}$, $n_{2}\in\{ 0,1\}$. If $n_{1}=1$, then $<g_{1},0>$ is a double in $\textrm{Jac}\left(\widetilde{\CbeC}\right)\left(\C\left( x\right)\right)$. Since this is not the case (see Proposition \ref{T_double?_part_1}), the integer $n_{1}$ is equal to $0$ and thus $2T$ is a multiple of $[4]<s-d,8d^{3}>$. In particular $T$ is the sum of a $2$-torsion point and a multiple of $[2]<s-d,8d^{3}>$. In the same way, we can prove that $\textrm{Jac}\left(\widetilde{\CbeC}\right)\left(\C\left( x\right)\right)[2]$ is generated by $<g_{1},0>$ and $<s-d,8d^{3}>$.

Let $T$ be a $16$-torsion element of $\textrm{Jac}\left(\widetilde{\CbeC}\right)\left(\C\left( x\right)\right)$. The point $2T$ being an $8$-torsion point, there exists two integers $n_{1},$ $n_{2}\in\N$ such that \linebreak$2T=n_{1}<g_{1},0>+n_{2}<s-d,8d^{3}>$. Applying $\pi_{\widetilde{\CbeC},1}$ to this equality gives $1\sim d^{6n_{1}}d^{n_{2}}$. By hypothesis this equivalence is possible only if $n_{2}$ is even i.e. only if $2T$ is a $4$-torsion point. In that case $T$ is an $8$-torsion point.

The $16$-torsion subgroup of $\textrm{Jac}\left(\widetilde{\CbeC}\right)\left(\C\left( x\right)\right)$ being equal to the $8$-torsion subgroup, the $2$-primary torsion subgroup is equal to the $8$-torsion subgroup. Thus the $2$-primary torsion subgroup is finite.
\fin%
\begin{stheo}\label{existence-pt-anti-exemple} Let $B$ and $C$ be two elements of $\R (x).$ Let $\CbeC$ be the hyperelliptic curve defined on $\R (x)$ by the affine equation
$$z^{2}+(y^2+1)(y^{2}+C)(y^{4}+(1+C)y^{2}+B)=0.$$ 
We assume that $B,$ $C,$ $C(B-C),$ $B-C$, $(1-C)(B-C)$ and $(1+C)^{2}-4B$ are not squares in $\C (x).$ Then $\textrm{Jac}(\CbeC )(\R (x))$ has no antineutral torsion point.
\end{stheo}
\dem%
Denote by $\sigma_{\star\CbeC}$ the action of $\sigma$ on $\textrm{Jac}(\widetilde{\CbeC}')$ induced by the natural action of $\sigma$ on $\CbeC '=\CbeC\times_{\R (x)}\C (x)$. We use the notation \ref{notations-sec-25}. %
As in the proof of Proposition \ref{coro-torsion-finie}, the polynomial $(y^2+1)(y^{2}+C)(y^{4}+(1+C)y^{2}+B)$ is squarefree 
%
and the $4$-torsion subgroup of $\textrm{Jac}\left(\widetilde{\CbeC}\right)\left(\C\left( x\right)\right)$ 
is generated by $<g_{1},0>$ and $<(s-d)^{2},16d^{2}s-8d^{3}>.$ 

From Proposition \ref{section-25-8-torsion} the point $<(s-d)^{2},16d^{2}s-8d^{3}>$ is not \linebreak$\sigma_{\star\CbeC}$-invariant, but using Proposition \ref{theo-caracterisation-antineutre-final} we can check the $\sigma_{\star\CbeC}$-invariance of the points $<g_{1}(s),0>$ and $<g_{2},0>=[2]<(s-d)^{2},16d^{2}s-8d^{3}>$. 
If a \linebreak$4$-torsion point $n_{1}<g_{1},0>+n_{2}<(s-d)^{2},16d^{2}s-8d^{3}>$ is \linebreak$\sigma_{\star\CbeC}$-invariant, then $n_{2}<(s-d)^{2},16d^{2}s-8d^{3}>$ is $\sigma_{\star\CbeC}$-invariant i.e. $n_{2}$ is even. As a consequence, the $4$-torsion subgroup of $\textrm{Jac}\left(\widetilde{\CbeC}\right)\left(\C\left( x\right)\right)^{\sigma_{\star\CbeC}}$ (the group of $\sigma_{\star\CbeC}$-invariant elements of $\textrm{Jac}\left(\widetilde{\CbeC}\right)\left(\C\left( x\right)\right)$) is equal to the $2$-torsion subgroup and thus to the $2$-primary torsion subgroup of $\textrm{Jac}\left(\widetilde{\CbeC}\right)\left(\C\left( x\right)\right)^{\sigma_{\star\CbeC}}$. 
In particular the $2$-primary torsion of $\textrm{Jac}\left(\widetilde{\CbeC}\right)\left(\C\left( x\right)\right)^{\sigma_{\star\CbeC}}$ is generated by $<g_{1},0>$ and $<g_{2},0>$. 

The degrees of $g_{1}$ and $g_{2}$ are strictly less than $3.$ Following Proposition \ref{theo-caracterisation-antineutre-final} (more precisely the assertion $2.$), 
the jacobian variety $\textrm{Jac}\left(\widetilde{\CbeC}\right)\left(\C\left( x\right)\right)$ has no antineutral ($2$-primary) torsion point.
%
\fin%

\end{subsection}

\begin{subsection}{Proposition \ref{formule-4-tors-non-triviale}: the heuristic.}\label{ref-heuristic-prod-sum3squares}

We use the notation of Proposition \ref{formule-4-tors-non-triviale} and  Proposition \ref{theo-caracterisation-antineutre-final}:
\begin{itemize}
\item[* ] we denote by $\CbeC$ the hyperelliptic curve defined over $\R (x)$ by the affine equation $z^{2}+P(x,y)= 0;$ 
\item[* ] we note $d=(1-a)(1-b)(1-c)$; 
\item[* ] For each $\alpha\in\C (x)$  we consider the polynomial
$$g_{\alpha}(s):=s^{2}+2d\frac{1+\alpha}{1-\alpha}s+d^{2};$$
\item[* ] we denote by $\widetilde{\CbeC}$ the hyperelliptic curve defined over $\R (x)$ by the affine equation $t^{2}=sg_{a}(s)g_{b}(s)g_{c}(s).$
\end{itemize}
The $2$-torsion elements of $\textrm{Jac}(\CbeC )(\C (x))=\textrm{Jac}(\widetilde{\CbeC})(\C (x))$ are the point with Mumford representation $<u,0>$ where $u\in\C (x)[y]$ is a monic divisor of $sg_{a}(s)g_{b}(s)g_{c}(s).$ 

Using Proposition \ref{def-Cas-Sch+im-in-ker} and Proposition \ref{caract_carres_ext_quad} 
we show the point \linebreak$<g_{a}(s),0>$ is a double if and only if the three following assumptions are satisfied:
\begin{enumerate}
\item either $(b-a)(c-a)\in\C (x)^{\times 2}$ or $a(b-a)(c-a)\in\C (x)^{\times 2},$
\item either $(b-a)(1-a)\in\C (x)^{\times 2}$ or $b(b-a)(1-a)\in\C (x)^{\times 2},$
\item either $(c-a)(1-a)\in\C (x)^{\times 2}$ or $c(c-a)(1-a)\in\C (x)^{\times 2}.$
\end{enumerate}
Assume the existence of $\widetilde{\beta},\widetilde{\gamma}\in\C (x)$ such that $(b-a)=\widetilde{\beta}^{2}(1-a)$ and $(c-a)=\widetilde{\gamma}^{2}(1-a)$. Then $<g_{a}(s),0>$ is the double of $<u,v>$ with 
\begin{itemize}
\item[* ] $u:=\frac{(\widetilde{\beta}\widetilde{\gamma}(s-d)^{2} +(1+\widetilde{\beta}+\widetilde{\gamma})g_{a})(s-d)}{(1+\widetilde{\beta})(1+\widetilde{\gamma})},$
\item[* ] $q:=\frac{(1-a)^{2}(g_{a}+(\widetilde{\beta}\widetilde{\gamma}+\widetilde{\beta}+\widetilde{\gamma})(s-d)^{2})}{2d}$
\item[* ] $v$ the remainder of the euclidean division of $qg_{a}$ by $u.$
\end{itemize}
From Proposition \ref{theo-caracterisation-antineutre-final} we get conditions on $a,b,c$  under which $<u,v>$ is antineutral. Those conditions are satisfied under the assumptions of Proposition \ref{formule-4-tors-non-triviale}.

Moreover, under the hypotheses of Proposition \ref{formule-4-tors-non-triviale}, the evaluation $u(0) =-d^{3}$ is a square in $\R (x)$ and Proposition \ref{theo-caracterisation-antineutre-final} gives a way to write $-1$ as a sum of two squares in $\R (x)(\CbeC )$ i.e. a way to write $P(x,y)$ as a sum of three squares  in $\R (x,y)$ (in fact $P(x,y)$ is a sum of two squares in $\R (x,y)$ if and only if $-1$ is a sum of two squares in $\R(x)(\CbeC )$; for an explicit version of this equivalence see \cite{Lam}).

\end{subsection}
\end{section}

\begin{section}{Symplifying some Mordell-Weil's ranks computations.}

\begin{subsection}{An application of Lang-N\'eron theorem.}


\begin{stheo}[Lang, N\'eron]\label{Lang_Neron}
Let $k$ be a field. Let $F$ be the function field of a variety defined over $k.$ Let $A$ be an abelian variety defined over $F.$ We assume that no abelian subvariety $B$ of $A$ 
can be obtained by scalar extension from an abelian variety defined over $k$ and of dimension at least $1$. 
Then the abelian group $A(F)$ is finitely generated.
\end{stheo}
\dem %
See \cite{Lang} Page $27$ Theorem $4.2 .$
\fin%

\begin{scoro}\label{finitude-jac-CbeC-final} Let $B$ and $C$ be two elements of $\R (x).$ Let $\CbeC$ be the hyperelliptic curve defined over $\R (x)$ by the affine equation 
$$z^{2}+(y^2+1)(y^{2}+C(x^{2}))(y^{4}+(1+C(x^{2}))y^{2}+B(x^{2}))=0.$$ 
Assume the polynomials $B(x^{2})$, $C(x^{2})$, $B(x^{2})-C(x^{2}),$ $C(x^{2})(B(x^{2})-C(x^{2}))$, $(B(x^{2})-C(x^{2}))(1-C(x^{2}))$ and $(1+C(x^{2}))^{2}-4B(x^{2})$ are not squares in $\C (x)$. %
Then the abelian group $\textrm{Jac}(\CbeC )(\C (x))$ is finitely generated.
\end{scoro}
\dem%
No abelian subvariety $A$ of $\textrm{Jac}(\CbeC )$ of dimension at least $1$ can be defined by scalar extension from an abelian variety defined over $\C$. In fact, if such a variety did exist, then the order of the $2$-primary torsion subgroup of $A(\C )$ would be infinite (since $\C$ is algebraically closed) and we would have a contradiction with Corollary \ref{coro-torsion-finie}. 
Thus the hypotheses of Theorem \ref{Lang_Neron} are satisfied. 
\fin%

\end{subsection}

\begin{subsection}{Involutions and Mordell-Weil's ranks.}

\begin{slemm}\label{4-tors-Pic-F-is-finite}
Let $K$ be a field. Let $F$ be a function field over $K$ (i.e. a finite extension of $K(\alpha )$ for some transcendent element $\alpha\in F$). Then the order of the $4$-torsion of $\textrm{Pic}^{0}(F)$ is finite.
\end{slemm}
\dem%
By definition the field $F$ is the function field of 
some smooth projective curve $\CbeD$ defined over $K$. There is an injection $\textrm{Pic}^{0}(F)$ to $\textrm{Jac}(\CbeD )(K)$. Lemma \ref{4-tors-Pic-F-is-finite} follows from the finiteness of the order of $\textrm{Jac}(\CbeD )(K)[4]$.
\fin
\begin{snota}
We use the definition and notation of \cite{Stichtenoth}. Let $F$ be a function field with full constant field $k$ and let $F_{2}$ be a finite extension of $F$ with full constant field $k_{2}$. 
If $\mathfrak{p}$ is a place of $F/k$, let $Cn_{F_{2}/F}(\mathfrak{p})$ be the divisor $\displaystyle\sum_{\mathcal{P}|\mathfrak{p}}e(\mathcal{P}|\mathfrak{p})\mathcal{P}\in\textrm{Div}(F_{2}/k_{2})$. %
By linearity we get a homomorphism \linebreak$Cn_{F_{2}/F}:\textrm{Div}^{0}(F/k)\longrightarrow\textrm{Div}^{0}(F_{2}/k_{2})$. %
The homomorphism $Cn_{F_{2}/F}$ induces a group homomorphism $CN_{F_{2}/F}$ from the quotient $\textrm{Pic}^{0}(F/k)$ to $\textrm{Pic}^{0}(F_{2}/k_{2})$.

If $\rho$ is an automorphism of a field $F$ we denote by $F^{\rho}$ 
the subfield of $\rho$-invariant 
elements of $F$.
\end{snota}

\begin{sprop}\label{decomp-pic-theo-gene} Let $k$ be a characteristic $0$ field. Let $P(T)\in k[T]$ be a polynomial and $\CbeH$ be the hyperelliptic curve defined over $k$ by the affine equation $z^{2}+P(y)=0.$ Let $\iota :
k(\CbeH )\longrightarrow k(\CbeH )
$ be the hyperelliptic involution and $\rho:k(\CbeH )\longrightarrow k(\CbeH )$ be an involution distinct from the identity map and from $\iota .$ Assume the commutativity of $\iota$ and $\rho$. 

Then the homomorphism $\varphi:=+\circ \left(CN_{k(\CbeH )/k(\CbeH )^{\iota\circ\rho}}\times CN_{k(\CbeH )/k(\CbeH )^{\rho}}\right)$ 
has a finite kernel and its image contains $2\textrm{Pic}^{0}(k(\CbeH )/k).$
\end{sprop}
\dem%
\begin{description}
\item[Step $1$. ] Let $\mathfrak{p}$ be a place of $k(\CbeH )^{\rho}/k^{\rho}$. The extension $k(\CbeH )/k(\CbeH )^{\rho}$ is a degree~$2$ Galois extension with Galois group $\{\textrm{Id},\rho\}$. Thus:
\begin{enumerate}
\item for every place $\mathcal{P}$ of $k(\CbeH )/k$ above $\mathfrak{p}$ the ramification indexes $e(\mathcal{P}|\mathfrak{p})$ and $e(\rho (\mathcal{P})|\mathfrak{p})$ are equal (see \cite{Stichtenoth} Corollary III.7.2);
\item $\rho$ induces a bijection from the set of places of $k(\CbeH )/k$ above $\mathfrak{p}$ into itself.
\end{enumerate}
In particular $Cn_{k(\CbeH )/k(\CbeH )^{\rho}}(\mathfrak{p})$ is $\rho$-invariant. Since this is true for every place $\mathfrak{p}$, we deduce from the linearity of $Cn_{k(\CbeH )/k(\CbeH )^{\rho}}$ that every element of the image of $Cn_{k(\CbeH )/k(\CbeH )^{\rho}}$ is $\rho$-invariant.

\item[Step $2$. ] if $F/k$ is a function field and $f$ is an element of $F$, we denote by $\textrm{div}_{F/k}(f)$ the principal divisor associated to $f$.

Let $D$ be a divisor of $k(\CbeH )^{\rho}/k^{\rho}$ whose linear equivalence class lies in $\textrm{Ker}(CN_{k(\CbeH )/k(\CbeH )^{\rho}})$. The divisor $Cn_{k(\CbeH )/k(\CbeH )^{\rho}}(D)$ is the principal divisor $\textrm{div}_{k(\CbeH )/k}(f)$ associated to some function $f\in k(\CbeH )$. 
Following Step~$1$ the divisor $Cn_{k(\CbeH )/k(\CbeH )^{\rho}}(D)+\rho (Cn_{k(\CbeH )/k(\CbeH )^{\rho}}(D))=\textrm{div}_{k(\CbeH )/k}(f\rho (f))$ is equal to $Cn_{k(\CbeH )/k(\CbeH )^{\rho}}(2D)$. The function $f\rho (f)$ belongs to $k(\CbeH)^{\rho}$; the divisor $\textrm{div}_{k(\CbeH )/k}(f\rho (f))$ is the image under $Cn_{k(\CbeH )/k(\CbeH )^{\rho}}$ 
of the principal divisor $\textrm{div}_{k(\CbeH)^{\rho}/k^{\rho}}(f\rho (f))$. 
The map \linebreak$Cn_{k(\CbeH )/k(\CbeH )^{\rho}}$ being injective, the divisors $2D$ and $\textrm{div}_{k(\CbeH)^{\rho}/k^{\rho}}(f\rho (f))$ are equal. 

In the same way we prove that the kernel of $CN_{k(\CbeH )/k(\CbeH )^{\iota\circ\rho}}$ is included in the $2$-torsion of $\textrm{Pic}^{0}(k(\CbeH )^{\iota\circ\rho}/k^{\iota\circ\rho})$.

\item[Step $3$. ] Let $(\alpha_{\rho},\alpha_{\iota\circ\rho})$ be an element of $\textrm{Ker}(\varphi )$. Then the divisor \linebreak$CN_{k(\CbeH )/k(\CbeH )^{\iota\circ\rho}}(\alpha_{\iota\circ\rho})=-CN_{k(\CbeH )/k(\CbeH )^{\rho}}(\alpha_{\rho})$ is both $\rho$-invariant and \linebreak$\iota\circ\rho$-invariant (see Step~$1$). In particular $CN_{k(\CbeH )/k(\CbeH )^{\iota\circ\rho}}(\alpha_{\iota\circ\rho})$ is \linebreak$\iota$-invariant i.e. it is a $2$-torsion element of $\textrm{Pic}^{0}(k(\CbeH )/k)$. Hence the orders of $\alpha_{\rho}$ and $\alpha_{\iota\circ\rho}$ are at most $4$ (see Step~$2$). We apply Lemma \ref{4-tors-Pic-F-is-finite}: there are only finitely many choices for $\alpha_{\rho}$ and $\alpha_{\iota\circ\rho}$. 

\item[Step $4$. ] Let $\mathcal{P}$ be a place of $k(\CbeH )/k$. Denote the place $\mathcal{P}\cap k(\CbeH )^{\iota\circ\rho}$ by $\mathfrak{p}$. The places of $k(\CbeH )/k$ above $\mathfrak{p}$ are $\mathcal{P}$ and $\iota\circ\rho (\mathcal{P})$. 
\begin{enumerate}
\item If $\mathcal{P}$ and $\iota\circ\rho (\mathcal{P})$ are distinct then $\mathcal{P}$ and $\iota\circ\rho (\mathcal{P})$ are not ramified; in that case $\mathcal{P}+\iota\circ\rho (\mathcal{P})$ is equal to $Cn_{k(\CbeH )/k(\CbeH )^{\iota\circ\rho}}(\mathcal{P})$;
\item  If $\mathcal{P}=\iota\circ\rho (\mathcal{P})$ then $\mathcal{P}+\iota\circ\rho (\mathcal{P})$ is equal to 
$$2\mathcal{P}=e(\mathcal{P}|\mathfrak{p})f(\mathcal{P}|\mathfrak{p})\mathcal{P}=Cn_{k(\CbeH )/k(\CbeH )^{\iota\circ\rho}}(f(\mathcal{P}|\mathfrak{p})\mathcal{P})$$ 
(where $f(\mathcal{P}|\mathfrak{p})$ denotes the residual degree of $\mathcal{P}$ above $\mathfrak{p}$).
\end{enumerate}
In both cases $\mathcal{P}+\iota\circ\rho (\mathcal{P})$ is in the image of $Cn_{k(\CbeH )/k(\CbeH )^{\iota\circ\rho}}$. In the same way we show that $\mathcal{P}+\rho (\mathcal{P})$ is in the image of $Cn_{k(\CbeH )/k(\CbeH )^{\rho}}$. By linearity this implies that the image of $\varphi$ contains the linear equivalence class of each divisor $D+\rho (D)+D+\iota\circ\rho (D)$ (where $D$ denotes a degree~$0$ divisor of $k(\CbeH )/k$). Since the linear equivalence class of $D+\iota (D)$ is trivial for every degree~$0$ divisor $D$ of $k(\CbeH )/k$, the image of $\varphi$ contains the double of each element of $\textrm{Pic}^{0}(k(\CbeH )/k)$. 
\fin%
\end{description}

\begin{slemm}\label{lemm_decomp-jac_pasge_C-R}
Let $\CbeD$ be a smooth projective geometrically integral curve defined over $\R (x)$. Assume that $\CbeD$ has a $\C (x)$-point. Then the following inclusions hold
$$2\textrm{Jac}(\CbeD )(\R (x))\subset\textrm{Pic}^{0}(\R (x)(\CbeD ))\subset\textrm{Jac}(\CbeD )(\R (x)) .$$
\end{slemm}
\dem%
Let $\CbeD ':=\CbeD\times_{\R (x)}\C (x)$ be the complexified of $\CbeD$. Denote by $\Sigma$ the Galois group $\textrm{Gal}(\C (x)/\R (x))=\textrm{Gal}(\C /\R )$. Following Lemma \ref{suite_exacte_antineutralite} we have an exact sequence
$$\xymatrix{ 0 \ar[r] & \textrm{Pic}(\R (x)(\CbeD )) \ar[r]^-{p^{*}} & \textrm{Pic}(\C (x)(\CbeD '))^{\Sigma} \ar[r]^-{\delta} & H^{1}(\Sigma , \C (x)(\CbeD ')^{\times}/\C (x)^{\times}) \ar[r] & 0. \\}$$
Using $p^{*}$ we identify $\textrm{Pic}(\R (x)(\CbeD ))$ with a subgroup of $\textrm{Pic}(\C (x)(\CbeD '))^{\Sigma}$. The exponent of $H^{1}(\Sigma , \C (x)(\CbeD ')^{\times}/\C (x)^{\times})$ is $2$. Thus $\textrm{ker}(\delta )=\textrm{Pic}(\R (x)(\CbeD ))$ contains $2\textrm{Pic}(\C (x)(\CbeD '))^{\Sigma}$. To conclude we notice that $\textrm{Jac}(\CbeD )(\R (x))=\textrm{Pic}^{0}(\C (x)(\CbeD '))^{\Sigma}$.
\fin%

\begin{sprop}\label{prop-decomp-jac}
Let $P(T)\in\R (x)[T]$ be a polynomial and $\CbeC$ be the hyperelliptic curve defined over $\R (x)$ by the affine equation $z^{2}+P(y^{2})=0.$ Assume that $\CbeC$ has a $\C (x)$-rational point and that $\textrm{Jac}(\CbeC)(\R (x))$ is finitely generated. 
Consider the two following $\R (x)$-hyperelliptic curves 
$$
\CbeC^{+}:t^{2}+sP(s)=0\textrm{ and }
\CbeC^{-}:\beta^{2}+P(\alpha )=0.
$$
Then the Mordell-Weil rank of $\textrm{Jac}(\CbeC )(\R (x))$ is the sum of Mordell-Weil ranks of $\textrm{Jac}(\CbeC^{+})(\R (x))$ and $\textrm{Jac}(\CbeC^{-})(\R (x))$.
%
\end{sprop}
\dem%
Consider the two involutions $$\iota:\begin{array}[t]{ccc}
\R (x)(\CbeC )&\longrightarrow &\R (x)(\CbeC )\\
A(y,z)&\longmapsto &A(y,-z)\\
\end{array}
\textrm{ and }\rho :\begin{array}[t]{ccc}
\R (x)(\CbeC )&\longrightarrow &\R (x)(\CbeC )\\
A(y,z)&\longmapsto &A(-y,z)\\
\end{array}.$$ 
The map $\iota$ is the hyperelliptic involution. Consider the two morphisms $$\phi^{+}:\begin{array}[t]{ccc}
\R (x)(\CbeC^{+})&\longrightarrow &\R (x)(\CbeC )\\
A(s,t)&\longmapsto &A(y^{2},yz)\\
\end{array}\textrm{ and }%
\phi^{-}:\begin{array}[t]{ccc}
\R (x)(\CbeC^{-})&\longrightarrow &\R (x)(\CbeC )\\
A(s,t)&\longmapsto &A(y^{2},z)\\
\end{array}.$$
Since $\iota\circ\rho\circ\phi^{+}$ and $\phi^{+}$ are equal, the image of $\phi^{+}$ is a subfield of $\R (x)(\CbeC )^{\iota\circ\rho}.$ Moreover $\R (x)(\CbeC )$ is a degree~$2$ extension of $\textrm{Im}(\phi^{+})$. Thus the image of $\phi^{+}$ is $\R (x)(\CbeC )^{\iota\circ\rho}.$ In the same way, we prove that the image of $\phi^{-}$ is $\R (x)(\CbeC )^{\rho}.$ Proposition \ref{decomp-pic-theo-gene} asserts the existence of a group homomorphism 
$$\varphi :\textrm{Pic}^{0}(\R (x)(\CbeC^{+}))\times\textrm{Pic}^{0}(\R (x)(\CbeC^{-}))\longrightarrow\textrm{Pic}^{0}(\R (x)(\CbeC ))$$
with finite kernel and whose image contains $2\textrm{Pic}^{0}(\R (x)(\CbeC )).$ From Lemma \ref{lemm_decomp-jac_pasge_C-R} we deduce the following inclusions:
\begin{itemize}
\item[* ] $2\textrm{Jac}(\CbeC )(\R (x))\subset\textrm{Pic}^{0}(\R (x)(\CbeC ))\subset\textrm{Jac}(\CbeC )(\R (x))$,
\item[* ] $2\textrm{Jac}(\CbeC^{+} )(\R (x))\subset\textrm{Pic}^{0}(\R (x)(\CbeC^{+} ))\subset\textrm{Jac}(\CbeC^{+} )(\R (x))$ and
\item[* ] $2\textrm{Jac}(\CbeC^{-} )(\R (x))\subset\textrm{Pic}^{0}(\R (x)(\CbeC^{-} ))\subset\textrm{Jac}(\CbeC^{-} )(\R (x))$.
\end{itemize}
In particular the map $\varphi$ induces a group homomorphism from \linebreak$2\textrm{Jac}(\CbeC^{+})(\R (x))\times 2\textrm{Jac}(\CbeC^{-})(\R (x))$ into $\textrm{Jac}(\CbeC )(\R (x))$ with finite kernel and whose image contains $8\textrm{Jac}(\CbeC )(\R (x))$. Hence the Mordell-Weil rank of \linebreak$\textrm{Jac}(\CbeC )(\R (x))$ is the sum of Mordell-Weil ranks of $\textrm{Jac}(\CbeC^{+})(\R (x))$\linebreak and $\textrm{Jac}(\CbeC^{-})(\R (x))$.
\fin%



\begin{sprop}\label{pol-pair-theo-final} Let $k$ be a characteristic $0$ field and $f(x,y)\in k(x)[y]$ be a polynomial of odd degree in $y$. Denote by $\CbeC$ the hyperelliptic curve defined over $k(x)$ by the affine equation $z^{2}=f(x^{2},y).$ For each $\delta\in k(x)^{\times}$ denote by $\CbeC_{\delta}$ the $k(x)$-hyperelliptic curve given by the affine equation \linebreak$t^{2}=\delta^{\deg_{y}(f)}f\left( x,\frac{s}{\delta}\right) .$

Then the Mordell-Weil rank of $\textrm{Jac}(\CbeC )(k(x))$ is the sum of Mordell-Weil ranks of $\textrm{Jac}(\CbeC_{1})(k(x))$ and $\textrm{Jac}(\CbeC_{x})(k(x))$.
\end{sprop}
\dem%
Since the degree of $f$ in $y$ is odd, the curve $\CbeC$ and all the curves $\CbeC_{\delta}$ have a $k(x)$-rational point at infinity. In particular $\textrm{Jac}(\CbeC )(k(x))$ is isomorphic to $\textrm{Pic}^{0}(k(x)(\CbeC ))$ and $\textrm{Jac}(\CbeC_{\delta})(k(x))$ is isomorphic to $\textrm{Pic}^{0}(k(x)(\CbeC_{\delta})).$ 

Let $\iota$ be the hyperelliptic involution and $\rho$ be the involution of $k(x)(\CbeC )$ that preserves $k$, $y$ and $z$, and sends $x$ on $-x$. As in the proof of Proposition \ref{prop-decomp-jac} we prove that 
the maps 
$$\phi_{1}:\begin{array}[t]{ccc}
k(x)(\CbeC_{1})&\longrightarrow &k(x)(\CbeC )^{\rho}\\
A(x,s,t)&\longmapsto &A\left( x^{2},y,z\right) \\
\end{array}\textrm{ and }
\phi_{x}:\begin{array}[t]{ccc}
k(x)(\CbeC_{x})&\longrightarrow &k(x)(\CbeC )^{\iota\circ\rho}\\
A(x,s,t)&\longmapsto &A\left( x^{2},x^{2}y,x^{\deg_{y}(f)}z\right) \\
\end{array}$$ %
are isomorphisms. %
To prove Proposition \ref{pol-pair-theo-final} we apply Proposition \ref{decomp-pic-theo-gene}; the isomorphism $\phi_{1}$ (respectively $\phi_{x}$) helps us to identify $\textrm{Pic}^{0}(k(x)(\CbeC )^{\rho})$ and $\textrm{Jac}(\CbeC_{1})(k(x))$ (respectively $\textrm{Pic}^{0}(k(x)(\CbeC )^{\iota\circ\rho})$ and $\textrm{Jac}(\CbeC_{x})(k(x))$). 
\fin%
\begin{sprop}\label{coro-final-pol-pairs} Let $B,$ $C\in\R (x)$ be two rational fractions. Let $\CbeC$ be the hyperelliptic curve defined over $\R (x)$ by the affine equation 
$$z^{2}+(y^{2}+1)(y^{2}+C(x^{2}))(y^{4}+(1+C(x^{2}))y^{2}+B(x^{2}))=0.$$ 
Assume the polynomials $B(x^{2})$, $C(x^{2})$, $B(x^{2})-C(x^{2}),$ $C(x^{2})(B(x^{2})-C(x^{2}))$, $(B(x^{2})-C(x^{2}))(1-C(x^{2}))$ and $(1+C(x^{2}))^{2}-4B(x^{2})$ are not squares in $\C (x)$. For each $\delta\in\R (x)^{\times}$ we consider the two following $\R (x)$-hyperelliptic curves: 
$$\begin{array}{l}
\CbeC^{+}_{\delta}:z^{2}=y(y-\delta )(y-\delta C(x))(y^{2}-\delta [1+C(x)]y+\delta^{2}B(x))\textrm{ and}\\
\CbeC^{-}_{\delta}:z^{2}=y(y^{2}-\delta [(1-C(x))^{2}-2(B(x)-C(x))]y+\delta^{2}(B(x)-C(x))^{2}).\\
\end{array}$$
Then the Mordell-Weil rank of $\textrm{Jac}(\CbeC )(\R (x))$ is the sum of Mordell-Weil ranks of $\textrm{Jac}(\CbeC^{+}_{1})(\R (x))$, $\textrm{Jac}(\CbeC^{+}_{x})(\R (x))$, $\textrm{Jac}(\CbeC^{-}_{1})(\R (x))$ and $\textrm{Jac}(\CbeC^{-}_{x})(\R (x))$.
\end{sprop}
\dem %
Consider the following two $\R (x)$-hyperelliptic curves 
$$\begin{array}{l}
\CbeH :t^{2}+s(s+1)(s+C(x^{2}))(s^{2}+(1+C(x^{2}))s+B(x^{2}))=0\textrm{ and }\\
\CbeE :t^{2}+(s+1)(s+C(x^{2}))(s^{2}+(1+C(x^{2}))s+B(x^{2}))=0.\\
\end{array}$$
Following Proposition \ref{prop-decomp-jac}, the Mordell-Weil rank of $\textrm{Jac}(\CbeC )(\R (x))$ is the sum of Mordell-Weil ranks of $\textrm{Jac}(\CbeH )(\R (x))$, $\textrm{Jac}(\CbeE )(\R (x))$. 

Using the change of variables given by $s=-\alpha$ and $t=\beta$, we get an isomorphism between the curve $\CbeH$ and the $\R (x)$-hyperelliptic curve
$$\CbeC^{+}:\beta^{2}=\alpha (\alpha -1)(\alpha -C(x^{2}))\left(\alpha^{2}-[1+C(x^{2})]\alpha +B(x^{2})\right).$$
Considering the change of variables given by 
$$\left\{\begin{array}{l}
s:=-(B(x^{2})-C(x^{2}))\left[ (C(x^{2})-1)\frac{1}{y+1}+1\right] .\\
t:=(B(x^{2})-C(x^{2}))(1-C(x^{2}))\frac{z}{(y+1)^{2}}\\
\end{array}\right.$$ 
we get an isomorphism between the curve $\CbeE$ and the $\R (x)$-elliptic curve
$$\CbeC^{-}:\beta^{2}=\alpha\left(\alpha^{2}-[(1-C(x^{2}))^{2}-2(B(x^{2})-C(x^{2}))]\alpha +(B(x^{2})-C(x^{2}))^{2}\right) .$$
To conclude we apply Proposition \ref{pol-pair-theo-final} to $\CbeC^{+}$ and then to $\CbeC^{-}$. 
\fin%

\end{subsection}

\begin{subsection}{Doing a $2$-descent.}

\begin{subsubsection}{Christie's Lemma.}

\begin{ssprop}\label{def_K}%
Let $k_{0}$ be a subfield of $\C .$ Let $f\in k_{0}(x)[y]$ be a squarefree polynomial of odd degree and $\CbeC$ be the hyperelliptic curve defined over $k_{0}(x)$ by the affine equation $z^{2}=f(y)$. 
We assume that the $2$-primary torsion of $\textrm{Jac}(\CbeC )(\C (x))$ is finite. %
Then $\textrm{Jac}(\CbeC )(\C (x))$ is equal to $\textrm{Jac}(\CbeC )(K (x))$ for some finite extension $K$ of $k_{0}$.
\end{ssprop}
\dem%
For every $\C (x)$-point $P$ of $\textrm{Jac}(\CbeC )$ denote by $K_{P}$ the smallest subfield of $\C$ containing $k_{0}$ and such that $P$ is defined over $K_{P}(x)$. 

If $K_{P}$ is not a finite extension of $k_{0}$, then $K_{P}$ is a finite extension of $k(t_{1},\cdots ,t_{n})$ with $t_{1},\cdots ,t_{n}$ algebraically independant over $k_{0}$. In that case, by specialising $t_{1},\cdots ,t_{n}$ over $\C$, the point $P$ gives uncountably many \linebreak$\C$-points of $\textrm{Jac}(\CbeC )$. This is a contradiction because $\textrm{Jac}(\CbeC )(\C (x))$ is finitely generated (as in the proof of Corollary \ref{finitude-jac-CbeC-final}, apply Theorem \ref{Lang_Neron}). Thus $K_{P}$ is a finite extension of $k_{0}$.

The group $\textrm{Jac}(\CbeC )(\C (x))$ is generated by a finite family $(P_{i})_{i=1}^{r}$ (see Corollary \ref{finitude-jac-CbeC-final}). The smallest subfield $K$ of $\C$ containing all the $K_{P_{i}}$ is a finite extension of $k_{0}$ and $\textrm{Jac}(\CbeC )(K (x))$ contains all the points $P_{i}$. In particular $\textrm{Jac}(\CbeC )(\C (x))$ and $\textrm{Jac}(\CbeC )(K (x))$ are equal.
\fin\\

For a better understanding of the field $K$ defined by Proposition \ref{def_K}, we use a lemma from Christie (see 
\cite{Christie}):
\begin{ssprop}\label{lem_base_action_mod_2}
Let $\Gamma$ be a finite group and $\CbeA$ be a finitely generated free abelian group on which $\Gamma$ acts. Assume the triviality of the action of $\Gamma$ on $\CbeA /2\CbeA$. 

Then $\CbeA$ has a basis $(a_{i})_{i=1}^{t}$ such that $\tau (a_{i})\in\{ -a_{i},a_{i}\}$ for every $\tau\in\Gamma$.
\end{ssprop}

\begin{ssprop}\label{prop-lemme-descente-res-final-322} Let $k$ be a subfield of $\R .$ Let $f\in k(x)[y]$ be a polynomial of degree $2g+1$. Let $\CbeC$ be the hyperelliptic curve defined over $k(x)$ by the affine equation $z^{2}=f(y)$. Denote by $J$ the jacobian variety associated to $\CbeC .$ 
%
%
We assume that 
\begin{enumerate}
\item the $2$-primary torsion of $J(\C (x))$ is finite and
\item the action of $\textrm{Gal}(\C /k)$ on $J(\C (x))/2J(\C (x))$ is trivial.  %
\end{enumerate}
For each $d\in k^{\times}$ denote by $\CbeC_{d}$ the hyperelliptic curve defined over $k(x)$ by the affine equation $z^{2}=d^{2g+1}f(\frac{y}{d}).$ %

Then the Mordell-Weil rank of $J(\R (x))$ is $0$ if and only if for every positive element $d\in k^{\times}$ the $k(x)$-Mordell-Weil rank of $\textrm{Jac}(\CbeC_{d})$ is $0$.
\end{ssprop}
\dem%
Since the $2$-primary torsion of $J(\C (x)$ is finite, Corollary \ref{def_K} asserts the existence of a finite extension $K$ of $k$ such that $J(\C (x))=J(K(x))$. 

Assume that the Mordell-Weil rank of $J(\R (x))$ is different from $0$. 
The group $\Gamma :=\textrm{Gal}(K/k)$ is finite. Following Corollary \ref{finitude-jac-CbeC-final}, the group \linebreak$\CbeA :=J(K(x))/J(K(x))_{tors}$ is a finitely generated free abelian group. 
This means that the hypotheses of Proposition \ref{lem_base_action_mod_2} are satisfied. This proposition gives a basis $(\alpha_{i})_{i=1}^{t}$ of $\CbeA$ such that $\tau (\alpha_{i})\in\{ -\alpha_{i},\alpha_{i}\}$ for every $\tau\in\Gamma$.%

Let $\sigma$ be the complex conjugation. The field $k$ being a subfield of $\R ,$ the group $\Gamma =\textrm{Gal}(\C /k)$ contains $\sigma$. Since $\tau (\alpha_{i})\in\{ -\alpha_{i},\alpha_{i}\}$ for every $\tau\in\Gamma$, the action of $\Gamma$ commutes with $\sigma$. Thus $\Gamma$ acts on the subgroup $\CbeA^{\sigma}$ of \linebreak$\sigma$-invariants elements of $\CbeA$. The action of $\Gamma$ on $\CbeA /2\CbeA$ is trivial. The group $\CbeA$ being a free abelian group, the intersection $\CbeA^{\sigma}\cap 2\CbeA$ is equal to $2\CbeA^{\sigma}$. 
This implies the triviality of the action of $\Gamma$ on $\CbeA^{\sigma}/2\CbeA^{\sigma}$. Hence, the hypotheses of Proposition \ref{lem_base_action_mod_2} are satisfied and its application to $\CbeA^{\sigma}$ and $\Gamma$ gives a basis $(a_{i})_{i=1}^{t}$ of $\CbeA^{\sigma}$ such that $\tau (a_{i})\in\{ -a_{i},a_{i}\}$ for every $\tau\in\Gamma$.

Let $P_{i}\in J(K(x))$ be an element of the class $a_{i}.$ Let $m$ be the exponent of $J(K(x))_{tors}.$ 
The point $mP_{i}$ is fixed by a subgroup $\Gamma_{i}$ of $\Gamma$ of index at most $2$. The degree of the field $K^{\Gamma_{i}}$ of elements in $K$ invariants under the action of $\Gamma_{i}$ is an extension of $k$ of degree at most $2$ i.e. $K^{\Gamma_{i}}=k(\sqrt{d_{i}})$ for some $d_{i}\in k^{\times}$. 
Since the complex conjugation $\sigma$ belongs to $\Gamma_{i}$, the field $k(\sqrt{d_{i}})$ is contained in $\R$. 
In particular $d_{i}$ is positive. 

If $\Gamma_{i}=\Gamma$, 
then $mP_{i}$ 
is an element of $J(k(x))=\textrm{Jac}(\CbeC_{1})(k(x))$ of infinite order. 
Assume the existence of $\tau_{i}\in\Gamma$ such that $\tau_{i}(a_{i})=-a_{i}$. 
Then $d_{i}$ is not a square in $k$. The degree of $f$ being odd, the curves $\CbeC$ and $\CbeC_{d_{i}}$ have a $k(x)$-rational point above the point at infinity of $\Pun.$ This implies the following isomorphims: 
\begin{itemize}
\item[* ] $\textrm{Jac}(\CbeC )(K^{\Gamma_{i}}(x))\simeq\textrm{Pic}^{0}(K^{\Gamma_{i}}(x)(\CbeC )),$ 
\item[* ] $\textrm{Jac}(\CbeC_{1})(k(x))\simeq\textrm{Pic}^{0}(k(x)(\CbeC_{1}))$ and 
\item[* ] $\textrm{Jac}(\CbeC_{d_{i}})(k(x))\simeq\textrm{Pic}^{0}(k(x)(\CbeC_{d_{i}})).$
\end{itemize}%
Let $\iota$ be the hyperelliptic involution. The subfield of $\tau_{i}$-invariants elements of $K^{\Gamma_{i}}(x)(\CbeC )$ is $k(x)(\CbeC ).$ The image of 
$\phi:\begin{array}[t]{ccc}
k(x)(\CbeC_{d_{i}})&\longrightarrow &K^{\Gamma_{i}}(x)(\CbeC )\\
A(s,t)&\longmapsto &A(d_{i}y,d_{i}^{g}\sqrt{d_{i}}z)\\
\end{array}$ is the subfield of $\iota\circ\tau_{i}$-invariants elements of $K^{\Gamma_{i}}(x)(\CbeC )$. 
Thus, as $2mP_{i}$ is of infinite order and belongs to $2\textrm{Jac}(\CbeC )(K^{\Gamma_{i}}(x))$, 
Proposition \ref{decomp-pic-theo-gene} asserts the existence of an element of infinite order in 
$\textrm{Jac}(\CbeC_{1})(k(x))\times\textrm{Jac}(\CbeC_{d_{i}})(k(x))$ %
i.e. that the Mordell-Weil rank of either $\textrm{Jac}(\CbeC_{1})(k(x))$ or $\textrm{Jac}(\CbeC_{d_{i}})(k(x))$ is at least $1$.

Conversely Assume the existence of $d\in k^{\times}$ positive such that \linebreak$\textrm{Jac}(\CbeC_{d})(k(x))$ has an infinite order element $a$. From Proposition \ref{decomp-pic-theo-gene} we know the existence of a morphism  
$\varphi :\textrm{Jac}(\CbeC_{d})(k(x))\longrightarrow\textrm{Jac}(\CbeC )(k(\sqrt{d})(x))$ %
of finite kernel. 
The image $\varphi (a)$ is an infinite element of $\textrm{Jac}(\CbeC )(k(\sqrt{d})(x))$ and thus of $\textrm{Jac}(\CbeC )(\R (x))$ 
(notice that $k(\sqrt{d})\subset\R$ since $d$ is positive).
\fin

\end{subsubsection}

\begin{subsubsection}{A first study of the image of $\pi_{\CbeC}$.}

\begin{ssnota}\label{nota-image-phi}
\rm Let $k$ be a characteristic $0$ field. For each monic polynomial $P(y)\in k(x)[y]$ 
Denote by $K_{P}$ the algebra $k(x)[y]/(P(y))$ and by $y_{P}$ the class of $y$ in $K_{P}$.
\end{ssnota}

\begin{ssnota}\label{nota-image-phi-bis}
\rm Let $f\in k[x][y]$ be a squarefree monic polynomial of odd degree 
and let $\CbeC$  be the hyperelliptic curve defined over $k(x)$ by the affine equation $\CbeC :z^{2}=f(y).$ %
%
%
\rm Let $f(y)=\displaystyle\prod_{l\in\widetilde{I}}\mu_{l}(y)$ be the decomposition of $f(y)$ into monic prime elements of $k(x)[y]$. 
For each $l\in\widetilde{I}$ we assume that $\mu_{l}$ belongs to $k[x][y]$. 
\rm Let $f'(y)$ be the usual derivative of $f(y)$. 
For each $l\in\widetilde{I}$ denote by $T_{l}$ the class $f'(y_{\mu_{l}})$ of $f'(y)$ in $K_{\mu_{l}}=k(x)[y]/(\mu_{l}(y)).$
\end{ssnota}


\begin{ssprop}\label{caract-image-CS-v-ideaux}\label{image_Phi1}
We use the notation \ref{nota-image-phi} and \ref{nota-image-phi-bis}. Let $l$ be an element of $\widetilde{I}.$ Let $\textrm{div}(u,v)\in\textrm{Div}^{0}(\CbeC )(k(x))$ be a semi-reduced divisor such that $u$ is coprime to $f$. 

Then the finite places of $K_{\mu_{l}}$ at which $u(y_{\mu_{l}})$ has odd valuation are 
in the support $\textrm{Supp}_{K_{\mu_{l}}}(T_{l})$ of $\textrm{div}(T_{l}).$ 
\end{ssprop}

\begin{ssnota}\label{nota-image-phi-2}


\rm We use the notation of Proposition \ref{caract-image-CS-v-ideaux}. Let \linebreak$u=\displaystyle\prod_{i\in I}p_{i}^{n_{i}}$ be the decomposition of $u$ into monic prime elements of $K_{\mu_{l}}[y]$ (it exists since $u$ is monic). 
\rm Consider an index $i\in I$. Denote by $K_{p_{i},\mu_{l}}$ the field $K_{\mu_{l}}[y]/(p_{i}(y))$ 
and by $y_{p_{i}}$ 
the class of $y$ in $K_{p_{i},\mu_{l}}$. 
%
\end{ssnota}

\begin{sslemm}\label{lemme1} We use the notation \ref{nota-image-phi}, \ref{nota-image-phi-bis} and \ref{nota-image-phi-2}. Let $\mathfrak{p}$ be a finite place of $K_{\mu_{l}}$ such that $v_{\mathfrak{p}}(p_{i}(y_{\mu_{l}}))\neq 0.$ Let $\mathcal{P}$ be a place of $K_{p_{i},\mu_{l}}$ above $\mathfrak{p}$ such that $v_{\mathcal{P}}(T_{l})=0.$ Let $f_{j}$ be the coefficient of $f(y)=f((y-y_{\mu_{l}})+y_{\mu_{l}})$ relative to the monomial $(y-y_{\mu_{l}})^{j+1}$. %
We assume that $v_{\mathcal{P}}(y_{p_{i}}-y_{\mu_{l}})\neq 0$. %

Then $v_{\mathcal{P}}\left( T_{l}+(y_{p_{i}}-y_{\mu_{l}})^{2g}+\displaystyle\sum_{j=1}^{2g-1}f_{j}(y_{p_{i}}-y_{\mu_{l}})^{j}\right)$ is even.
\end{sslemm}
\dem%
By definition of $y_{\mu_{l}}$, the minimal polynomial of $y_{\mu_{l}}$ over $k(x)$ is $\mu_{l}$. Since this polynomial belongs to $k[x][y]$, the element $y_{\mu_{l}}$ is integral over $k[x]$. Hence $v_{\mathcal{P}}(y_{\mu_{l}})$ is positive or equal to $0$ (see \cite{Stichtenoth} Proposition III.3.1).


The coefficients $f_{j}$ are polynomials in $y_{\mu_{l}}$ and in some elements of $k[x]$ (the coefficients of $f(y)$ relatives to the monomials $y^{r}$). Thus $v_{\mathcal{P}}(f_{j})$ is positive or equal to $0$. 
\begin{description}
\item[Case 1: If $v_{\mathcal{P}}(y_{p_{i}}-y_{\mu_{l}})>0$. ] Since $v_{\mathcal{P}}(f_{j})\ge 0$ and $v_{\mathcal{P}}(T_{l})=0$, the triangle inequality show that 
$$v_{\mathcal{P}}\left( T_{l}+(y_{p_{i}}-y_{\mu_{l}})^{2g}+\displaystyle\sum_{j=1}^{2g-1}f_{j}(y_{p_{i}}-y_{\mu_{l}})^{j}\right) =v_{\mathcal{P}}(T_{l})=0.$$ 
\item[Case 2: If $v_{\mathcal{P}}(y_{p_{i}}-y_{\mu_{l}})<0$. ] Since $v_{\mathcal{P}}(f_{j})\ge 0$ and $v_{\mathcal{P}}(T_{l})=0$, the triangle inequality show that
$$\begin{array}{rcl}
v_{\mathcal{P}}\left( T_{l}+(y_{p_{i}}-y_{\mu_{l}})^{2g}+\displaystyle\sum_{j=1}^{2g-1}f_{j}(y_{p_{i}}-y_{\mu_{l}})^{j}\right) &=&v_{\mathcal{P}}((y_{p_{i}}-y_{\mu_{l}})^{2g})\\
&=&2gv_{\mathcal{P}}(y_{p_{i}}-y_{\mu_{l}}).\textrm{ \fin}\\
\end{array}$$
\end{description}

\begin{sslemm}\label{lemme_parite}We use the notation 
\ref{nota-image-phi}, \ref{nota-image-phi-bis} and \ref{nota-image-phi-2}. Let $\mathfrak{p}$ be a finite place of $K_{\mu_{l}}$ such that $v_{\mathfrak{p}}(p_{i}(y_{\mu_{l}}))\neq 0.$ Let $\mathcal{P}$ be a place of $K_{p_{i},\mu_{l}}$ above $\mathfrak{p}$ such that $v_{\mathcal{P}}(T_{l})=0.$ 
Then the valuation $v_{\mathcal{P}}(y_{p_{i}}-y_{\mu_{l}})$ is even.
\end{sslemm}
\dem%
Assume $v_{\mathcal{P}}(y_{p_{i}}-y_{\mu_{l}})$ is different from $0$ (if $v_{\mathcal{P}}(y_{p_{i}}-y_{\mu_{l}})=0$ the result is direct). 
From our hypotheses we know the coprimality of $p_{i}$ and $f$. Following the definition of Mumford's representation we have \linebreak$f(y)\equiv v(y)^{2}\bmod p_{i}(y)$. In particular 
$v_{\mathcal{P}}(f(y_{p_{i}}))$ is even. Since $\mu_{l}(y)$ divides $f(y),$ the element $f(y_{\mu_{l}})$ is equal to $0$. Taylor's formula gives 
\begin{equation}\label{congruence-base-lemme_parite}
f(y_{p_{i}})= (y_{p_{i}}-y_{\mu_{l}})\left( T_{l}+(y_{p_{i}}-y_{\mu_{l}})^{2g}+(\displaystyle\sum_{j=1}^{2g-1}f_{j}(y_{i}-y_{\mu_{l}})^{j})\right) 
\end{equation}
Lemma \ref{lemme_parite} is obtained by applying the parity of $v_{\mathcal{P}}(f(y_{p_{i}}))$ and of \linebreak$v_{\mathcal{P}}( T_{l}+(y_{p_{i}}-y_{\mu_{l}})^{2g}+\displaystyle\sum_{j=1}^{2g-1}f_{j}(y_{p_{i}}-y_{\mu_{l}})^{j})$ (see Lemma \ref{lemme1}) to Equation \ref{congruence-base-lemme_parite}. 
\fin\\

\noindent\textbf{Proof of Proposition \ref{image_Phi1}. } %
We use the notation \ref{nota-image-phi-2}. 
Let $\mathfrak{p}$ be a place of $K_{\mu_{l}}$ at which $p_{i}(y_{\mu_{l}})$ has odd valuation. Assume that $v_{\mathfrak{p}}(T_{l})$ is equal to $0$. 
For each place $\mathcal{P}$ of $K_{p_{i},\mu_{l}}$ above $\mathfrak{p}$, the valuation $v_{\mathcal{P}}(T_{l})=e(\mathcal{P}|\mathfrak{p})v_{\mathfrak{p}}(T_{l})$ is equal to $0$. In particular the hypotheses of Lemma \ref{lemme_parite} are satisfied. Following this lemma, for each place $\mathcal{P}$ above $\mathfrak{p}$, the valuation $v_{\mathcal{P}}(y_{p_{i}}-y_{\mu_{l}})$ is even. 
A classical computation shows that $v_{\mathfrak{p}}(N_{K_{p_{i},\mu_{l}}/K_{\mu_{l}}}(y_{p_{i}}-y_{\mu_{l}}))$ is equal to $\displaystyle\sum_{\mathcal{P}\textrm{ place of }K_{p_{i},\mu_{l}},\mathcal{P}|\mathfrak{p}}f(\mathcal{P}|\mathfrak{p})v_{\mathcal{P}}(y_{p_{i}}-y_{\mu_{l}})$ (see \cite{Zariski-Samuel} for the Dedekind rings case). Thus $v_{\mathfrak{p}}(p_{i}(y_{\mu_{l}}))$ is even. This is a contradiction with the choice of $\mathfrak{p}.$
\fin

\end{subsubsection}

\begin{subsubsection}{Our choice for the constant field.}

\begin{ssprop}\label{Im-phi-R-k-descente} Let $k$ be a subfield of $\R$. 
Let $f(y)\in k[x][y]$ be a squarefree monic polynomial of odd degree $2g+1$. Denote by $\CbeC$ the hyperelliptic curve defined over $k(x)$ by the affine equation $z^{2}=f(y).$ We assume the existence of $2g$ elements $e_{1},\cdots ,e_{2g-1},$ $H\in k[x]$ and of a polynomial $\mu (y)\in k[x][y]$ of degree $2$ such that $f(y)=\mu (y)\displaystyle\prod_{i=1}^{2g-1}(y-He_{i})$. Let us also assume that:
\begin{itemize}
\item[* ] the discriminant $\Delta (f)$ of $f(y)$ splits into linear factors over $k$,
\item[* ] the discriminant $\Delta (\mu )$ of $\mu$ is equal to $H^{2}Q^{2}D$ with $D\in k[x]$ a polynomial of degree $1$ and $Q\in k[x]$,
\item[* ] 
$\Delta (f)=Q^{2}Q_{1}$ with $Q_{1}\in k[x]$ coprime to $Q$, and 
\item[* ] $D(\alpha )$ is a square in $k$ for every root $\alpha\in k$ of $H$. 
\end{itemize}
Let $L$ be the algebra $\C (x)[t]/(f(t))$. Let $\pi_{\CbeC}:\textrm{Jac}(\CbeC )(\C (x))\longrightarrow\displaystyle L^{\times}/L^{\times 2}$ be the morphism defined by Proposition \ref{def-Cas-Sch+im-in-ker} (and relative $\textrm{Jac}(\CbeC )(\C (x))$). 
%
Then the action of $\textrm{Gal}(\C /k)$ on the image of $\pi_{\CbeC}$ is trivial. 
\end{ssprop}

\begin{sslemm}\label{Im-phi-R-k-descente-lemme-factoriel} We keep the notation and hypotheses of Proposition \ref{Im-phi-R-k-descente}. We assume that $\mu$ is irreducible. 
Denote by $K_{\mu ,\C}$ the algebra $\C (x)[y]/(\mu (y))$ and by $y_{\mu}$ the class of $y$ in $K_{\mu ,\C}$. Denote by $s$ the element $\frac{\mu'(y_{\mu})}{2HQ}$. 

Then the minimal polynomial of $s$ over $\C (x)$ is $y^{2}-D(x)$. Thus $\C [x,s]$ is a unique factorization domain and its fractions field is $K_{\mu ,\C}$.
\end{sslemm}

\begin{sslemm}\label{image-phi-sigma-invariance-HD} We keep the notation and hypotheses of Lemma \ref{Im-phi-R-k-descente-lemme-factoriel}. 
Let $\alpha\in k$ be a root of the resultant $\textrm{Res}_{T}(f'(T),\mu (T))$ such that $Q(\alpha )\neq 0$. Let $\beta$ be a prime element of $\C [x,s]$ such that $N_{K_{\mu ,\C}/\C (x)}(\beta )=\lambda (x-\alpha)$ for some constant $\lambda\in \C$. %

Then the valuation $v_{\beta}$ is invariant under the action of $\textrm{Gal}(\C /k)$. In particular for every semi-reduced divisor $\textrm{div}(u,v)\in\textrm{Div}^{0}(\C (x)(\CbeC ))$ with $u$ coprime to $f$ and for every $\sigma\in\textrm{Gal}(\C /k)$ the valuation $v_{\beta}(u(y_{\mu})\sigma (u(y_{\mu})))$ is even.
\end{sslemm}
\dem%
Since it belongs to $\C[x,s]$, the element $\beta$ can be written as \linebreak$\beta =\beta_{1}s+\beta_{0}$ with $\beta_{0},$ $\beta_{1}\in\C [x]$. The degree of $D$ is $1$ and $\lambda (x-\alpha )$ is equal to $\N_{\C (x)(s)/\C (x)}(\beta )=\beta_{0}^{2}-\beta_{1}^{2}D$. Thus $\beta_{0}$ and $\beta_{1}$ are in $\C$ and $\beta_{0}^{2}=\beta_{1}^{2}D(\alpha )$. 

The resultant $\textrm{Res}_{T}(f'(T),\mu (T))$ is equal to $\Delta (\mu )\displaystyle\prod_{i=1}^{2g-1}\mu (He_{i})$ with \linebreak$\Delta (\mu )=H^{2}Q^{2}D$. Hence $\alpha$ is either a root of $H$ or a root of $D$ or a root of $\displaystyle\prod_{i=1}^{2g-1}\mu (He_{i})$.
\begin{description}
\item[Case $1$: if $H (\alpha )\neq 0$ and $\mu (He_{i})(\alpha )= 0$ for some $i\in\{ 1,\cdots ,2g-1\}$. ]\textrm{ }\linebreak Applying Taylor's formula to $\mu (T)$ at $He_{i}$ we get the equality 
$$\mu (T)=(T-He_{i})^{2}+(T-He_{i})\mu '(He_{i})+\mu (He_{i}).$$
From this equality we deduce $\Delta (\mu )$ : it is equal to $\left(\mu '(He_{i})\right)^{2}-4\mu (He_{i})$.\linebreak In particular, since $\mu (He_{i})(\alpha )=0$, we deduce from it that \linebreak$\left( H(\alpha )Q(\alpha )\right)^{2}D(\alpha ) =\Delta (\mu )(\alpha )$ is a square in $k$. 
The element $H(\alpha )Q(\alpha )$ being nontrivial, this means that $D(\alpha )$ is a square in $k$ i.e. that $\frac{\beta_{0}}{\beta_{1}}$ is in $k$. In particular $v_{\beta}$ is invariant under the action of $\textrm{Gal}(\C /k)$.
\item[Case 2: $\alpha$ is a root of $D$. ] Then $\beta_{0}=0$ is in $k$ and thus $v_{\beta}$ is invariant under the action of $\textrm{Gal}(\C /k)$.
\item[Case 3: $\alpha$ is a root of $H$. ] Then, by hypothesis, $D(\alpha )$ is a square in $k$ i.e. $\frac{\beta_{0}}{\beta_{1}}$ belongs to $k$, and thus $v_{\beta}$ is invariant under the action of $\textrm{Gal}(\C /k)$.\fin%
\end{description}

\begin{sslemm}\label{image-phi-sigma-invariance-Q} We keep the notation and hypotheses of Lemma \ref{Im-phi-R-k-descente-lemme-factoriel}. Let $\alpha\in k$ be a root of $Q$ and $\beta$ be a prime element of $\C [x,s]$ such that $N_{K_{\mu ,\C}/\C (x)}(\beta )=\lambda (x-\alpha)$ for some constant $\lambda\in \C$. 

Let $\textrm{div}(u,v)\in\textrm{Div}^{0}(\C (x)(\CbeC ))$ be a semi-reduced divisor with $u$ coprime to $f$ and let $\sigma$ be an element of $\textrm{Gal}(\C /k)$. 
Then 
the valuation $v_{\beta}(u(y_{\mu})\sigma (u(y_{\mu})))$ is even. 
\end{sslemm}
\dem%
The polynomial $\mu$ admits $y_{\mu}$ as a root in $K_{\mu ,\C}$ and its degree is $2$. Thus $\mu$ is totally split in $K_{\mu ,\C}$. Denote by $\iota$ the unique \linebreak$\C(x)$-automorphism of $K_{\mu ,\C}=\C (x)(y_{\mu})$ sending $y_{\mu}$ on the other root of $\mu$.  

Since $\beta\iota (\beta )=N_{\C (x,s)/\C (x)}(\beta )=\lambda (x-\alpha)$ the set of prime factors of $\lambda (x-\alpha )$ is $\{\beta ,\iota (\beta )\}$. Moreover $x-\alpha$ is $\sigma$-invariant. Hence the set of prime factors of $\lambda (x-\alpha )$ is also $\{\sigma^{-1}(\beta ),\sigma^{-1}\circ\iota (\beta )\}$. In particular $v_{\sigma^{-1}(\beta )}$ is equal to either $v_{\beta}$ or $v_{\iota (\beta )}$. When $v_{\beta}=v_{\sigma^{-1}(\beta )}$ the result is straightforward. Assume that $v_{\sigma^{-1}(\beta )}=v_{\iota (\beta )}=v_{\iota^{-1}(\beta )}$.

The polynomial $f'(He_{i})$ divides $Q_{1}$ and $Q_{1}$ is coprime to $Q$. As a consequence $f'(He_{i})$ is coprime to $x-\alpha$. Following Proposition \ref{caract-image-CS-v-ideaux} this implies that the valuation $v_{x-\alpha}(u(He_{i})$ is even. In particular the valuation $v_{\beta}(u(He_{i}))=e(\beta |x-\alpha )v_{x-\alpha}(u(He_{i}))$ is even. 

Denote by $K_{u,\mu ,\C}$ the algebra $K_{\mu ,\C}[y]/(u(y))$. By definition of Mumford's representation, $f(y)=(y-y_{\mu})(y-\iota (y_{\mu}))\displaystyle\prod_{i=1}^{2g-1}(y-He_{i})$ is a square mo\-du\-lo $u$. In particular $N_{K_{u,\mu ,\C}/K_{\mu ,\C}}(f(y))=(-1)^{\deg (u)}u(y_{\mu})u(\iota (y_{\mu}))\displaystyle\prod_{i=1}^{2g-1}u(He_{i})$ is a square in $K_{\mu ,\C}$. In particular its valuation at $\beta$ is even. As $v_{\beta}(u(He_{i}))$ is even, we get that $v_{\beta}(u(y_{\mu})u(\iota (y_{\mu}))$ is even. This is enough to conclude since $v_{\sigma^{-1}(\beta )}=v_{\iota^{-1}(\beta )}$. In fact we have 
$$\begin{array}{rcl}
v_{\beta}(u(y_{\mu})u(\iota (y_{\mu}))&=&v_{\beta}(u(y_{\mu}))+v_{\iota^{-1}(\beta )}(u(y_{\mu}))\\
&=&v_{\beta}(u(y_{\mu}))+v_{\sigma^{-1}(\beta )}(u(y_{\mu}))\\
&=&v_{\beta}(u(y_{\mu})\sigma (u(y_{\mu}))\textrm{ \fin}\\
\end{array}$$

\noindent\textbf{Proof of Proposition \ref{Im-phi-R-k-descente}.}
For every prime factor $p\in\C (x)[y]$ of $f$ denote by $K_{p ,\C}$ the field $\C (x)[y]/(p(y))$ and by $y_{p}$ the class of $y$ in $K_{p,\C}$. Under hypotheses of Proposition \ref{Im-phi-R-k-descente} the field $K_{p,\C}$ is the fraction field of a unique factorization domain $\mathcal{O}_{p,\C}$ (see Lemma \ref{Im-phi-R-k-descente-lemme-factoriel}). %

By definition of $\pi_{\CbeC}$ Proposition \ref{Im-phi-R-k-descente} is proved if we show that for every semi-reduced divisor $\textrm{div}(u,v)\in\textrm{Div}^{0}(\C (x)(\CbeC ))$ with $u$ coprime to $f$, for every prime factor $p$ of $f$ and for every $\sigma\in\textrm{Gal}(\C /k)$ the class of $u(y_{p})\sigma (u(y_{p}))$ in $K_{p,\C}$ is a square. 

Let $\textrm{div}(u,v)\in\textrm{Div}^{0}(\C (x)(\CbeC ))$ be a semi-reduced divisor with $u$ coprime to $f$, let $p$ be a prime factor of $f$ and let $\sigma$ be an element of $\textrm{Gal}(\C /k)$. Since every element of $\C$ is a square and since $\mathcal{O}_{p,\C}$ is a unique factorization domain we prove Proposition \ref{Im-phi-R-k-descente} if we can show that $v_{\beta}(u(y_{p})\sigma (u(y_{p})))$ is even for every prime $\beta\in\mathcal{O}_{p,\C}$.

Assume the existence of a prime element $\beta\in\mathcal{O}_{p,\C}$ such that \linebreak$v_{\beta}(u(y_{p})\sigma (u(y_{p})))=v_{\beta}(u(y_{p}))+v_{\sigma^{-1}(\beta )}(u(y_{p}))$ is odd. Eventually replacing $\beta$ by $\sigma^{-1}(\beta )$ we can assume that $v_{\beta}(u(y_{p}))$ is odd and that $v_{\sigma^{-1}(\beta )}(u(y_{p}))$ is even. Following Proposition \ref{caract-image-CS-v-ideaux} the norm $N_{K_{p,\C}/\C (x)}(\beta )$ is a divisor of $N_{K_{p,\C}/\C (x)}(f'(y_{p}))=\textrm{res}_{T}(f'(T),p(T))$. In particular $N_{K_{p,\C}/\C (x)}(\beta )$ divides $\Delta (f)$. Since $\Delta (f)$ splits into linear factors in $k[x]$ the norm $N_{K_{p,\C}/\C (x)}(\beta )$ is equal to $\lambda (x-\alpha )$ for some $\lambda\in\C^{\times}$ and some $\alpha\in k$.
\begin{description}
\item[Case 1: if the degree of $p$ is $1$. ] Then $\beta$ is equal to $N_{K_{p,\C}/\C (x)}(\beta )$ i.e. to $\lambda (x-\alpha )$. In particular the valuation $v_{\beta}$ and $v_{\sigma (\beta )}$ are equals. This is in contradiction with the definition of $\beta$.
\item[Case 2: if $p=\mu$ is irreducible of degree $2$. ] We apply Lemma \ref{image-phi-sigma-invariance-HD} and Lemma \ref{image-phi-sigma-invariance-Q}. 
\fin%
\end{description}

\begin{sscoro}\label{prop-2desc-coeff-inter} We use the notation and hypotheses of Proposition \ref{coro-final-pol-pairs}. 
For each $\delta\in k(x)$ consider the two $k (x)$-hyperelliptic curves: 
$$\begin{array}{l}
\CbeC^{+}_{\delta}:z^{2}=y(y-\delta )(y-\delta C(x))(y^{2}-\delta [1+C(x)]y+\delta^{2}B(x))\textrm{ and}\\
\CbeC^{-}_{\delta}:t^{2}=s(s^{2}-\delta [(1-C(x))^{2}-2(B(x)-C(x))]s+\delta^{2}(B(x)-C(x))^{2}).\\
\end{array}$$
%
Let $k$ be a subfield of $\R$ containg the coefficients of $B$ and $C$. We assume 
\begin{itemize}
\item[* ] that $B(x),$ $C(x),$ $B(x)-C(x)$ and $1-C(x)$ split over $k$ into linear factors,
\item[* ] that $(1+C(x))^{2}-4B(x)$ has degree $1$ and that its value at $0$ is a square in the field $k$, 
\item[* ] that $1-C$ is coprime to $x$, $B$, $B-C$ and $(1+C)^{2}-4B$.
\end{itemize}
Then the $\R (x)$-Mordell-Weil rank of $\textrm{Jac}(\CbeC^{+}_{1}),$ $\textrm{Jac}(\CbeC^{+}_{x}),$ $\textrm{Jac}(\CbeC^{-}_{1})$ and $\textrm{Jac}(\CbeC^{-}_{x})$ are trivials if and only if for every positive element $\zeta\in k^{\times}$ the $k(x)$-Mordell-Weil ranks of $\textrm{Jac}(\CbeC^{+}_{\zeta}),$ $\textrm{Jac}(\CbeC^{+}_{\zeta x}),$ $\textrm{Jac}(\CbeC^{-}_{\zeta})$ and $\textrm{Jac}(\CbeC^{-}_{\zeta x})$ are trivials.
\end{sscoro}
\dem Following Proposition \ref{coro-final-pol-pairs} the groups $\textrm{Jac}(\CbeC^{+}_{1})(\C (x)),$ $\textrm{Jac}(\CbeC^{+}_{x})(\C (x)),$ $\CbeC^{-}_{1}(\C (x))$ et $\CbeC^{-}_{x}(\C (x))$ are finitely generated. In particular their $2$-primary torsion subgroups are finite. For every $\delta$ the discriminant of 
$$y^{2}+\delta ((1-C(x))^{2}-2(B(x)-C(x)))y+\delta^{2}(B(x)-C(x))^{2}$$
is $\delta^{2}(1-C(x))^{2}((1+C(x))^{2}-4B(x))$ and the discriminant of 
$$y(y^{2}+\delta ((1-C(x))^{2}-2(B(x)-C(x)))y+\delta^{2}(B(x)-C(x))^{2})$$
is $\delta^{6}(B(x)-C(x))^{4}(1-C(x))^{2}((1+C(x))^{2}-4B(x)).$ Thus Proposition \ref{Im-phi-R-k-descente} can be applied to the curves $\CbeC^{-}_{1}$ and $\CbeC^{-}_{x}$. Using Proposition \ref{Im-phi-R-k-descente} we prove that the hypotheses for the application of Proposition \ref{prop-lemme-descente-res-final-322} to the curves $\CbeC^{-}_{1}$ and $\CbeC^{-}_{x}$ are satisfied. 

In the same way we prove that the hypotheses for the application of Proposition \ref{prop-lemme-descente-res-final-322} to the curves $\CbeC^{+}_{1}$ and $\CbeC^{+}_{x}$ are satisfied by using Proposition \ref{Im-phi-R-k-descente} and noticing that for every $\delta$ the discriminant of the polynomial $y^{2}-\delta (1+C(x))y+\delta^{2}B(x)$ is $\delta^{2}[(1+C(x))^{2}-4B(x)]$ and the discriminant of the polynomial $y(y-\delta )(y-\delta C(x))(y^{2}-\delta (1+C(x))y+\delta^{2}B(x))$ is $\delta^{20}B(x)^{2}C(x)^{2}(C(x)-1)^{2}(B(x)-C(x))^{4}((1+C(x))^{2}-4B(x)))$.

Corollary \ref{prop-2desc-coeff-inter} is a direct consequence of Proposition \ref{prop-lemme-descente-res-final-322}.
\fin%

\begin{sscoro}\label{stheo-3-3-final-desc-finie} 
We use the notation \ref{Notations_Principales_21} 
and we assume that all the following elements are different from $0$:
\begin{itemize}
\item[* ] $\eta$, $\omega ,$ $\rho$ and $\omega^{2}-\eta^{2}$ and 
\item[* ] $2b_{1}-2+\omega^{2}-\eta^{2}$,
\item[* ] $\omega^{2}-\eta^{2}-2+2\eta$ and $\omega^{2}-\eta^{2}-2-2\eta$,
\item[* ] $\omega^{2}-\eta^{2}+2\omega$ and $\omega^{2}-\eta^{2}-2\omega$, 
\item[* ] $\omega^{2}-\eta^{2}-1+2\eta$ and  $\omega^{2}-\eta^{2}-1-2\eta$,
\item[* ] $2b_{1}+\omega^{2}-\eta^{2}-1,$ $b_{1}+\eta$, $b_{1}-\eta$, $b_{1}-1+\omega$ and $b_{1}-1-\omega$
\end{itemize}
%
%
For each $\delta\in k(x)^{\times}$ denote respectively by $\CbeC^{+}_{\delta}$ and $\CbeC^{-}_{\delta}$ the two following $k (x)$-hyperelliptic curves: 
$$\begin{array}{l}
\CbeC^{+}_{\delta}:t^{2}=y(y-\delta )(y-\delta C(x))(y^{2}-\delta [1+C(x)]y+\delta^{2}B(x))\textrm{ and}\\
\CbeC^{-}_{\delta}:t^{2}=s(s^{2}-\delta [(1-C(x))^{2}-2(B(x)-C(x))]s+\delta^{2}[B(x)-C(x)]^{2}).\\
\end{array}$$ %
Then the $\R (x)$-Mordell-Weil rank of $\textrm{Jac}(\CbeC)$ is $0$ if and only if for every positive element $\zeta\in k^{\times}$ the $k(x)$-Mordell-Weil ranks of $\textrm{Jac}(\CbeC^{+}_{\zeta}),$ $\textrm{Jac}(\CbeC^{+}_{\zeta x}),$ $\textrm{Jac}(\CbeC^{-}_{\zeta})$ and $\textrm{Jac}(\CbeC^{-}_{\zeta x})$ are trivial. 
\end{sscoro}
\dem%
Corollary \ref{prop-2desc-coeff-inter} can be applied since 
\begin{enumerate}
\item $B=(x+b_{1})^{2}-\eta^{2}$, $C=2(x+b_{1})+\omega^{2}-\eta^{2}-1$, $1-C$ and \linebreak$B-C=(x+b_{1}-1)^{2}-\omega^{2}$ split into linear fatcors;
\item $(1+C(x))^{2}-4B(x)=4(\omega^{2}-\eta^{2})x+4\rho^{2}$ has degree $1$;
\item $1-C$ is coprime to $xB(x)(B(x)-C(x))((1+C(x))^{2}-4B(x))$;
\item the following polynomials are not squares in $\C (x)$: 
\begin{enumerate}
\item $B(x^{2})$ (since $\eta$, $b_{1}-\eta$ and $b_{1}+\eta$ are nontrivial), 
\item $B(x^{2})-C(x^{2})$ (since $\omega$, $b_{1}-1+\omega$ and $b_{1}-1+\omega$ are nontrivial), 
\item $C(x^{2})$ (since $2b_{1}+\omega^{2}-\eta^{2}-1$ is nontrivial) and 
\item $(1+C(x^{2}))^{2}-4B(x^{2})=4(\omega^{2}-\eta^{2})x^{2}+4\rho^{2}$ (since $\rho\neq 0$). 
\end{enumerate}
Since $\textrm{Gcd}(C(x^{2}),B(x^{2}))=\textrm{Gcd}(C(x^{2}),B(x^{2})-C(x^{2})=1$, the polynomials $B(x^{2})C(x^{2})$ and $C(x^{2})(B(x^{2})-C(x^{2}))$ are not squares in $\C (x)$. In the same way, $(1-C(x^{2})(B(x^{2})-C(x^{2}))$ is not squares in $\C (x)$ because $\textrm{Gcd}(1-C(x^{2},B(x^{2})-C(x^{2}))=1$. 
\fin
\end{enumerate}

\end{subsubsection}

\end{subsection}

\begin{subsection}{Richelot's isogenies.}

\begin{sprop}\label{genre2-prop-base} 
Let $K$ be a characteristic $0$ field. 
Let $J$ and $\widehat{J}$ be two abelian varieties defined over $K$. We assume that $J(K)$ is finitely generated. 
We assume the existence of two isogenies $\varphi : J\longrightarrow\widehat{J}$ and $\widehat{\varphi}:\widehat{J}\longrightarrow J$ such that $\varphi\circ\widehat{\varphi}=[2]_{\widehat{J}}$ and $\widehat{\varphi}\circ\varphi =[2]_{J}$. 
Then the $K$-Mordell-Weil rank of $J$ is $0$ if and only if 
$J(K)/\widehat{\varphi}(\widehat{J}(K))=J(K)_{tors}/\widehat{\varphi}(\widehat{J}(K))$ and \linebreak
$\widehat{J}(K)/\varphi (J(K))=\widehat{J}(K)_{tors}/\varphi (J(K)).$
\end{sprop}
\dem%
Assume that $J(K)/\widehat{\varphi}(\widehat{J}(K))=J(K)_{tors}/\widehat{\varphi}(\widehat{J}(K))$ and \linebreak
$\widehat{J}(K)/\varphi (J(K))=\widehat{J}(K)_{tors}/\varphi (J(K)).$ The group $J(K)/J(K)_{tors}$ is free and finitely generated. Assume its rank is greater than $0$ and let $\alpha$ be an element of one of its $\Z$-bases. Notice that $\alpha$ is not a double in $J(K)/J(K)_{tors}$. Let $\beta$ be an element of the class $\alpha$. 

By hypothesis $\beta$ is equal to $T_{1}+\widehat{\varphi}(\beta_{1})$ for some torsion point $T_{1}\in J(K)$ and some point $\beta_{1}\in\widehat{J}(K)$. In the same way we get $\beta_{1}=T_{2}+\varphi (\beta_{2})$ for some $T_{2}\in\widehat{J}(K)_{tors}$ and some $\beta_{2}\in J(K)$. As a consequence $\beta$ is equal to $T_{1}+\widehat{\varphi}(T_{2})+2\beta_{2}$ and thus $\alpha$ is a double in $J(K)/J(K)_{tors}$. This is a contradiction with the choice of $\alpha$.
\fin%

\begin{snota}\label{nota-genre2-richelot-cas-ell-def-gamma}

\rm Let $K$ a characteristic $0$ field. Let $\alpha$ and $\beta$ be two elements of $K$ such that $\beta (\alpha^{2}-4\beta )\neq 0$. Let $\CbeD$ be the elliptic curve defined over $K$ by the affine equation $z^{2}=y(y^{2}+\alpha y+\beta )$. 
We define a group homomorphism \linebreak$\gamma_{\CbeD}:\CbeD (K)\longrightarrow K^{\times}/K^{\times 2}$ by mapping a $K$-point $(y,z)\in\CbeD (K)$ to 
\begin{itemize}
\item[* ]\rm the class of $y$ in $K^{\times}/K^{\times 2}$ when $y$ is nontrivial, 
\item[* ]\rm the class of $\beta$ in $K^{\times}/K^{\times 2}$ when $(y,z)=(0,0)$  
\item[* ]\rm and by mapping the neutral element of $\CbeD (K)$ to the class of $1$.
\end{itemize}
\end{snota}

\begin{sprop}\label{formulation-Pb-courbe-ell} 
We use the notation 
\ref{nota-genre2-richelot-cas-ell-def-gamma}. We assume that $\CbeD (K)$ is finitely generated. 
Let $\widehat{\CbeD}$ be the elliptic curved defined over $K$ by the affine equation $z^{2}=y(y^{2}-2\alpha y+\alpha^{2}-4\beta ).$ 

Then the Mordell-Weil rank of $\CbeD (K)$ is $0$ if and only if \linebreak$\gamma_{\CbeD}(\CbeD (K))=\gamma_{\CbeD}(\CbeD (K)_{tors})$ and $\gamma_{\widehat{\CbeD}}(\widehat{\CbeD}(K))=\gamma_{\widehat{\CbeD}}(\widehat{\CbeD}(K)_{tors})$. 
\end{sprop}
\dem%
There are two isogenies $\varphi : \CbeD\longrightarrow\widehat{\CbeD}$ and $\widehat{\varphi}:\widehat{\CbeD}\longrightarrow\CbeD$ such that $\varphi\circ\widehat{\varphi}=[2]_{\widehat{\CbeD}}$ and $\widehat{\varphi}\circ\varphi =[2]_{\CbeD}$, and such that $\textrm{Ker}(\gamma_{\CbeD})=\widehat{\varphi}(\widehat{\CbeD}(K))$ and $\textrm{Ker}(\gamma_{\widehat{\CbeD}})=\varphi (\CbeD (K))$ (for more details see \cite{Silverman-Tate} section III.4, page 76 and section III.5, page 85). Proposition \ref{formulation-Pb-courbe-ell} is proved by applying Proposition \ref{genre2-prop-base} to the isogenies $\varphi$ and $\widehat{\varphi}$.
\fin%

\begin{snota}\label{nota_genre2}\rm Let $K$ be a characteristic $0$ field. Let $G_{1}(y),G_{2}(y),G_{3}(y)$ be three elements of $K[y]$ of degree at most $2.$ We assume that the degree $G_{1}G_{2}G_{3}$ is $5$. Denote by $K_{i,\CbeH}$ the $K$-algebra $K[T]/(G_{i}(T))$. Let $\CbeH$ be the hyperelliptic curve defined over $K$ by the affine equation
$$\CbeH : z^{2}=G_{1}(y)G_{2}(y)G_{3}(y).$$

We define a group homomorphism $\Pi_{\CbeH}:\textrm{Jac}(\CbeH )(K)\longrightarrow (K^{\times}/K^{\times 2})^{3}$ in the following way: when $\textrm{div}(u,v)\in\textrm{Div}^{0}(K(\CbeH ))$ is a semi-reduced divisor with $u$ coprime to $G_{i}$ and if $\alpha$ is its linear equivalence class, 
then the i-th coordinate of $\Pi_{\CbeH}(\alpha)$ is the class of $N_{K_{i,\CbeH}/K}((-1)^{\deg (u)}u(T))$ in $K^{\times}/K^{\times 2}$. 
\end{snota}



\begin{sprop}\label{genre2-reformulation-Richelot} We use the notation \ref{nota_genre2}. We assume that $\textrm{Jac}(\CbeH )(K)$ is finitely generated. 
If $P$, $Q\in K[y]$ are two polynomials then we denote by $[P,Q]$ the polynomial $P'(y)Q(y)-P(y)Q'(y)$. Let $\widehat{\CbeC}$ be the hyperelliptic curve defined over $K$ by the affine equation 
$$\widehat{\CbeC} : \Delta \widehat{z}^{2}=L_{1}(\widehat{y})L_{2}(\widehat{y})L_{3}(\widehat{y}).$$ 
where $L_{1}:=[G_{2},G_{3}],$ $L_{2}:=[G_{3},G_{1}]$, $L_{3}:=[G_{1},G_{2}]$ and $\Delta :=\textrm{det}(g_{i,j})$ (we denote by $g_{i,j}$ the coefficient of $G_{i}$ relative to the monomial $y^{j}$). 

Then the Mordell-Weil rank of $\textrm{Jac}(\CbeH )(K)$ is $0$ if and only if \linebreak$\Pi_{\CbeH}(\textrm{Jac}(\CbeH )(K))$ is equal to $\Pi_{\CbeH}(\textrm{Jac}(\CbeH )(K)_{tors})$ and $\Pi_{\widehat{\CbeH}}(\textrm{Jac}(\widehat{\CbeH})(K))$ is equal to $\Pi_{\widehat{\CbeH}}(\textrm{Jac}(\widehat{\CbeH})(K)_{tors})$. 
\end{sprop}
\dem%
We consider Richelot's isogenies $\varphi :\textrm{Jac}(\CbeH )\longrightarrow\textrm{Jac}(\widehat{\CbeH})$ and $\widehat{\varphi}:\textrm{Jac}(\widehat{\CbeH})\longrightarrow\textrm{Jac}(\CbeH )$ (see \cite{Cassels-Flynn} chapter 9). 
Following \cite{Cassels-Flynn} Theorem 9.9.1, we have $\textrm{Ker}(\Pi_{\CbeH})=\widehat{\varphi}(\textrm{Jac}(\widehat{\CbeH})(K))$ and $\textrm{Ker}(\Pi_{\widehat{\CbeH}})=\varphi (\textrm{Jac}(\CbeH )(K))$. 
Proposition \ref{formulation-Pb-courbe-ell} is proved by applying Proposition \ref{genre2-prop-base} to the isogenies $\varphi$ and $\widehat{\varphi}$.
\fin%
\begin{stheo}\label{theo-critere-nullite-rang} 

We use the notation \ref{Notations_Principales_21}. 
We assume that the hypotheses of Corollary \ref{stheo-3-3-final-desc-finie} are satisfied. 
For each $\delta\in k(x)^{\times}$ we consider the hyperelliptic curves $\CbeC^{+}_{\delta},$ $\widehat{\CbeC}^{+}_{\delta},$ $\CbeC^{-}_{\delta},$ $\widehat{\CbeC}^{-}_{\delta}$ defined over $k(x)$ by the affine equations affines
$$\begin{array}{l}
\CbeC^{+}_{\delta}:z^{2}=\left( y+\frac{\delta(1+C(x))}{2}\right)\left( y^{2}-\left(\frac{\delta (1-C(x))}{2}\right)^{2}\right)
\left( y^{2}-\frac{\delta^{2}[(1+C(x))^{2}-4B(x)]}{4}\right)\\
\widehat{\CbeC}^{+}_{\delta}:z^{2}=(y+\delta (1+C(x)))(y^{2}-4\delta^{2}B(x))(y^{2}-4\delta^{2}C(x))\\
\CbeC^{-}_{\delta}:z^{2}=y(y^{2}-\delta [(1-C(x))^{2}-2(B(x)-C(x))]y+\delta^{2}(B(x)-C(x))^{2})\\
\widehat{\CbeC}^{-}_{\delta}:t^{2}=y(y+\delta (1-C(x))^{2})(y+\delta ((1-C(x))^{2}-4(B(x)-C(x)))).\\
\end{array}$$
\rm We use the notation \ref{nota-genre2-richelot-cas-ell-def-gamma} relative to $\CbeC^{-}_{\delta}$ and $\widehat{\CbeC}^{+}_{\delta}$, and the notation \ref {nota_genre2}
\begin{itemize}
\item[* ] with $G_{1}(y):=y+\frac{\delta(1+C(x))}{2}$ and $G_{2}(y):=y^{2}-\left(\frac{\delta (1-C(x))}{2}\right)^{2}$ and $G_{3}(y):=y^{2}-\frac{\delta^{2}[(1+C(x))^{2}-4B(x)]}{4}$ when $\CbeH=\CbeC^{+}_{\delta}$;
\item[* ] with $G_{1}(y):=y-\delta\left( 1+C\left( x\right)\right)$ and $G_{2}(y):=y^{2}-4\delta^{2}B\left( x\right)$ and $G_{3}(y):=y^{2}-4\delta^{2}C\left( x\right)$ when $\CbeH=\widehat{\CbeC}^{+}_{\delta}$.
\end{itemize}
%
%
Then the $\R (x)$-Mordell-Weil rank of $\textrm{Jac}(\CbeC )$ is $0$ if and only if for every positive element $\zeta\in k^{\times}$ the images of 
$$
\gamma_{\CbeC^{-}_{\zeta}},\textrm{ } \gamma_{\CbeC^{-}_{\zeta x}},\textrm{ } \gamma_{\widehat{\CbeC}^{-}_{\zeta}},\textrm{ } \gamma_{\widehat{\CbeC}^{-}_{\zeta x}},\textrm{ } \Pi_{\CbeC^{+}_{\zeta}},\textrm{ } \Pi_{\CbeC^{+}_{\zeta x}},\textrm{ } \Pi_{\widehat{\CbeC}^{+}_{\zeta}}\textrm{ and }\Pi_{\widehat{\CbeC}^{+}_{\zeta x}}
$$ 
are respectively the images of the $k(x)$-rational torsion subgroups of 
$$
\CbeC^{-}_{\zeta},\textrm{ } \CbeC^{-}_{\zeta x},\textrm{ } \widehat{\CbeC}^{-}_{\zeta},\textrm{ } \widehat{\CbeC}^{-}_{\zeta x},\textrm{ } \textrm{Jac}(\CbeC^{+}_{\zeta}),\textrm{ } \textrm{Jac}(\CbeC^{+}_{\zeta x}),\textrm{ } \textrm{Jac}(\widehat{\CbeC}^{+}_{\zeta})\textrm{ and }\textrm{Jac}(\widehat{\CbeC}^{+}_{\zeta x}).
$$ 
\end{stheo}
\dem%
Following Corollary \ref{stheo-3-3-final-desc-finie}, the $\R (x)$-Mordell-Weil rank of $\textrm{Jac}(\CbeC)$ is $0$ if and only if for every positive element $\zeta\in k^{\times}$ the $k(x)$-Mordell-Weil ranks of $\textrm{Jac}(\CbeC^{+}_{\zeta}),$ $\textrm{Jac}(\CbeC^{+}_{\zeta x}),$ $\CbeC^{-}_{\zeta}$ and $\CbeC^{-}_{\zeta x}$ are trivials (notice that even if the curve $\CbeC^{+}_{\delta}$ is not the same as in Corollary \ref{stheo-3-3-final-desc-finie}, both curves are isomorphic over $k(x)$: replace $y$ by $y-\frac{\delta(1+C(x))}{2}$ in the equation defining $\CbeC^{+}_{\delta}$). Theorem \ref{theo-critere-nullite-rang} is obtained by applying Proposition \ref{formulation-Pb-courbe-ell} to the elliptic curves $\CbeC^{-}_{\zeta}$ and $\CbeC^{-}_{\zeta x}$ and Proposition \ref{genre2-reformulation-Richelot} to the hyperelliptic curves $\CbeC^{+}_{\zeta}$ and $\CbeC^{+}_{\zeta x}$.
\fin%

\end{subsection}

\end{section}

\begin{section}{A familly of positive polynomials that are not a sum of three squares in $\R (x,y)$.}\label{chapitre-Rang-MW}

\begin{subsection}{How we will proceed.}
\begin{snota}\label{nota-section-4-3-carac-im-pi-H1}
\rm Let $k$ be a characteristic $0$ field. Let $f(y)\in k(x)[y]$ be a monic squarefree polynomial with odd degree and $\CbeH$ be the hyperelliptic curve defined on $k(x)$ by the affine equation $z^{2}=f(y).$  

\rm Let $f=\displaystyle\prod_{i=1}^{r}P_{i}(y)$ be the decomposition of $f$ into prime elements of $k(x)[y].$ Denote the algebra $k(x)[y]/f(y)$ by $M$, the field $k(x)[y]/P_{i}(y)$ by $K_{i}$ and the class of $y$ in $K_{i}$ by $y_{i}$. 

Let $\pi_{\CbeH}:\textrm{Jac}(\CbeH )(k(x))\longrightarrow M^{\times}/M^{\times 2}$ be the morphism obained by applying Proposition \ref{def-Cas-Sch+im-in-ker} to $\CbeH$ and $k(x)$. Using the chinese remainder theorem, we can think of $\pi_{\CbeH}$ as a morphism $\pi_{\CbeH}:\textrm{Jac}(\CbeH )(k(x))\longrightarrow\displaystyle\prod_{i=1}^{r}K_{i}^{\times}/K_{i}^{\times 2}$.
\end{snota}

\begin{snota}\label{nota-section-4-3-carac-im-pi-H2} 
\rm The norm map $N_{K_{i}/k(x)}$ asociated to the field extension $K_{i}/k(x)$ induces a homomorphism $N_{K_{i}/k(x)}:K_{i}^{\times}/K_{i}^{\times 2}\longrightarrow k(x)^{\times}/k(x)^{\times 2}.$ 

Denote by $\pi_{\CbeH ,i}:\textrm{Jac}(\CbeH )(k(x))\longrightarrow K_{i}^{\times}/K_{i}^{\times 2}$ the i-th coordinate of $\pi_{\CbeH}$ and by $\Xi_{\CbeH,i}$ the composite map $N_{K_{i}/k(x)}\circ\pi_{\CbeH ,i}.$ In this subsection we study the homorphism $\Xi_{\CbeH}:\textrm{Jac}(\CbeH )(k(x))\longrightarrow\displaystyle\prod_{i=1}^{r}k(x)^{\times}/k(x)^{\times 2}$ with i-th coordinate $\Xi_{\CbeH,i}$ (for ervery $i\in\{ 1,\cdots ,r\}$). 
\end{snota}

\begin{sprop}\label{caracterisationimagePi} We use the notation \ref{nota-section-4-3-carac-im-pi-H1} and \ref{nota-section-4-3-carac-im-pi-H2}. Let $i\in\{ 1,\cdots ,r\}$ be an interger. For each index $j\neq i,$ denote by $d_{i,j}$ the rational fraction
$$d_{i,j}:=\pgcd\left( N_{K_{i}/k(x)}\left( P_{i}'\left( y\right)\displaystyle\prod_{k\neq i}P_{k}\left( y\right)\right) ,N_{K_{j}/k(x)}\left( P_{j}'\left( y\right)\displaystyle\prod_{k\neq j}P_{k}\left( y\right)\right)\right) .$$
Then each element of the image of $\Xi_{\CbeH}$ is a class $\left(\left[\displaystyle\prod_{j\neq 1}\mu_{1,j}\right] ,\cdots ,\left[\displaystyle\prod_{j\neq r}\mu_{r,j}\right]\right)$ 
for some family $(\mu_{i,j})_{\begin{array}{l}1\le i\le r,\\j\neq i\\\end{array}}$ of squarefree elements of $k[x]$ such that $\mu_{i,j}=\mu_{j,i}$ and such that 
the prime factors of $\mu_{i,j}$ are prime factors of  $d_{i,j}$. %
\end{sprop}
\dem %
Let $\beta$ be an element of $\textrm{Jac}(\CbeH )(k(x))$ and let $\textrm{div}(u,v)$ be a semi-reduced divisor on $\CbeH$ with linear equivalence class $\beta$. 
The field $K_{i}$ is the fraction field of a Dedekind ring $\mathcal{O}_{i}$. 
Let 
\begin{equation}\label{eqt-lemme-premiere-etude-image-Xi-411}
\left( (-1)^{\deg (u)}u(y_{i})\right) =\displaystyle\prod_{l\in I}\mathcal{P}_{i,l}^{2n_{i,l}+1}\times\displaystyle\prod_{l\in\widetilde{I}}\mathcal{Q}_{i,l}^{2m_{i,l}}
\end{equation}
be the decomposition of $\left( (-1)^{\deg (u)}u(y_{i})\right)$ into prime ideals of $\mathcal{O}_{i}$ 
(i.e. $(\mathcal{P}_{i,l})_{l\in I}$ and $(\mathcal{Q}_{i,l})_{l\in\widetilde{I}}$ are two families of prime ideals of $\mathcal{O}_{i}$). 

Applying $N_{K_{i}/k(x)}$ to Decomposition \ref{eqt-lemme-premiere-etude-image-Xi-411} we get the existence of $\alpha_{i}\in k[x]$ squarefree whose prime factors are prime factors of $\displaystyle\prod_{l\in I}N_{K_{i}/k(x)}(\mathcal{P}_{i,l})$ and the existence of $\gamma_{i}\in k(x)^{\times}$ such that 
$N_{K_{i}/k(x)}((-1)^{\deg (u)}u(y_{i}))=\alpha_{i}\gamma_{i}^{2}$. From the definitions of $\pi_{\CbeH}$ and $\Xi_{\CbeH}$ 
we know that $\Xi_{\CbeH,i}(\beta )$ 
is the class of $\alpha_{i}$ in $k(x)^{\times}/k(x)^{\times 2}$. 

\begin{description}
\item[Step $1$: applying Proposition \ref{caract-image-CS-v-ideaux}. ]  
Following Proposition \ref{caract-image-CS-v-ideaux}, all the valuation $v_{\mathcal{P}_{i,l}}(f'(y_{i}))$ are nontrivial. Thus, for each $l\in I$, the prime factors of $N_{K_{i}/k(x)}(\mathcal{P}_{i,l})$ are prime factors of $N_{K_{i}/k(x)}(f'(y_{i}))=N_{K_{i}/k(x)}\left( P_{i}'\left( y_{i}\right)\displaystyle\prod_{j\neq i}P_{j}\left( y_{i}\right)\right)$. In particular the prime factors of $\alpha_{i}$ are prime factors of $N_{K_{i}/k(x)}\left( P_{i}'\left( y_{i}\right)\displaystyle\prod_{j\neq i}P_{j}\left( y_{i}\right)\right)$.
\item[Step $2$: applying Proposition \ref{def-Cas-Sch+im-in-ker}. ] For $\mu_{i,j}$ we take $\pgcd (\alpha_{i},\alpha_{j})$; 
the leading coefficient of $\mu_{i,j}$ can be choosen such that $\alpha_{i}$ and $\displaystyle\prod_{j\neq i}\mu_{i,j}$ have the same leading coefficient. Write $\displaystyle\prod_{j\neq i}\mu_{i,j}=\Lambda_{i}\Gamma_{i}^{2}$ 
with $\Lambda_{i}\in k[x]$ squarefree and $\Gamma_{i}\in k[x]$ monic. 
\begin{enumerate}
\item If $p$ divides $\alpha_{i}$, then $v_{p}(\mu_{i,j})=v_{p}(\alpha_{j})$ for each $j\neq i$ (notice that $\alpha_{j}$ and $\mu_{i,j}$ are squarefree). In particular if $p$ divides $\alpha_{i}$ then $v_{p}(\displaystyle\prod_{j\neq i}\mu_{i,j})=v_{p}(\displaystyle\prod_{j\neq i}\alpha_{i})$.
\item Following Proposition \ref{def-Cas-Sch+im-in-ker} the product $\displaystyle\prod_{i=1}^{r}\Xi_{\CbeH ,i}(\beta )$ is the trivial element of $k(x)^{\times}/k(x)^{\times 2}$. Since $\Xi_{\CbeH ,i}(\beta )$ is the class of $\alpha_{i}$ in $k(x)^{\times}/k(x)^{\times 2}$, the polynomial $\displaystyle\prod_{i=1}^{r}\alpha_{i}$ is a square in $k(x)$. 
In particular, for every prime $p$, the valuations $v_{p}(\displaystyle\prod_{j\neq i}\alpha_{j})$ and $v_{p}(\alpha_{i})$ are congruent modulo $2$.
\end{enumerate}
From those two remarks we deduce that for every prime divisor $p$ of $\alpha_{i}$ the valuation $v_{p}(\displaystyle\prod_{j\neq i}\mu_{i,j})$ is odd. Since $\alpha_{i}$ is squarefree this implies the divisibility of $\Lambda_{i}$ by $\alpha_{i}$. 

Moreover the definition of $\Lambda_{i}$ implies that $\Lambda_{i}$ divides $\alpha_{i}$. Since they have the same leading coefficient, $\Lambda_{i}$ and $\alpha_{i}$ are equal. Thus $\Xi_{\CbeH ,i}(\beta )$ is the class of $\displaystyle\prod_{j\neq i}\mu_{i,j}$ in $k(x)^{\times}/k(x)^{\times 2}$. 
From Step 1 we know that the prime factors of $\mu_{i,j}=\gcd (\alpha_{i},\alpha_{j})$ are prime factors of $d_{i,j}$.
\fin%


\end{description}

\begin{slemm}\label{lem-equation-pol-cbe-ell}
Let $k$ be a characteristic $0$ field. Let $S,$ $T\in k[x]$ be such that $T(S^{2}-4T)$ is nontrivial. Denote by $\CbeD$ the elliptic curve defined over $k(x)$ by the Weierstrass equation
$$\CbeD : \beta^{2}=\alpha (\alpha^{2}+S\alpha +T).$$
Let $(\alpha ,\beta )\neq (0,0)$ be a $k(x)$-point of $\CbeD$. 
Then there is a constant $\epsilon\in k$, a monic squarefree divisor $\mu \in k[x]$ of $T$ and two polynomials $\theta\in k[x]$ and $\psi\in k[x]$ such that %
$\alpha =\epsilon\mu\frac{\theta^{2}}{\psi^{2}}$ and $\mu\theta$ is coprime to $\psi $ and %
$$\epsilon\mu\nu^{2}=\epsilon^{2}\mu^{2}\theta^{4}+S\epsilon\mu\theta^{2}\psi^{2}+T\psi^{4}.$$ 
\end{slemm}


\begin{snota}\label{nota-equiv-mod-carre}
\rm Two rational fractions $\alpha\in k(x)^{\times}$ and $\beta\in k(x)^{\times}$ are said to be equivalent modulo squares (and we write $\alpha\sim\beta$) if there exists $\gamma\in k(x)^{\times}$ such that $\alpha =\gamma^{2}\beta$. The equivalence class of $a\in k(x)^{\times}$ relative to the relation $\sim$ is denoted by $[a]$. %
%
\end{snota}

\begin{snota}\label{nota-etude_y2-A_en_p}
\rm Let $\mathcal{P}$ be a place of $k(x)$ and $\mathcal{O}_{\mathcal{P}}$ be the associated valuation ring. Let $\alpha ,$ $\beta$ be elements of $\mathcal{O}_{\mathcal{P}}^{\times}.$ %

\rm The element $\alpha$ is equivalent to $\beta$ modulo $\mathcal{P}$ and modulo squares (and we write $\alpha\sim\beta\bmod\mathcal{P}$) if there exists $\gamma\in\mathcal{O}_{\mathcal{P}}^{\times}$ such that $\alpha$ are $\beta\gamma^{2}$ are congruent modulo $\mathcal{P}$.
\end{snota}


The two following propositions are a first example of the specialisation results we use in the four next sections to study of the image of $\Xi_{\CbeH ,i}=N_{K_{i}/k(x)}\circ\pi_{\CbeH ,i}$. They are stated using common notation.

\begin{snota}\label{etude_y2-A_en_p-enonce-nota-comm}
Let $k$ be a characteristic $0$ field. Let $A$ be an element of $ k(x)$. Denote $k(x)[T]/(T^{2}-A)$ by $K$ and the class of $T$ in $K$ by $t$.%

Let $\mathcal{P}$ be a place of $k(x),$ $\mathcal{O}_{\mathcal{P}}$ be the associated valuation ring, $v_{\mathcal{P}}$ be the valuation at $\mathcal{P}$, and $p$ be an uniformizor of $\mathcal{P}.$
\end{snota}

\begin{sprop}\label{etude_y2-A_en_p}  We use the notation \ref{nota-etude_y2-A_en_p} and \ref{etude_y2-A_en_p-enonce-nota-comm}. We denote the element $p^{-v_{\mathcal{P}}(A)}A\in\mathcal{O}_{\mathcal{P}}^{\times}$ by $\widetilde{A}$. Let $u:=u_{0}(y^{2})+yu_{1}(y^{2})\in k(x)[y]$ be a polynomial and denote $p^{-v_{\mathcal{P}}(N_{K/k(x)}(u(t)))}N_{K/k(x)}(u(t))\in\mathcal{O}_{\mathcal{P}}^{\times}$ by $\alpha$. Assume that $v_{\mathcal{P}}(A)$ is odd.
\begin{enumerate}
\item If $v_{\mathcal{P}}(N_{K/k(x)}(u(t)))$ is even, then $\alpha\sim 1\bmod\mathcal{P}$ ;
\item if $v_{\mathcal{P}}(N_{K/k(x)}(u(t)))$ is odd, then $\alpha\sim -\widetilde{A}\bmod\mathcal{P}$ and 
$$\frac{v_{\mathcal{P}}(A)+1}{2}+v_{\mathcal{P}}(u_{1}(A))\le v_{\mathcal{P}}(u_{0}(A)).$$
\end{enumerate}
\end{sprop}
\dem%
The valuation $v_{\mathcal{P}}((u_{0}(A))^{2})$ is even and the valuation $v_{\mathcal{P}}(A(u_{1}(A))^{2})$ is odd. These valuations are differents. As a consequence, the valuation at $\mathcal{P}$ of $N_{K/k(x)}(u(t))=(u_{0}(A))^{2}-A(u_{1}(A))^{2}$ is $\textrm{Min}(v_{\mathcal{P}}((u_{0}(A))^{2}),v_{\mathcal{P}}(A(u_{1}(A))^{2})).$

\begin{description}
\item[If $v_{\mathcal{P}}(N_{K/k(x)}(u(t)))$ is even.] Then $v_{\mathcal{P}}(N_{K/k(x)}(u(t)))$ is equal to \linebreak$v_{\mathcal{P}}((u_{0}(A))^{2})$ and is strictly less than $v_{\mathcal{P}}(A(u_{1}(A))^{2})$. 
The class of
$$\alpha =p^{-v_{\mathcal{P}}(N_{K/k(x)}(u(t)))}(u_{0}(A))^{2}-p^{-v_{\mathcal{P}}(N_{K/k(x)}(u(t)))}A(u_{1}(A))^{2}$$ 
in the residual field at $\mathcal{P}$ is nontrivial and equal to the class of %
\linebreak$\left( p^{-v_{\mathcal{P}}(u_{0}(A))}u_{0}(A)\right)^{2}.$

\item[If $v_{\mathcal{P}}(N_{K/k(x)}(u(t)))$ is odd.] Then $v_{\mathcal{P}}(N_{K/k(x)}(u(t)))$ is equal to \linebreak$v_{\mathcal{P}}(A(u_{1}(A))^{2})$ and is strictly less than $v_{\mathcal{P}}((u_{0}(A))^{2})$. 
The class of %
$$\alpha =p^{-v_{\mathcal{P}}(N_{K/k(x)}(u(t)))}(u_{0}(A))^{2}-p^{-v_{\mathcal{P}}(N_{K/k(x)}(u(t)))}A(u_{1}(A))^{2}$$ 
in the residual field at $\mathcal{P}$ is nontrivial and equal to the class of %
\linebreak$-\widetilde{A}\left( p^{-v_{\mathcal{P}}(u_{1}(A))}u_{1}(A)\right)^{2}.$ Moreover, since $v_{\mathcal{P}}((u_{0}(A))^{2})$ is greater than or equal to $v_{\mathcal{P}}(A(u_{1}(A))^{2})+1$, the stated inequality holds. 
\fin%
\end{description}

\begin{sprop}\label{etude_y2-A_en_p-cas-pair} 

We use the notation \ref{nota-etude_y2-A_en_p} and and \ref{etude_y2-A_en_p-enonce-nota-comm}. We assume that
\begin{itemize}
\item[* ] the valuation $v_{\mathcal{P}}(A)$ is even,
\item[* ] that $p^{-v_{\mathcal{P}}(A)}A$ is not a square in the residual field $\mathcal{O}_{\mathcal{P}}/\mathcal{P}$. 
\end{itemize}
Then, for every polynomial $u:=u_{0}(y^{2})+yu_{1}(y^{2})\in k(x)[y]$, the valuation $v_{\mathcal{P}}(N_{K/k(x)}(u(t)))$ is even.
\end{sprop}
\dem%
Assume that the valuation $v_{\mathcal{P}}(N_{K/k(x)}(u(t)))$ is odd. Denote by $r$ the minimum $\textrm{Min}\left( v_{\mathcal{P}}(u_{0}(A)),\frac{v_{\mathcal{P}}(A)}{2}+v_{\mathcal{P}}(u_{1}(A))\right) .$ When the classes of two elements $\alpha$ and $\beta$ in $\mathcal{O}_{\mathcal{P}}$ have the same classes in the residual field  $\mathcal{O}_{\mathcal{P}}/\mathcal{P}$ we write $\alpha\equiv\beta\bmod\mathcal{P}$.

From the definition of $r$ we know that the valuation at $\mathcal{P}$ of 
$$p^{-2r}N_{K/k(x)}(u(t))=\left( p^{-r}u_{0}(A)\right)^{2}-p^{-v_{\mathcal{P}}(A)}A\left( p^{(v_{\mathcal{P}}(A)/2-r)}u_{1}(A)\right)^{2}$$ 
is nonnegative. We have more : since it is odd, this valuation is positive. In particular we have  %
\begin{equation}\label{congruence-etude_y2-A_en_p-cas-pair}
(p^{-r}u_{0}(A))^{2}\equiv p^{-v_{\mathcal{P}}(A)}A\left( p^{(v_{\mathcal{P}}(A)/2-r)}u_{1}(A)\right)^{2}\bmod\mathcal{P}.
\end{equation} 

From the definition of $r$ and from the congruence \ref{congruence-etude_y2-A_en_p-cas-pair} we deduce that $v_{\mathcal{P}}\left( p^{(v_{\mathcal{P}}(A)/2-r)}u_{1}(A)\right) =v_{\mathcal{P}}\left( u_{0}(A)\right) =0$. 
In particular, the class of \linebreak$p^{(v_{\mathcal{P}}(A)/2-r)}u_{1}(A)$ in $\mathcal{O}_{\mathcal{P}}/\mathcal{P}$ is invertible. Thus the congruence \ref{congruence-etude_y2-A_en_p-cas-pair} shows that $p^{-v_{\mathcal{P}}(A)}A$ is a square in $\mathcal{O}_{\mathcal{P}}/\mathcal{P}$. This is a contradiction. 
\fin%
\end{subsection}

\begin{subsection}{A study of the elliptic curve $\CbeC^{-}_{\delta}.$}\label{sous-section-E}
\begin{snota}\label{nota-section42} We use the notation \ref{Notations_Principales_21}. For every $\delta\in k(x)^{\times}$ let $\CbeC^{-}_{\delta}$ be the elliptic curve defined on $k(x)$ by the affine equation :
$$\CbeC^{-}_{\delta}:t^{2}=s\left( s^{2}-\delta\left(\left( 1-C\right)^{2}-2\left( B-C\right)\right) s+\delta^{2}\left( B-C\right)^{2}\right) .$$
In our statements and their proofs we also use the notation \ref{nota-genre2-richelot-cas-ell-def-gamma} and \ref{nota-equiv-mod-carre}.
\end{snota}
\rmk The polynomials $B-C,$ $1-C,$ $(1-C)^{2}-4(B-C)$ and $\delta$ are nontrivial. 
Thus $s\left( s^{2}-\delta\left(\left( 1-C\right)^{2}-2\left( B-C\right)\right) s+\delta^{2}\left( B-C\right)^{2}\right)$ is squarefree.

\begin{snota}\label{nota-section42bis}
To simplify our statements we denote $1-C$ by $e$ and $B-C$ by $d$. Notice that $e$ has degree $1$, $d$ is monic of even degree and that $e^{2}-4d$ has odd degree. 

Let $(\alpha ,\beta)$ be a $k(x)$-point of $\CbeC^{-}_{\delta}$. Following proposition \ref{lem-equation-pol-cbe-ell}, there is a constant $\epsilon\in k$, a monic squarefree divisor $\mu \in k[x]$ of $\delta d$ and two polynomials $\theta\in k[x]$ et $\psi\in k[x]$ such that $\alpha =\epsilon\mu\frac{\theta^{2}}{\psi^{2}}$ and $\mu\theta$ is coprime to $\psi $ and 
\begin{equation}\label{eq-lemm-42-C--delta}
\begin{array}{rcl}
\epsilon\mu\nu^{2}&=&\epsilon^{2}\mu^{2}\theta^{4}-\delta(e^{2}-2d)\epsilon\mu\theta^{2}\psi^{2}+\delta^{2}d^{2}\psi^{4}\\
&=&(\epsilon\mu\theta^{2}-\frac{\delta(e^{2}-2d)}{2}\psi^{2})^{2}-\frac{\delta^{2}e^{2}(e^{2}-4d)}{4}\psi^{4}.
\end{array}
\end{equation}
For every positive element $\zeta\in k^{\times}$ we study the image of $\gamma_{\CbeC^{-}_{\zeta}}$ and $\gamma_{\CbeC^{-}_{\zeta x}}$ (see the notation \ref{nota-genre2-richelot-cas-ell-def-gamma}). 
\end{snota}
\begin{slemm}\label{Lemma-C--delta42} 
%
We use the notation \ref{nota-section42}. Let $\delta\in k[x]$ be a nontrivial polynomial. We assume that the leading coefficient of $\delta$ and the leading coefficient of $(1-C)^{2}-4(B-C)$ are strictly greater than $0$.

Then every element of the image of $\gamma_{\CbeC^{-}_{\delta}}$ is the class in $k(x)^{\times}/k(x)^{\times 2}$ of a monic squarefree divisor $\mu$ of $\delta (B-C)$ with even degree. 


Moreover, for every divisor $p\in k[x]$ of $1-C$ coprime to $\delta (B-C)$, we have either $\mu\sim 1\bmod p$ or $\mu\sim -\delta d\bmod p$.
\end{slemm}
\dem%
We use the notation \ref{nota-section42bis}. Let $\lambda$ be the leading coefficient of $e^{2}-4d$. Applying Proposition \ref{etude_y2-A_en_p} with $\mathcal{P}$ the place above the infinite place of $k(x)$, $A:=\frac{\delta^{2}e^{2}(e^{2}-4d)}{4}$ and $u(y)=y-\alpha +\frac{\delta (e^{2}-2d)}{2}$ we get the class in $k^{\times}/k^{\times 2}$ of $\epsilon$: 
$$\epsilon\sim 1\textrm{ if }\deg(\mu)\textrm{ is even, and }\epsilon\sim -\lambda\textrm{ if }\deg(\mu)\textrm{ is odd}$$
(notice that following Proposition \ref{def-Cas-Sch+im-in-ker}, and if $K_{A}$ denotes the algebra $k(x)[t]/(t^{2}-A)$, the element $\alpha N_{K_{A}/k(x)}(-u(t))$ is a square in $k(x)$; thus $\epsilon$ has the same class in $k^{\times}/k^{\times 2}$ as the leading coefficient of $N_{K_{A}/k(x)}(-u(t))$).

Assume that $\deg(\mu)$ is odd. 
From our application of Proposition \ref{etude_y2-A_en_p} we know that:
\begin{equation}\label{etude-C-delta41-deg-infini}
\deg \left(\alpha -\frac{\delta (e^{2}-2d)}{2}\right)\le\frac{-1+\deg (\delta^{2}e^{2}(e^{2}-4d))}{2}. %
\end{equation}

As the degree of $e^{2}-4d$ is odd and the degree of $d$ is even, we have $\deg (e^{2}-2d)=\deg (e^2)>\deg (e^{2}-4d)$. Using those inequalities one can show
\begin{equation}\label{etude-C-delta41-deg-infini-contradictoire}
2\deg\left(\delta\left( e^{2}-2d\right)\right)>\deg (\delta^{2}e^{2}(e^{2}-4d)).
\end{equation}

The inequalities \ref{etude-C-delta41-deg-infini} and \ref{etude-C-delta41-deg-infini-contradictoire} are compatible only when $\alpha\psi^{2}=\epsilon\mu\theta^{2}$ and $\frac{\delta (e^{2}-2d)}{2}\psi^{2}$ have the same leading term. 

The leading coefficient of $e^{2}-2d=e^2-4d+2d$ is $2$ and $\mu$ is monic. Hence we know that $\epsilon\sim\zeta$ (with $\zeta$ the leading coefficient of $\delta$). In particular $\epsilon$ is positive. This is in contradiction with $\epsilon\sim -\lambda$ (because $\lambda$ is positive). Thus the degree of $\mu$ is even. \\

Let $p\in k[x]$ be a divisor of $e$. Assume that $p$ is coprime to $\delta d$. From Equation \ref{eq-lemm-42-C--delta} we deduce a congruence $\mu\nu^{2}\equiv (\mu\theta^{2}+\delta d\psi^{2})^{2}\bmod p$. If $p$ does not divide $\nu$, then $\mu\sim 1\bmod p$. In the other case, we use the divisibility of $\mu\theta^{2}+\delta d\psi^{2}$ by $p$, the coprimality of $\theta$ and $\psi$, and the coprimality of $\delta d$ and $p$ to show that $\mu\sim -\delta d\bmod p$.
\fin

\begin{slemm}\label{Lemma-C--delta42-supp} We use the notation \ref{nota-section42}. Let $\delta\in k[x]$ be a nontrivial polynomial. Let $p$ be a divisor of $\delta$ such that $v_{p}(\delta )=1$ and $v_{p}(1-C)=v_{p}(B-C)=0$. We assume that the degree $(1-C)^{2}-4(B-C)$ is $1$. 

Then every element of the image of $\gamma_{\CbeC^{-}_{\delta}}$ is the class in $k(x)^{\times}/k(x)^{\times 2}$ of a squarefree polynomial $\mu$ such that 
\begin{enumerate}
\item $\mu\sim 1\bmod p$ whenever $v_{p}(\mu)=0$ and
\item $p^{-1}\mu\sim p^{-1}\delta\bmod p$ whenever $v_{p}(\mu )=1$.
\end{enumerate}
\end{slemm}
\dem%
We use the notation \ref{nota-section42bis}. Assume that the valuation $v_{p}(\mu )$ is trivial.
\begin{description}
  \item[If $p$ is coprime to $\theta$. ]Then Equation \ref{eq-lemm-42-C--delta} gives $\epsilon\mu\nu^{2}\equiv\epsilon^{2}\mu^{2}\theta^{4}\bmod p$. 
  \item[If $p$ divides $\theta$: ] Then the polynomial 
$$p^{3}\epsilon\mu\left(\frac{\theta}{p}\right)^{2}\left(\epsilon\mu p\left(\frac{\theta}{p}\right)^{2}-\frac{\delta}{p}(e^{2}-2d)\psi^{2}\right) =\epsilon\mu\nu^{2}-\delta^{2}d^{2}\psi^{4}$$ 
is divisible by $p^{3}$. 
As $\psi$ and $\theta$ are coprime and $v_{p}(\delta )=1$ and $v_{p}(d)=0$, the valuation of $v_{p}(\delta^{2}d^{2}\psi^{4})$ at is $2$. Thus $v_{p}(\nu )$ is equal to $1$. Let $\widetilde{\nu}$ denote $p^{-1}\nu$ and $\widetilde{\delta}$ denote $p^{-1}\delta$. Equation \ref{eq-lemm-42-C--delta} implies a congruence $\epsilon\mu\widetilde{\nu}^{2}\equiv\widetilde{\delta}^{2}d^{2}\psi^{4}\bmod p$. 
\end{description}
In both cases we have an equivalence $\epsilon\mu\sim 1\bmod x$.

Now assume that the valuation $v_{p}(\mu )$ is $1$. Then, from Equation \ref{eq-lemm-42-C--delta}, the prime $p$ divides $\nu$. Denote the polynomial $p^{-1}\epsilon\mu$ by $\widetilde{\mu}$, the polynomial $p^{-1}\delta$ by $\widetilde{\delta}$ and the polynomial $p^{-1}\nu$ by $\widetilde{\nu}$. Simplifying Equation \ref{eq-lemm-42-C--delta} by $p$ gives an equation $p\widetilde{\mu}\widetilde{\nu}^{2}=\left(\widetilde{\mu}\theta^{2}+\widetilde{\delta}d\psi^{2}\right)^{2}-\widetilde{\delta}e^{2}\widetilde{\mu}\theta^{2}\psi^{2}.$
In particular we get 
\begin{equation}\label{Etud-C-del-41-congruence-mod-delta-cas-impair}
\widetilde{\delta}e^{2}\widetilde{\mu}\theta^{2}\psi^{2}\equiv\left(\widetilde{\mu}\theta^{2}+\widetilde{\delta}d\psi^{2}\right)^{2}\bmod p.
\end{equation}
If $p$ divides $\theta$ then the congruence \ref{Etud-C-del-41-congruence-mod-delta-cas-impair} shows that $p$ divides $\widetilde{\delta}d\psi^{2}$ i.e. that $p$ divides $\psi$. This would contradicts the coprimality of $\theta$ and $\psi$. Thus $p$ is coprime to $\theta$. In the same way $p$ is coprime to $\psi$. Hence the congruence \ref{Etud-C-del-41-congruence-mod-delta-cas-impair} gives an equivalence $\widetilde{\mu}\sim\widetilde{\delta}\bmod p$.
\fin%




\begin{sprop}\label{prop-verif-C-zeta41final}We use the notation \ref{nota-section42}. Assume the inequality \linebreak$\omega^{2}>\eta^{2}$ and the hypothesis \ref{hypothesisA} (see Assumptions \ref{Asummption24}). 

Then for every positive element $\zeta\in k^{\times}$ the image of $\gamma_{\CbeC^{-}_{\zeta}}$ is trivial.
\end{sprop}
\dem%
Let $\zeta >0$ be an element of $k$. 
Let $(\alpha ,\beta)$ be a $k(x)$-point of the curve $\CbeC^{-}_{\zeta}$. The leading coefficient of $(1-C)^{2}-4(B-C)$ is $(\omega^{2}-\eta^{2})$. This coefficient is positive by assumption. 
Following Lemma \ref{Lemma-C--delta42}, $\gamma_{\CbeC^{-}_{\zeta}}(\alpha ,\beta )$ is either the trivial class in $k(x)^{\times}/k(x)^{\times 2}$ or the class of $B-C$. 

Assume that $\gamma_{\CbeC^{-}_{\zeta}}(\alpha ,\beta )$ is the class of $B-C$. The remainder of the Euclidean division of $B-C$ by $1-C$ is not a square in $k$ (see hypothesis \ref{hypothesisA}). In particular this remainder is nontrivial and the last statement of Lemma \ref{Lemma-C--delta42} applies : 
$$B-C\sim 1\bmod (1-C)\textrm{ or }B-C\sim -\zeta(B-C)\bmod (1-C).$$
The first equivalence is in contradiction with the hypothesis \ref{hypothesisA}; the second one contradicts the positivity of $\zeta$.
\fin%

\begin{sprop}\label{prop-verif-C-zetax41final}
We use the notation \ref{nota-section42}. Assume the inequality \linebreak$\omega^{2}-\eta^{2}>2|\omega|$, and the hypotheses \ref{hypothesisB}, \ref{hypothesisC}, \ref{hypothesisDsupplementaire}, \ref{hypothesisEsupplementaire} and \ref{hypothesisI} (at least with $n_{1}=1$ and $n_{2}=n_{3}=0$) (see Assumptions \ref{Asummption24}).

Then for every positive element $\zeta\in k^{\times}$ the image of $\gamma_{\CbeC^{-}_{\zeta x}}$ is trivial.
\end{sprop}
\dem%
Let $\zeta >0$ be an element of $k$. 
Let $\beta$ be an element of the image of $\gamma_{\CbeC^{-}_{\zeta x}}$. The leading coefficient of $(1-C)^{2}-4(B-C)$ is $(\omega^{2}-\eta^{2})$. It is positive by assumption. 
Following Lemma \ref{Lemma-C--delta42}: the class $\beta$ is the class of 
a monic squarefree divisor $\mu$ of $x(B-C)$ with even degree 
\begin{description}
  \item[Case 1: $\mu\sim 1$. ] Then $\gamma_{\CbeC^{-}_{\zeta x}}(\alpha ,\beta )$ is the trivial class.
  \item[Case 2: $\mu\sim B-C$. ] Following Lemma \ref{Lemma-C--delta42-supp} we have $B-C\sim 1\bmod x$ (to show that the asumptions of Lemma \ref{Lemma-C--delta42-supp} are satisfied, use hypotheses \ref{hypothesisB}, \ref{hypothesisDsupplementaire} and \ref{hypothesisEsupplementaire}). This is in contradiction with the hypothesis \ref{hypothesisI} (applied with $n_{1}=1$ and $n_{2}=n_{3}=0$).
  \item[Case 3: $\mu\sim x(x+b_{1}-1-\omega )$. ] The remainder of the Euclidean division of $x^{2}(B-C)$ by $1-C$ is not a square in $k$ (see hypothesis \ref{hypothesisB} and hypothesis \ref{hypothesisC}). In particular this remainder is nontrivial and the last statement of Lemma \ref{Lemma-C--delta42} applies : 
$$\mu\sim 1\bmod (1-C)\textrm{ or }\mu\sim -\zeta x(B-C)\bmod (1-C).$$
After computing the remainders of the Euclidean division of $\mu$ and $-x(B-C)\mu$ by $1-C$ we obtain
\begin{equation}\label{alt-cas3-C-zx-61}
(2b_{1}-2+\omega^{2}-\eta^{2})(\omega^{2}-\eta^{2}+2\omega )\sim 1\textrm{ or }\zeta\sim 2(\omega^{2}-\eta^{2}-2\omega ).
\end{equation}
The first equivalence is in contradiction with the hypothesis \ref{hypothesisC}. Following Lemma \ref{Lemma-C--delta42-supp} (its asumptions are satisfied: use hypotheses \ref{hypothesisB}, \ref{hypothesisDsupplementaire} and \ref{hypothesisEsupplementaire}) we have $(b_{1}-1-\omega)\sim \zeta$. Thus if the second of the equivalences \ref{alt-cas3-C-zx-61} holds then 
we get a contradiction with hypothesis \ref{hypothesisDsupplementaire}.
  \item[Case 4: $\mu\sim x(x+b_{1}-1+\omega )$. ] This case is analogous to Case $3$ (with $-\omega$ instead of $\omega$).\fin%
\end{description}

\end{subsection}

\begin{subsection}{A study of the elliptic curve $\widehat{\CbeC}^{-}_{\delta}.$}\label{Etude-gamma-Cchap--delta43}
\begin{snota}\label{nota-section43} We use the notation \ref{Notations_Principales_21}. 
For every $\delta\in k(x)^{\times}$ let $\widehat{\CbeC}^{-}_{\delta}$ be the elliptic curve defined on $k(x)$ by the affine equation :
$$\widehat{\CbeC^{-}_{\delta}}:t^{2}=s(s+\delta (1-C)^{2})(s+\delta ((1-C)^{2}-4(B-C))).$$
In our statements and their proofs we also use the notation \ref{nota-genre2-richelot-cas-ell-def-gamma} and \ref{nota-equiv-mod-carre}.
\end{snota}
\rmk%
The polynomials $B-C,$ $1-C,$ $\delta$ and $(1-C)^{2}-4(B-C)$ are nontrivial. Thus the polynomial $s(s+\delta (1-C)^{2})(s+\delta ((1-C)^{2}-4(B-C)))$ is squarefree.
\begin{snota}\label{nota-section43bis}
Let $(\alpha ,\beta)$ be a $k(x)$-point of $\widehat{\CbeC}^{-}_{\delta}$. To simplify our statements we denote $1-C$ by $e$ and $B-C$ by $d$.

Following Proposition \ref{lem-equation-pol-cbe-ell}, there is a squarefree divisor $\mu \in k[x]$ of $\delta (1-C)((1-C)^{2}-4(B-C))$ and two polynomials $\theta\in k[x]$ et $\psi\in k[x]$ such that $\alpha =\mu\frac{\theta^{2}}{\psi^{2}}$, and $\mu\theta$ is coprime to $\psi $, and 
\begin{equation}\label{eq-lemm-43-Cchap--delta}
\mu\nu^{2}=(\mu\theta^{2}+\delta e^{2}\psi^{2})(\mu\theta^{2}+\delta (e^{2}-4d)\psi^{2}).
\end{equation}
\end{snota}

\begin{slemm}\label{ell_Echap1} We use the notation \ref{nota-section43}. Let $\delta\in k[x]$ be a nontrivial polynomial.
%
%
We assume that $-\delta (B-C)$ is not a square modulo $(1-C)$. 

Then every element of the image of $\gamma_{\widehat{\CbeC}^{-}_{\delta}}$ is the class in $k(x)^{\times}/k(x)^{\times 2}$ of a squarefree divisor of $\delta ((1-C)^{2}-4(B-C))$.
\end{slemm}
\dem%
We use the notation \ref{nota-section43bis}. The degree of $e=1-C$ is $1$. Assume the existence of a polynomial $\mu_{1}$ such that $\mu=e\mu_{1}=(1-C)\mu_{1}$. Then $\mu_{1}$ is a squarefree divisor of $\delta (e^{2}-4d)$. Denote by $\mu_{2}$ the unique polynomial such that $\delta (e^{2}-4d)=\mu_{1}\mu_{2}$. As $-\delta d$ is not a square modulo $e$, the polynomials $e$ (which is prime) and $\delta (e^{2}-4d)$ are coprime. Hence $\mu_{2}$ and $e$ are coprime. After simplification, Equation \ref{eq-lemm-43-Cchap--delta} becomes
\begin{equation}\label{eqn_Echap+2}
\nu^{2}=(\mu_{1}\theta^{2}+\delta e\psi^{2})(e\theta^{2}+\mu_{2}\psi^{2}).
\end{equation}
From this equation, we deduce a congruence $\nu^{2}\equiv -4\delta d\theta^{2}\psi^{2}\bmod e$. But $-\delta d$ is not a square modulo $e$. So $e$ divides $\nu$ and either $\theta$ or $\psi$.
\begin{description}
  \item[If $e$ divides $\theta$. ] Then $\nu^{2}-\mu\theta^{4}-\delta e^{2}\theta^{2}\psi^{2}-\delta (e^{2}-4d)\theta^{2}\psi^{2}=\delta e\mu_{2}\psi^{4}$ is divisible by $e^2$. This is a contradiction : $e$ and $\mu_{2}$ are coprime, and $e$ and $\psi$ are coprime (because $\theta$ and $\psi$ are coprime). 
  \item[If $e$ divides $\psi$. ] Then $\nu^{2}-\delta e^{2}\theta^{2}\psi^{2}-\delta (e^{2}-4d)\theta^{2}\psi^{2}-\delta e\mu_{2}\psi^{4}=\mu\theta^{4}$ is divisible by $e^2$. This is a contradiction : $\mu$ is squarefree, and $\theta$ and $\psi$ are coprime.
\end{description}
Both case are impossible; the irreducible polynomial $e=1-C$ does not divides $\mu$.
\fin%

\begin{slemm}\label{lemm_ell_Echap2}We use the notation \ref{nota-section43}. Let $\delta\in k[x]$ be a nontrivial polynomial. Let $p$ be a divisor of $B-C$ coprime to $\delta$. 

Then every element of the image of $\gamma_{\widehat{\CbeC}^{-}_{\delta}}$ is the class in $k(x)^{\times}/k(x)^{\times 2}$ of a squarefree divisor of $\delta (1-C)((1-C)^{2}-4(B-C))$ with either $\mu\sim 1\bmod p$ or $\mu\sim -\delta\bmod p$.
\end{slemm}
\dem%
We use the notation \ref{nota-section43bis}. From Equation \ref{eq-lemm-43-Cchap--delta} we get the congruence $\mu\nu^{2}\equiv (\mu\theta^{2}+\delta e^{2}\psi^{2})^{2}\bmod p$.
\begin{description}
  \item[If $p$ does not divide $\mu\theta^{2}+\delta e^{2}\psi^{2}$. ] Then $\nu$ is invertible modulo $p$ and thus $\mu\sim 1\bmod p$.
  \item[Else ] As $\theta$ and $\psi$ are coprime (and $e$ and $p$ are corpime), the congruence $\mu\theta^{2}\equiv -\delta e^{2}\psi^{2}\bmod p$ gives an equivalence $\mu\sim -\delta\bmod p$.\fin%
\end{description}

\begin{sprop}\label{prop-verif-Cchap-zeta42final} We use the notation \ref{nota-section43}. 
%
We assume that the element \linebreak
$\left(\omega^{2}-\eta^{2}\right)^{2}-4\omega^{2}$ 
is positive.

Then, for every positive element $\zeta\in k^{\times}$, the image of $\gamma_{\widehat{\CbeC}^{-}_{\zeta}}$ is the subgroup of $k(x)^{\times}/k(x)^{\times 2}$ generated by the two classes $\left[ -\zeta\right]$ and $\left[ (1+C)^{2}-4B\right]$. 
\end{sprop}
\dem%
Let $\zeta >0$ be an element of $k$. 
Let $(\alpha ,\beta)$ be a $k(x)$-point of the curve $\widehat{\CbeC}^{-}_{\zeta}$. By assumption, the remainder $-\zeta ((\omega^{2}-\eta^{2})^{2}-4\omega^{2})$ of the Euclidean division of $-\zeta (B-C)$ by $1-C$ is negative. As $k$ is a subfield of $\R$, this remainder is not a square in $k$ and Lemma \ref{ell_Echap1} applies : $\gamma_{\widehat{\CbeC}^{-}_{\zeta}}(\alpha ,\beta )$ is the class in $k(x)^{\times}/k(x)^{\times 2}$ of a squarefree divisor $\mu$ of $(1-C)^{2}-4(B-C)$. Since $[(1-C)^{2}-4(B-C)]=\gamma_{\widehat{\CbeC}^{-}_{\zeta}}(0,0)$, we can assume that $\mu$ is an element of $k$. 

To conclude we use the lemma \ref{lemm_ell_Echap2} : since $\mu$ is constant we have either $\mu\sim 1$ or $\mu\sim -\zeta$.
\fin%

\begin{sprop}\label{prop-verif-Cchap-zetax42final} We use the notation \ref{nota-section43}. 
%
We assume that $\rho$ is different from $0$. We also assume that the hypotheses \ref{hypothesisA}, \ref{hypothesisD}, \ref{hypothesisE}, \ref{hypothesisF} and \ref{hypothesisI} (at least with $n_{1}=1$ and $n_{2}=n_{3}=0$) are satisfied.

Then, for every positive element $\zeta\in k^{\times}$, the image of $\gamma_{\widehat{\CbeC}^{-}_{\zeta x}}$ is the subgroup of $k(x)^{\times}/k(x)^{\times 2}$ generated by the two classes $\left[ -\zeta x\right]$ and $\left[ (1+C)^{2}-4B\right]$. 
\end{sprop}
\dem%
Let $\zeta >0$ be an element of $k$. 
Let $(\alpha ,\beta)$ be a $k(x)$-point of the curve $\widehat{\CbeC}^{-}_{\zeta x}$. From lemma \ref{lem-equation-pol-cbe-ell}, the image $\gamma_{\widehat{\CbeC}^{-}_{\zeta}}(\alpha ,\beta )$ is the class in $k(x)^{\times}/k(x)^{\times 2}$ of a squarefree divisor $\mu$ of $x(1-C)((1-C)^{2}-4(B-C))$. Since $[-\zeta x]=\gamma_{\widehat{\CbeC}^{-}_{\zeta}}(-\zeta x(1-C)^{2},0)$ and $[(1-C)^{2}-4(B-C)]=\gamma_{\widehat{\CbeC}^{-}_{\zeta}}(0,0)$, we can assume that $\mu$ is a divisor of $1-C$.%

The polynomial $B-C$ has only two prime factors: $p_{1}:=x+b_{1}-1-\omega$ and $p_{2}:=x+b_{1}-1+\omega$. They are coprime to $x$ (see hypothesis \ref{hypothesisI}). From Lemma \ref{lemm_ell_Echap2} there exists two integers $m_{1}$ and $m_{2}$ such that
\begin{equation}\label{eqtion-verif-Cchap-42-zeta-bis-cas-pas-carre}
\begin{array}{l}
\mu\sim\left( -\zeta x\right)^{m_{1}}\bmod p_{1}\sim\left( -\zeta\left(1-b_{1}+\omega\right)\right)^{m_{1}}\bmod p_{1}\textrm{  and}\\
\mu\sim\left( -\zeta x\right)^{m_{2}}\bmod p_{2}\sim\left( -\zeta\left(1-b_{1}-\omega\right)\right)^{m_{2}}\bmod p_{2}.\\
\end{array}
\end{equation}

The remainder of the Euclidean division of $-16\zeta x(B-C)$ by $1-C$ is $r:=-2\zeta (1-C(0))((\omega^{2}-\eta^{2})^{2}-4\omega^{2})$. 
\begin{description}
  \item[Case 1: If $r\notin k^{\times 2}$.] Then Lemma \ref{ell_Echap1} applies : $\mu$ is constant.

If $m_{1}$ or $m_{2}$ is even, then one of the two equivalences \ref{eqtion-verif-Cchap-42-zeta-bis-cas-pas-carre} is $\mu\sim\epsilon\sim 1$. Else, by taking the product of the two equivalences \ref{eqtion-verif-Cchap-42-zeta-bis-cas-pas-carre}, we obtain $\epsilon^{2}\sim\left( -\zeta\right)^{2}\left(\left( 1-b_{1}\right)^{2}-\omega^{2}\right)$ i.e. a contradiction with hypothesis \ref{hypothesisI}.
  \item[Case 2: If $r$ is a square i.e. if $\zeta\sim-2(1-C(0))((\omega^{2}-\eta^{2})^{2}-4\omega^{2})$. ] \textrm{ }
  \begin{description}
    \item[Subcase 1 : If $\mu$ is a constant. ] This case is analogous to Case 1.
    \item[Subcase 2 : $\mu =\epsilon (1-C)$ with $\epsilon\in k^{\times}$. ] 
The equivalences \ref{eqtion-verif-Cchap-42-zeta-bis-cas-pas-carre} are 
$$\begin{array}{l}
\epsilon\left(\eta^{2}-\omega^{2}-2\omega\right)\sim\left( -\zeta\left(1-b_{1}+\omega\right)\right)^{m_{1}}\textrm{ and }\\
\epsilon\left(\eta^{2}-\omega^{2}+2\omega\right)\sim\left( -\zeta\left(1-b_{1}-\omega\right)\right)^{m_{2}}.\\
\end{array}$$
Taking the product of theses equivalences, we obtain
$$\left(\left(\eta^{2}-\omega^{2}\right)^{2}-4\omega^{2}\right)\sim \left( -\zeta\right)^{m_{1}+m_{2}}\left(1-b_{1}+\omega\right)^{m_{1}}\left(1-b_{1}-\omega\right)^{m_{2}}.$$
Since $\zeta\sim-2(1-C(0))((\omega^{2}-\eta^{2})^{2}-4\omega^{2})$ we have a contradiction with hypothesis \ref{hypothesisA}, \ref{hypothesisD}, \ref{hypothesisE} and \ref{hypothesisF}. \fin%
\end{description}
\end{description}

\end{subsection}

\begin{subsection}{The image of $\Pi_{\CbeC^{+}_{\delta}}.$}
\begin{snota}\label{nota-sect-44-gene1} Let $k$ be a characteristic $0$ field. Let $E,$ $D,$ $\delta\in k[x]$ be three nontrivial polynomials such that $1-2E,$ $(1-E)^{2}-D$ and $E^{2}-D$ are nontrivial. We consider the hyperelliptic curve $\CbeH$ defined on $k(x)$ by the affine equation 
$$\CbeH :z^{2}=(y+\delta (1-E))(y-\delta E)(y+\delta E)(y^{2}-\delta^{2}D).$$
\end{snota}
\rmk%
The hypotheses on $E$ and $D$ implies that the polynomial \linebreak$(y+\delta (1-E))(y-\delta E)(y+\delta E)(y^{2}-\delta^{2}D)$ is squarefree.

\begin{snota}\label{nota-sect-44-gene2} We use the notation \ref{nota-sect-44-gene1}. Consider the four polynomials $f_{1}(y):=y+\delta (1-E)$, $f_{2}(y):=y-\delta E$, $f_{3}:=y+\delta$ and $f_{4}(y):=y^{2}-\delta^{2}D$. 

We assume that $f_{4}(y)$ is irreducible (i.e. that $D$ is a square in $k(x)$). Denote by $K_{i}$ the field $k(x)[T]/(f_{i}(T))$ and by $y_{i}$ the class of $T$ in $K_{i}$. Let $\pi_{\CbeH}:\textrm{Jac}(\CbeH )(k(x))\longrightarrow\displaystyle\prod_{i=1}^{4}K_{i}^{\times}/K_{i}^{\times 2}$ be morphism defined by Proposition \ref{def-Cas-Sch+im-in-ker} and $\pi_{\CbeH ,i}:\textrm{Jac}(\CbeH )(k(x))\longrightarrow\displaystyle K_{i}^{\times}/K_{i}^{\times 2}$ be its $i$-th component. In this section we study the homorphism $\Xi_{\CbeH}:\textrm{Jac}(\CbeH )(k(x))\longrightarrow\displaystyle\prod_{i=1}^{4}k(x)^{\times}/k(x)^{\times 2}$ with $i$-th component the homomorphism $\Xi_{\CbeH,i}:=N_{K_{i}/k(x)}\circ\pi_{\CbeH ,i}$.
\end{snota}

\begin{snota}\label{nota-sect-44-gene3} We use the notation \ref{nota-sect-44-gene1} and \ref{nota-sect-44-gene2}. We assume that %
$1-2E$ and $(1-E)^{2}-D$ are coprime, %
$(1-E)^{2}-D$ and $E$ are coprime, %
$1-2E$ and $D$ are coprime, and %
$E$ and $D$ are coprime.

Let $\alpha$ be an element of the image $\textrm{Im}(\Xi_{\CbeH})$. From proposition \ref{caracterisationimagePi}, we know that $\alpha =([\mu_{1,2}\mu_{1,3}\mu_{1,4}],[\mu_{1,2}\mu_{2,3}\mu_{2,4}],[\mu_{1,3}\mu_{2,3}\mu_{3,4}],[\mu_{1,4}\mu_{2,4}\mu_{3,4}])$ with
\begin{itemize}
\item[* ] $\mu_{1,2}\in k[x]$ a squarefree divisor of $\delta$,
\item[* ] $\mu_{1,3}\in k[x]$ a squarefree divisor of $\delta (1-2E),$
\item[* ] $\mu_{1,4}\in k[x]$ a squarefree divisor of $\delta ((1-E)^{2}-D),$
\item[* ] $\mu_{2,3}\in k[x]$ a squarefree divisor of $\delta E(E^{2}-D),$
\item[* ] $\mu_{2,4}\in k[x]$ a squarefree divisor of $\delta (E^{2}-D)$ and
\item[* ] $\mu_{3,4}\in k[x]$ a squarefree divisor of $\delta (E^{2}-D)$.
\end{itemize}
Consider the following polynomials:
\begin{itemize}
\item[* ] $\widetilde{\mu}_{1,2}=1$, $\widetilde{\mu}_{1,3}=\mu_{1,3}\mu_{1,2}$, $\widetilde{\mu}_{1,4}=\mu_{1,4}$,
\item[* ] $\widetilde{\mu}_{2,3}=\mu_{2,3}\mu_{1,2}\mu_{2,4}$, $\widetilde{\mu}_{2,4}=1$ and $\widetilde{\mu}_{3,4}=\mu_{3,4}\mu_{2,4}$,
\end{itemize}
Replacing $\mu_{i,j}$ by a squarefree polynomial in the same class in $k(x)^{\times}/k(x)^{\times 2}$ than $\widetilde{\mu}_{i,j}$ we assume without loss of generality that $\mu_{1,2}=\mu_{2,4}=1$.
\end{snota}

\begin{sprop}\label{J++-mod-1-C} We use the notation \ref{nota-sect-44-gene1}, \ref{nota-sect-44-gene2} and \ref{nota-sect-44-gene3}. Let $p$ be a prime factor of $E$ such that $v_{p}(E)=1$ and $v_{p}(\delta )=v_{p}(D)=0.$

Then for every $\alpha$ in the image of $\Xi_{\CbeH}$ there exists four polynomials \linebreak$\alpha_{1},$ $\alpha_{2}$, $\alpha_{3},$ $\alpha_{4}\in k[x]$ such that
$\alpha =\left(\left[\alpha_{1}\right] ,\left[\alpha_{2}\right] ,\left[\alpha_{3}\right] ,\left[\alpha_{4}\right]\right)$ and  
\begin{enumerate}
\item $\alpha_{2}\alpha_{3}\sim 1\bmod p$ whenever $v_{p}(\alpha_{2})$ is even; 
\item $\alpha_{2}\alpha_{3}\sim -\delta D\bmod p$ whenever $v_{p}(\alpha_{2})$ is odd.
\end{enumerate}
\end{sprop}

\begin{slemm}\label{lemm-J++-mod-separation-cas-1-C} We keep the notation and assumptions from Proposition \ref{J++-mod-1-C}. 
%
Let $\beta\in\textrm{Jac}(\CbeH )(k(x))$ be a $k(x)$-point of $\textrm{Jac}(\CbeH )$. 


Then $\beta$ can be written in the form $\beta =<u,v>+T$ with $T\in\textrm{Jac}(\CbeH )(k(x))[2]$ a $2$-torsion point and $u$ a polynomial of degree at most $2$, coprime to $f_{1}$, $f_{2}$, $f_{3}$ and $f_{4}$, such that  
\begin{enumerate} 
\item either $\deg (u)\le 2$ and $v_{p}(u(\delta E))$ is even, 
\item or $\deg (u)\le 1$ and $v_{p}(u(\delta E))$ is odd.
\end{enumerate}
\end{slemm}
\dem
Let $V_{p}:k(x)^{\times}/k(x)^{\times 2}\longrightarrow\Z/2\Z$ be the map induced by $v_{p}$. Let $n_{2}\in\{ 0,1\}$ represent the class $V_{p}(\Xi_{\CbeH ,2}(\beta ))$. By assumption the valuation $v_{p}(2E(E^{2}-D))$ is odd. The image of the $2$-torsion point $<y-\delta E,0>$ by $\Xi_{\CbeH ,2}$ being $[2E(E^{2}-D)]$, we have $V_{p}(\Xi_{\CbeH ,2}(\beta +n_{2}<y-\delta E,0>))=0$.

Let $(\widetilde{u},\widetilde{v})$ be Mumford's representation for $\beta +n_{2}<y-\delta E,0>$. As the image of $<\widetilde{u},\widetilde{v}>$ under $\Xi_{\CbeH ,2}$ is the class $[(-1)^{\deg (\widetilde{u})}\widetilde{u}(\delta E)]$, the valuation $v_{p}(\widetilde{u}(\delta E))$ is even. The point $<\widetilde{u},\widetilde{v}>$ does not satisfy all the conditions of the lemma : $\widetilde{u}$ and $f_{1}f_{2}f_{3}f_{4}$ may have a common factor. To solve this problem we consider three polynomials: the polynomial $u_{f}:=\textrm{gcd}\left(\widetilde{u},f_{1}f_{2}f_{3}f_{4}\right)$, the polynomial $u:=\widetilde{u}/u_{f}$, and the remainder $v$ of the Euclidean division of $\widetilde{v}$ by $u$. 

If $deg (u) =2$, then the point $<u,v>=<\widetilde{u},\widetilde{v}>$ satisfies condition $1.$ Otherwise the degree $deg (u)$ is at most $1$ and the point $<u,v>$ satisfies either condition $1$ or condition $2$ (the correct case depends only on the parity of $v_{p}(u(\delta E))$).
\fin\\
\\
\textbf{Proof of Proposition \ref{J++-mod-1-C}. }Let $V_{p}:k(x)^{\times}/k(x)^{\times 2}\longrightarrow\Z/2\Z$ be the map induced by $v_{p}$. We use the notation \ref{nota-sect-44-gene3}. Let $\beta$ be a $k(x)$-point of $\textrm{Jac}(\CbeH )$ with $\alpha =\Xi_{\CbeH}(\beta )$. The main idea of the proof is to use the equality 
\begin{equation}\label{eqtn-base-lemm-J++-mod-B-C-1-C-44}
u(-\delta E)=u(\delta E)-2\delta Eu_{1}
\end{equation} 
(which holds for every polynomial $u(y)=u_{2}y^{2}+u_{1}y+u_{0}\in k(x)[y]$ of degree at most $2$).
\begin{description}
        \item[Step 1. ] Let $\mathcal{S}$ be the set of all $\widetilde{\alpha}\in\left( k(x)^{\times}/k(x)^{\times 2}\right)^{4}$ that can be written as $\widetilde{\alpha}=\left(\left[\widetilde{\alpha}_{1}\right] ,\left[\widetilde{\alpha}_{2}\right] ,\left[\widetilde{\alpha}_{3}\right] ,\left[\widetilde{\alpha}_{4}\right]\right)$ for some polynomials $\widetilde{\alpha}_{1},$ $\widetilde{\alpha}_{2}$, $\widetilde{\alpha}_{3},$ $\widetilde{\alpha}_{4}\in k[x]$ 
such that
\begin{enumerate}
\item $\widetilde{\alpha}_{2}\widetilde{\alpha}_{3}\sim 1\bmod p$ whenever $v_{p}(\widetilde{\alpha}_{2})$ is even, and
\item $\widetilde{\alpha}_{2}\widetilde{\alpha}_{3}\sim -\delta D\bmod p$ si $v_{p}(\widetilde{\alpha}_{2})$ is odd.
\end{enumerate}
From Lemma \ref{lemm-J++-mod-separation-cas-1-C}, $\beta =<u,v>+T$ with $T\in\textrm{Jac}(\CbeH )(k(x))[2]$ and $u\in k(x)[y]$ of degree at most $2$, coprime to $f_{1}$, $f_{2}$, $f_{3}$ and $f_{4}$, such that  
\begin{enumerate} 
\item either $\deg (u)\le 2$ and $v_{p}(u(\delta E))$ is even, 
\item or $\deg (u)\le 1$ and $v_{p}(u(\delta E))$ is odd.
\end{enumerate}
As $\Xi_{\CbeH}$ is a homomorphism and since the $2$-torsion of $\textrm{Jac}(\CbeH )(k(x))$ is a subset of $\mathcal{S}$, the image $\alpha=\Xi_{\CbeH}(\beta )$ belongs to $\mathcal{S}$ if and only if $\Xi_{\CbeH}(<u,v>)$ belongs to $\mathcal{S}$. In particular, without loss of generality, we can assume that $\beta$ is equal to $<u,v>$.
        \item[Step 2. ] Write $u(y)=\frac{u_{2}}{\lambda}y^{2}+\frac{u_{1}}{\lambda}y+\frac{u_{0}}{\lambda}$ with $u_{0},$ $u_{1},$ $u_{2}$, $\lambda\in k[x]$ four coprime polynomials. The polynomial $u$ being monic, $\lambda$ is the leading coefficient of $\widetilde{u}(y):=u_{2}y^{2}+u_{1}y+u_{0}.$ In particular $u_{0},$ $u_{1}$ and $u_{2}$ are coprime.

The image $\Xi_{\CbeH ,1}(<u,v>)$ is equal to the class $[(-1)^{\deg (u)}u(\delta (E-1))]$. It is also the class of the divisor $\mu_{1,3}\mu_{1,4}$ of $\delta^{2}(1-2E)((1-E^{2})-D)$ (see the notation \ref{nota-sect-44-gene3}). Thus, since $\delta^{2}(1-2E)((1-E^{2})-D)$ is coprime to $p$, the valuation $v_{p}u(\delta (E-1))$ is even. 


Assume that $v_{p}(\lambda )$ is odd. Then the valuation $v_{p}(\widetilde{u}(\delta (E-1)))$ is odd. The coefficient $u_{2}$ being either $0$ or $\lambda$, the prime $p$ divides $u_{2}$. As a consequence the congruence 
$-\delta u_{1}+u_{0}\equiv\widetilde{u}(\delta (E-1))\bmod p\equiv 0\bmod p$ holds. But $p$ is not a common divisor of $u_{0},$ $u_{1}$ and $u_{2}$ (they are coprime). Thus $p$ is coprime to $u_{0}$ and coprime to $u_{1}$. In particular $u_{2}$ is equal to $\lambda$ ($\lambda$ is divisible by $p$ and equal to $u_{0}$, $u_{1}$ or $u_{2}$). We have a contradiction :
\begin{itemize}
        \item[* ] the polynomial $u_{2}=\lambda$ being nontrivial, the degree of $u(y)$ is 2. From our hypothesis on $<u,v>$ (see Lemma \ref{lemm-J++-mod-separation-cas-1-C}), the valuation $v_{p}(u(\delta E))$ is even;
        \item[* ] because of the congruence $\widetilde{u}(\delta E)\equiv u_{0}\bmod p$, the polynomial $p$ does not divide $\widetilde{u}(\delta E)$; thus $v_{p}(\widetilde{u}(\delta E))$ is even i.e. $v_{p}(u(\delta E))$ is odd.
\end{itemize} 
From this contradiction we deduce that $v_{p}(\lambda )$ is even.
        \item[Step 3. ] 
Following Step $1$, we only have three cases to look at. 
        \begin{description}
                \item[Case 1: if $v_{p}(\widetilde{u}(\delta E))$ is trivial. ] Following Equation \ref{eqtn-base-lemm-J++-mod-B-C-1-C-44} we have \linebreak$\widetilde{u}(\delta E)\equiv \widetilde{u}(-\delta E)\bmod p.$ Moreover $\widetilde{u}(\delta E)$ is invertible modulo $p$. This implies the triviality of the image 
$$\begin{array}{rcl}
\Xi_{\CbeH ,2}(<u,v>)\Xi_{\CbeH ,3}(<u,v>)&=&\left[ u(\delta E)u(-\delta E)\right]\\
&=&\left[\widetilde{u}(\delta E)\widetilde{u}(-\delta E)\right].\\
\end{array}$$
                \item[Case 2: if $v_{p}(\widetilde{u}(\delta E))$ is even and nontrivial. ] Since it divides \linebreak$\delta^{2}(1-2E)(E^{2}-D)$ the polynomial $\mu_{1,3}\mu_{3,4}$ 
is coprime to $p$. Using the equalities $\left[(-1)^{\deg (u)}u(\delta E)\right] =\Xi_{\CbeH ,2}(<u,v>)=\left[\mu_{2,3}\right]$, we get the parity of $v_{p}(\mu_{2,3})$. Hence from the equality 
$$\left[ (-1)^{\deg (u)}u\left( -\delta E\right)\right] =\Xi_{\CbeH ,2}(<u,v>)=\left[\mu_{1,3}\mu_{2,3}\mu_{3,4}\right]$$ 
we deduce that $v_{p}(\widetilde{u}(-\delta E))=v_{p}(u(-\delta E))+v_{p}(\lambda )$ is even.

The polynomial $\widetilde{u}( -\delta E)=\widetilde{u}(\delta E)-2\delta Eu_{1}$ is divisible by $p$ (because $p$ divise $E$ and $\widetilde{u}(\delta E)$). Thus $\widetilde{u}( -\delta E)$ is divisible by $p^{2}$. 

Since $v_{p}(E)=1$ and since $\widetilde{u}(\delta E)$ and $\widetilde{u}(-\delta E)$ are divisible by $p^{2}$, the equation $2\delta Eu_{1}=\widetilde{u}(\delta E)-\widetilde{u}(-\delta E)$ implies the divisibility of $u_{1}$ by $p$. From the divisibility of $u_{0}=\widetilde{u}(\delta E)-u_{2}\delta^{2}E^{2}-u_{1}\delta E$ and $u_{1}$ by $p$ we deduce 
\begin{itemize}
\item[* ] the coprimality of $u_{2}$ and $p$  and thus the equality of $u_{2}$ and $\lambda$ (because $u_{0}$, $u_{1}$ and $u_{2}$ are coprime);
\item[* ] $(-1)^{\deg (\widetilde{u})}\lambda \widetilde{u}(\delta (E-1))\equiv\delta^{2}(E-1)^{2}u_{2}^{2}\bmod p$; 
\item[* ] $\begin{array}[t]{rcl}
N_{K_{4}/k(x)}\left( (-1)^{\deg (\widetilde{u})}\widetilde{u}\right)&=&\left(u_{0}+\delta^{2}Du_{2}\right)^{2}-\delta^{2}Du_{1}^{2}\\
&\equiv &\delta^{4}D^{2}u_{2}^{2}\bmod p.\\
\end{array}$
\end{itemize}
We can conclude thanks to Proposition \ref{def-Cas-Sch+im-in-ker}: the image
$$\Xi_{\CbeH ,2}(<u,v>)\Xi_{\CbeH ,3}(<u,v>)=\Xi_{\CbeH ,1}(<u,v>)\Xi_{\CbeH ,4}(<u,v>)$$
is equal to the class 
$$\left[ (-1)^{\deg (\widetilde{u})}\lambda\widetilde{u}(\delta (E-1))N_{K_{4}/k(x)}\left( (-1)^{\deg (\widetilde{u})}\widetilde{u}\right)\right] =\left[ 1\right].$$
                \item[Case 3: if $v_{p}(\widetilde{u}(\delta E))$ is odd i.e. if the degree of $\widetilde{u}(y)$ is at most $1$. ] Then $p$ divides $u_{0}=\widetilde{u}(\delta E)-\delta Eu_{1}$. Hence we have 
\begin{itemize}
        \item[* ] $N_{K_{4}/k(x)}((-1)^{\deg (u)}\widetilde{u})=u_{0}^{2}-\delta^{2}Du_{1}^{2}\equiv -\delta^{2}Du_{1}^{2}\bmod p$ and
        \item[* ] $\widetilde{u}(\delta (E-1))=\delta (E-1)u_{1}+u_{0}\equiv -\delta u_{1}\bmod p.$
\end{itemize}
Morevover $u_{0}$ and $u_{1}$ being coprime, $p$ does not divides $u_{1}$ and we have $u_{1}=\lambda$. To conclude we use Proposition \ref{def-Cas-Sch+im-in-ker}: the image
$$\Xi_{\CbeH ,2}(<u,v>)\Xi_{\CbeH ,3}(<u,v>)=\Xi_{\CbeH ,1}(<u,v>)\Xi_{\CbeH ,4}(<u,v>)$$
is equal to the class 
$$\left[ (-1)^{\deg (u)}\lambda\widetilde{u}(\delta (E-1))N_{K_{4}/k(x)}\left( (-1)^{\deg (u)}\widetilde{u}\right)\right] =\left[ -\delta D\right]\textrm{. \fin}$$
\end{description}
\end{description}
%

\begin{sprop}\label{J++-mod-B-C} We use the notation \ref{nota-sect-44-gene1}, \ref{nota-sect-44-gene2} and \ref{nota-sect-44-gene3}. Let $p\in k[x]$ be a prime such that $v_{p}(E^{2}-D)=1$ and $v_{p}(\delta )=v_{p}(E)=v_{p}(1-2E)=0.$

Then for every $\alpha$ in the image of $\Xi_{\CbeH}$ there exists four polynomials $\alpha_{1},$ $\alpha_{2}$, $\alpha_{3},$ $\alpha_{4}\in k[x]$ such that
$\alpha =\left(\left[\alpha_{1}\right] ,\left[\alpha_{2}\right] ,\left[\alpha_{3}\right] ,\left[\alpha_{4}\right]\right)$ and 
\begin{enumerate}
\item $\alpha_{1}\sim 1\bmod p$ or $\alpha_{1}\sim 1-2E\bmod p$ whenever $v_{p}(\alpha_{2}\alpha_{3})$ is even;
\item $\alpha_{1}\sim \delta\bmod p$ ou $\alpha_{1}\sim \delta (1-2E)\bmod p$ whenever $v_{p}(\alpha_{2}\alpha_{3})$ is odd.
\end{enumerate}
\end{sprop}
\dem%
Let $V_{p}:k(x)^{\times}/k(x)^{\times 2}\longrightarrow\Z/2\Z$ be the map induced by $v_{p}$. Let $\beta$ be a $k(x)$-point of $\textrm{Jac}(\CbeH )$ and denote by $\alpha$ the image $\Xi_{\CbeH}(\beta )$. We use the notation \ref{nota-sect-44-gene3} relative to $\alpha$. The proof is analogous to the proof of Proposition \ref{J++-mod-1-C}.
\begin{description}
        \item[Step 1. ] As in the proof of Proposition \ref{J++-mod-1-C}, and because of the three following facts
\begin{enumerate}
        \item the image of the $2$-torsion of $\textrm{Jac}(\CbeH )(k(x))$ by $\Xi_{\CbeH}$ contains
$$\begin{array}{l}
\left(\left[\left( 1-E\right)^{2}-D\right] ,\left[ E^{2}-D\right] ,\left[ E^{2}-D\right] ,\left[\left( 1-E\right)^{2}-D\right]\right)\\
=\Xi_{\CbeH}\left( <y^{2}-\delta^{2}D,0>\right)\textrm { and }\\
\\
\left(\left[\delta\right] ,\left[ 2E\left( E^{2}-D\right)\right] ,\left[ 2\delta E\right] ,\left[ E^{2}-D\right]\right) =\Xi_{\CbeH}\left( <y-\delta E,0>\right)\\
\end{array}$$
        \item $v_{p}(E)=0$ and $v_{p}(E^{2}-D)=1$, and
        \item if $<u,v>$ is a $k(x)$-point of $\textrm{Jac}(\CbeH )$ with $u$ and $f_{1}f_{1}f_{1}f_{1}$ coprime, then $\Xi_{\CbeH ,2}(<u,v>)=\left[ u\left(\delta E\right)\right]$ and $\Xi_{\CbeH ,3}(<u,v>)=\left[ u\left( -\delta E\right)\right]$,
\end{enumerate}
we can assume without loss of generality that $\alpha =\Xi_{\CbeH}(<u,v>)$ with %
$u$ a polynomial of degree at most $2$, coprime to $f_{1}$, $f_{2}$, $f_{3}$ and $f_{4}$, such that 
\begin{enumerate} 
\item either $\deg (u)\le 2$, and both $v_{p}(u(\delta E))$ and $v_{p}(u(-\delta E))$ are even, 
\item or $\deg (u)\le 1$ and $v_{p}(u(\delta E)u(-\delta E))$ is odd.
\end{enumerate}

        \item[Step 2. ] Write $u(y)=\frac{u_{2}}{\lambda}y^{2}+\frac{u_{1}}{\lambda}y+\frac{u_{0}}{\lambda}$ with $u_{0},$ $u_{1},$ $u_{2}$, $\lambda\in k[x]$ four coprime polynomials. The polynomial $u$ being monic, $\lambda$ is the leading coefficient of $\widetilde{u}(y):=u_{2}y^{2}+u_{1}y+u_{0}.$ In particular $u_{0},$ $u_{1}$ and $u_{2}$ are coprime.

Assume that $v_{p}(\lambda )$ is odd. 
By definition of the $\mu_{i,j}$, we have \linebreak$V_{p}([(-1)^{\deg (u)}u(\delta (E-1))])
=V_{p}([\mu_{1,3}\mu_{1,4}])$. Morevover $\mu_{1,3}\mu_{1,4}$ is coprime to $p$ because:
\begin{itemize}
        \item[* ] it is a divisor of $\delta^{2}(1-2E)((1-E)^{2}-D$, and 
        \item[* ] $\delta^{2}(1-2E)((1-E)^{2}-D)$ is coprime to $p$ (since $v_{p}(\delta )$ and $v_{p}(1-2E)$ are trivial and since $(1-E)^{2}-D=(1-2E)+(E^{2}-D)$). 
\end{itemize}
Hence $v_{p}(u(\delta (E-1)))$ is even and $v_{p}(\widetilde{u}(\delta (E-1))$ 
is odd. In particular $p$ divides $\widetilde{u}(\delta (E-1))$.

The coefficient $u_{2}$ being either $0$ or $\lambda$ the prime $p$ is a divisor of $u_{2}$.  
This divisibility gives a congruence 
\begin{equation}\label{eqt-646-step2-u0-u1-congruence-etoile}
\delta (E-1)u_{1}+u_{0}=\widetilde{u}(\delta (E-1))-u_{2}\delta^{2}(E-1)^{2}\bmod p\equiv 0\bmod p.
\end{equation}
The polynomial $p$ is not a common divisor of $u_{0},$ $u_{1}$ and $u_{2}$ (they are coprime). As a consequence $p$ is coprime to $u_{0}$ and coprime to $u_{1}$. In particular $u_{2}$ is equal to $\lambda$ ($\lambda$ is divisible by $p$ and equal to $u_{0}$, $u_{1}$ or $u_{2}$). We have a contradiction :
\begin{itemize}
        \item[* ] the polynomial $u_{2}=\lambda$ being nontrivial, the degree of $u(y)$ is 2. From our hypothesis on $<u,v>$ (see Step $1$), the valuation $v_{p}(u(\delta E))$ is even i.e. $v_{p}(\widetilde{u}(\delta E))$ is odd;
        \item[* ] since $\widetilde{u}(\delta E)\equiv\delta Eu_{1}+u_{0}\bmod p\equiv\delta u_{1}\bmod p$ (see Congruence \ref{eqt-646-step2-u0-u1-congruence-etoile}), the polynomial $p$ does not divide $\widetilde{u}(\delta E)$; thus $v_{p}(\widetilde{u}(\delta E))$ is even.
\end{itemize} 
The existence of a contraction means that $v_{p}(\lambda )$ is even.
        \item[Step 3. ] For every $w\in k(x)[y]$ we denote by $\Psi_{w}(T)$ the resultant 
$$\Psi_{w}(T):=\textrm{res}_{y}\left((-1)^{\deg (w)}w(y),y^{2}-T\right) .$$
This notation is motivated by the equalities 
$$\begin{array}{cl}
&\Xi_{\CbeH ,4}(<u,v>) =\left[\Psi_{u}\left(\delta^{2}D\right)\right] =\left[\Psi_{\widetilde{u}}\left(\delta^{2}D\right)\right]\\
\textrm{and }&\Xi_{\CbeH ,2}(<u,v>)\Xi_{\CbeH ,3}(<u,v>) =\left[\Psi_{\widetilde{u}}\left(\delta^{2}E^{2}\right)\right] .\\
\end{array}$$ 
In particular we deduce from Proposition \ref{def-Cas-Sch+im-in-ker} that 
\begin{equation}\label{Calcul_Xi_1-C+delta44-mod-(E2-D)}
\begin{array}{rcl}
\Xi_{\CbeH ,1}(<u,v>) &=&\displaystyle\prod_{i=2}^{4}\Xi_{\CbeH ,i}(<u,v>)\\
&=&\left[\Psi_{\widetilde{u}}\left(\delta^{2}E^{2}\right)\Psi_{\widetilde{u}}\left(\delta^{2}D\right)\right] .\\
\end{array}
\end{equation}
The proof of Proposition \ref{J++-mod-B-C} is based on Taylor's formula for $\Psi_{\widetilde{u}}$:
\begin{equation}\label{equation-C+delta44-mod-(E2-D)}
\begin{array}{rcl}
\Psi_{\widetilde{u}}\left( \delta^{2}D\right) -\Psi_{\widetilde{u}}\left(\delta^{2}E^{2}\right) &=&\delta^{2}\left( D-E^{2}\right)\Psi_{\widetilde{u}}^{'}\left(\delta^{2}E^{2}\right)\\
&&+u_{2}^{2}\delta^{4}\left( D-E^{2}\right)^{2}\\
\end{array}
\end{equation}
with $\Psi_{\widetilde{u}}'(T)=2u_{2}(u_{0}+Tu_{2})-u_{1}^{2}$ the usual derivative of $\Psi_{\widetilde{u}}(T)$. We will especially use the congruence $\Psi_{\widetilde{u}}\left( \delta^{2}D\right)\equiv\Psi_{\widetilde{u}}\left(\delta^{2}E^{2}\right)\bmod p$.
\begin{description}
        \item[Case 1: if $v_{p}(\widetilde{u}(\delta E))=v_{p}(\widetilde{u}(-\delta E))=0$. ] Then $\Psi_{\widetilde{u}}\left(\delta^{2}E^{2}\right)$ has an inverse modulo $p$. Hence the reduction modulo $p$ of Equation \ref{equation-C+delta44-mod-(E2-D)} gives 
an equivalence $\Psi_{\widetilde{u}}\left( \delta^{2}D\right)\sim\Psi_{\widetilde{u}}\left(\delta^{2}E^{2}\right)\bmod p$. We conclude by invoking Equality \ref{Calcul_Xi_1-C+delta44-mod-(E2-D)}: the class $\Xi_{\CbeH ,1}(<u,v>)$ is trivial.
        \item[Case 2: if $v_{p}(\widetilde{u}(\delta E)\widetilde{u}(-\delta E))$ is even and greater or equal to $1$. ] Since $\Psi_{\widetilde{u}}\left( \delta^{2}D\right)\equiv\Psi_{\widetilde{u}}\left(\delta^{2}E^{2}\right)\bmod p\equiv\widetilde{u}(\delta E)\widetilde{u}(-\delta E)\bmod p$, the polynomial $p$ divides $\Psi_{\widetilde{u}}\left( \delta^{2}D\right)$.

The valuation $v_{p}(\widetilde{u}(\delta (E-1))$ is even (see step $2$). From the equality $\left[\widetilde{u}\left(\delta\left( E-1\right)\right)\right] =\Xi_{\CbeH ,1}(<u,v>) =\left[\Psi_{\widetilde{u}}\left(\delta^{2}E^{2}\right)\Psi_{\widetilde{u}}\left(\delta^{2}D\right)\right]$ we deduce that $v_{p}\left(\Psi_{\widetilde{u}}\left(\delta^{2}E^{2}\right)\right)$ and $v_{p}\left(\Psi_{\widetilde{u}}\left(\delta^{2}D\right)\right)$ have the same parity. In case $2$ these valuations are even because $\Psi_{\widetilde{u}}\left(\delta^{2}E^{2}\right)=\widetilde{u}(\delta E)\widetilde{u}(-\delta E)$. In particular $p^{2}$ divides $\Psi_{\widetilde{u}}\left(\delta^{2}E^{2}\right)$ and $\Psi_{\widetilde{u}}\left(\delta^{2}D\right)$. This is compatible with Equation \ref{equation-C+delta44-mod-(E2-D)} only if $p$ divides $\Psi_{\widetilde{u}}'\left(\delta^{2}E^{2}\right)$ (the valuation $v_{p}(\delta^{2}(D-E^{2})$ is equal to $1$).
\begin{description}
        \item[subcase (a): if $p$ does not divide $\widetilde{u}(-\delta E)$. ] Then $\widetilde{u}(\delta E)$ is divisible by $p$. As a consequence we have a congruence 
\begin{equation}\label{valeur-u0-modE2-DJ++44-premiere}
u_{0}\equiv -\delta^{2}E^{2}u_{2}-\delta Eu_{1}\bmod p.
\end{equation}
Since $2\delta Eu_{1}=\widetilde{u}(\delta E)-\widetilde{u}(-\delta E)$, the polynomial $p$ does not divide $u_{1}$. Using Equation \ref{valeur-u0-modE2-DJ++44-premiere} we deduce from the divisibility of $\Psi_{\widetilde{u}}'\left(\delta^{2}E^{2}\right) =2u_{2}(u_{0}+\delta^{2}E^{2}u_{2})-u_{1}^{2}$ by $p$ that the prime $p$ divides $u_{1}(\delta Eu_{2}+u_{1})$. Using the coprimality of $p$ and $u_{1}$ this implies that 
\begin{equation}\label{valeur-u0-modE2-DJ++44-premiere-bis}
u_{1}\equiv -2\delta Eu_{2}\bmod p.
\end{equation}
As $p$ and $u_{1}$ are coprime, $p$ does not divide $u_{2}$. Since it is nontrivial, $u_{2}$ is equal to $\lambda$. 
Thus, using Equations \ref{valeur-u0-modE2-DJ++44-premiere} and \ref{valeur-u0-modE2-DJ++44-premiere-bis}, we get a congruence $(-1)^{\deg (\widetilde{u})}\lambda\widetilde{u}(\delta (E-1))\equiv u_{2}^{2}\delta^{2}\bmod p$ i.e. the triviality of $\Xi_{\CbeH ,1}(<u,v>)=\left[ (-1)^{\deg (\widetilde{u})}\lambda\widetilde{u}(\delta (E-1))\right]$.
        \item[subcase (b): if $p$ does not divide $\widetilde{u}(\delta E)$. ] Doing computa-\linebreak tions analogous to subcase (a), we get two congruences \linebreak$u_{0}\equiv -\delta^{2}E^{2}u_{2}+\delta Eu_{1}\bmod p$ and $u_{1}\equiv 2\delta Eu_{2}\bmod p$. 
Those congruences implies the coprimality of $u_{2}$ and $p$. In particular $u_{2}$ and $\lambda$ are equals and the degree of $\widetilde{u}$ is $2$. This means that $(-1)^{\deg (\widetilde{u})}\lambda\widetilde{u}(\delta (E-1))\equiv u_{2}^{2}\delta^{2}(1-2E)^{2}\bmod p$ i.e. that the class $\Xi_{\CbeH ,1}(<u,v>)$ is trivial.
        \item[subcase (c): if $p$ divides $\widetilde{u}(-\delta E)$ and $\widetilde{u}(\delta E)$. ] Then $p$ divides $2(\delta^{2}E^{2}u_{2}+u_{0})=\widetilde{u}(\delta E)+\widetilde{u}(-\delta E)$ and $u_{1}$ (because $2\delta Eu_{1}=\widetilde{u}(\delta E)-\widetilde{u}(-\delta E)$). This is possible only if $u_{2}$ is coprime to $p$ (since $p$ can not divide both $u_{0}$, $u_{1}$ and $u_{2}$). As a consequence $u_{2}$ and $\lambda$ are equal. 

Using the divisibility of $2(\delta^{2}E^{2}u_{2}+u_{0})$ and $u_{1}$ by $p$ and the equalities $\deg (\widetilde{u})=2$ and $\lambda =u_{2}$ we obtain an equivalence $(-1)^{\deg (\widetilde{u})}\lambda \widetilde{u}(\delta (E-1))\equiv u_{2}^{2}\delta^{2}(1-2E)\bmod p.$ Hence we have $\Xi_{\CbeH ,1}(<u,v>)=\left[(-1)^{\deg (\widetilde{u})}\lambda \widetilde{u}\left(\delta\left( E-1\right)\right)\right] =\left[ 1-2E\right]$.
\end{description}
        \item[Case 3: if $\deg (\widetilde{u})=1$ and $v_{p}(\widetilde{u}(\delta E)\widetilde{u}(-\delta E))$ is odd. ]\textrm{ }
        \begin{description}
                \item[Subcase (a): if $v_{p}(\widetilde{u}(\delta E)$ is even. ] Then $v_{p}(\widetilde{u}(-\delta E))$ is odd. In particular we have $u_{0}\equiv\delta Eu_{1}\bmod p$. As $u_{0}$ and $u_{1}$ are coprime, $p$ does not divide $u_{1}$. Since it is nontrivial, $u_{1}$ is equal to $\lambda$. From the congruence 
$$(-1)^{\deg (\widetilde{u})}\lambda \widetilde{u}(\delta (E-1))\equiv u_{1}^{2}\delta (1-2E)\bmod p$$ 
we deduce that $\Xi_{\CbeH ,1}(<u,v>)=\left[\delta (1-2E)\right]$.
                \item[Subcase (b): if $v_{p}(\widetilde{u}(\delta E)$ is odd. ] Since $u_{0}\equiv -\delta Eu_{1}\bmod p$ the prime $p$ does not divide $u_{1}$ (the polynomials $u_{0}$ and $u_{1}$ are coprime) and the congruence 
$$(-1)^{\deg (\widetilde{u})}\lambda \widetilde{u}(\delta (E-1))\equiv u_{1}^{2}\delta\bmod p$$ 
holds. In particular $\Xi_{\CbeH ,1}(<u,v>)$ is equal to $\left[\delta\right]$.
\fin%
        \end{description}
\end{description}
\end{description}

\begin{sprop}\label{verification-J++1} We use the notation \ref{Notations_Principales_21}. We also use the notation \ref{nota_genre2} and \ref{nota-equiv-mod-carre}. We assume that $|\omega |>1+|\eta |$ and $\omega^{2}-\eta^{2}>2|\omega |$. We also assume that the hypotheses \ref{hypothesis2}, \ref{hypothesis4}, \ref{hypothesisG}, \ref{hypothesisH} and \ref{hypothesisI} (at least with $n_{1}=0$) are satisfied (see Assumptions \ref{Asummption23} and \ref{Asummption24}). 

Then, for every positive element $\zeta\in k^{\times}$, the image of $\Pi_{\CbeC^{+}_{\zeta}}$ is generated by the images of the $2$-torsion points of $\textrm{Jac}(\CbeC^{+}_{\zeta})(k(x)).$
\end{sprop}
\dem%
The idea of the proof is to apply Propositions \ref{caracterisationimagePi}, \ref{etude_y2-A_en_p}, \ref{J++-mod-B-C} and \ref{J++-mod-1-C} with $E(x):=\frac{1-C(x)}{2}$ and $D=\frac{(1+C(x))^{2}-4B(x)}{4}.$ Notice that under our choice for $E$ and $D$ we have
$$1-E=\frac{1+C}{2}\textrm{, }1-2E=C\textrm{, }E^{2}-D=B-C\textrm{ and }(1-E)^2-D=B.$$

Let $\beta$ be a $k(x)$-point of $\textrm{Jac}(\CbeC^{+}_{\zeta})$. We use the notation \ref{nota-sect-44-gene1}, \ref{nota-sect-44-gene2} and \ref{nota-sect-44-gene3} (relative to $\Pi_{\CbeC^{+}_{\zeta}}(\beta)$). 
The assumptions made on $E$ and $D$ while stating those notation 
are satisfied. In fact $\delta =\zeta\in k^{\times}$ is nontrivial and:
\begin{enumerate}
        \item Hypothesis \ref{hypothesis2} states that the remainder of the Euclidean division of $(1-E)^{2}-D=B$ by $2E=1-C$ is different from zero; hence $(1-E)^{2}-D$ and $E$ are coprime;
        \item in the same way 
        \begin{enumerate}
                \item $(1-E)^{2}-D$ and $1-2E$ are coprime (see Hypothesis \ref{hypothesis4}); 
                \item the inequality $\omega^{2}-\eta^{2}>2|\omega |$ implies the coprimality of $D$ and $E$ and the coprimality $D$ and $1-2E=C$;
        \end{enumerate}
        \item Since $\omega^{2}>(1+|\eta |)^{2}$, the degree of the polynomial \linebreak$4D=4\left(\omega^{2}-\eta^{2}\right)\left( x+b_{1}\right) +\left(\omega^{2}-\eta^{2}\right)^{2}+4\eta^{2}$ is $1$; thus $D$ is not a square in $k(x)$;
\end{enumerate}
The polynomials $\mu_{1,3}$, $\mu_{1,4}$, $\mu_{2,3}$ and $\mu_{3,4}$ defined in the notation \ref{nota-sect-44-gene3} help us to understand $\alpha=\Xi_{\CbeC^{+}_{\zeta}}(\beta )$ and $\Pi_{\CbeC^{+}_{\zeta}}(\beta )$:
$$\begin{array}{l}
\alpha =([\mu_{1,3}\mu_{1,4}],[\mu_{2,3}],[\mu_{1,3}\mu_{2,3}\mu_{3,4}],[\mu_{1,4}\mu_{3,4}])\textrm{ and}\\
\Pi_{\CbeC^{+}_{\zeta}}(\beta ) =\left(\left[\mu_{1,3}\mu_{1,4}\right] ,\left[\mu_{1,3}\mu_{3,4}\right] ,\left[\mu_{1,4}\mu_{3,4}\right]\right) .\\
\end{array}$$

The polynomials $1-2E=C$ and $2E=1-C$ are irreducible, 
the polynomial $(1-E)^{2}-D=B$ is equal to $\left( x+b_{1}+\eta\right)\left( x+b_{1}-\eta\right)$ and %
the polynomial $E^{2}-D=B-C$ is equal to $\left( x+b_{1}-1+\omega\right)\left( x+b_{1}-1-\omega\right) .$ Moreover
$$\begin{array}{l}
\Xi_{\CbeC^{+}_{\zeta}}(<y-\zeta E,0>)=\left(\left[\zeta\right] ,\left[ 2E\left( E^{2}-D\right)\right] ,\left[ 2\zeta E\right] ,\left[ E^{2}-D\right]\right) ,\\
\Xi_{\CbeC^{+}_{\zeta}}(<y^{2}-\zeta^{2}E^{2},0>)=\left(\left[ 1-2E\right] ,\left[ -\zeta \left( E^{2}-D\right)\right] ,\left[ -\zeta \left( 1-2E\right)\left( E^{2}-D\right)\right],\left[ 1\right]\right)\\
\Xi_{\CbeC^{+}_{\zeta}}(<y^{2}-\zeta^{2}D,0>)=\left( \left[\left( 1-E\right)^{2}-D\right] ,\left[ E^{2}-D\right] ,\left[ E^{2}-D\right] ,\left[\left( 1-E\right)^{2}-D\right]\right)\\
\end{array}$$
are in the image of the $2$-torsion by $\Xi_{\CbeC^{+}_{\zeta}}$. Thus we can assume without loss of generality that 
\begin{itemize}
\item[* ] $\mu_{1,3}$ is a constant $\epsilon_{1,3}\in k^{\times}$,
\item[* ] $\mu_{1,4}=\epsilon_{1,4}\left( x+b_{1}+\eta\right)^{n_{1}}$ with $\epsilon_{1,4}\in k^{\times}$ and $n_{1}\in\{ 0,1\}$, and
\item[* ] $\mu_{3,4}=\epsilon_{3,4}\left( x+b_{1}-1+\omega\right)^{n_{2}}$ with $\epsilon_{3,4}\in k^{\times}$ and $n_{2}\in\{ 0,1\}$.
\end{itemize} 
\vspace{0.3cm}

\noindent The next step is to apply Proposition \ref{J++-mod-B-C}. Its assumptions are satisfied: $E^{2}-D$ and $1-2E$ are coprime (since $(1-E)^{2}-D$ and $1-2E$ are coprime), $E^{2}-D$ and $E$ are coprime (since $E$ and $D$ are coprime), $E^{2}-D=B-C=(x+b_{1}-1)^{2}-\omega^{2}$ is squarefree (since $\omega\neq 0$), and $\delta$ and $E^{2}-D$ are coprime (since $\delta\in k^{\times}$). %
Proposition \ref{J++-mod-B-C} asserts the existence of two integers $n_{3},$ $n_{5}\in\{ 0,1\}$ such that 
\begin{equation}\label{equiv-J++44-verification-alpha1-cas-impair-zeta}
\begin{array}{l}
\epsilon_{1,3}\epsilon_{1,4}\left( x+b_{1}+\eta\right)^{n_{1}}\sim\zeta^{n_{2}}(1-2E)^{n_{3}}\bmod\left( x+b_{1}-1+\omega\right)\textrm{ and}\\
\epsilon_{1,3}\epsilon_{1,4}\left( x+b_{1}+\eta\right)^{n_{1}}\sim (1-2E)^{n_{5}}\bmod\left( x+b_{1}-1-\omega\right) .\\
\end{array}
\end{equation}
In other words we have
\begin{equation}\label{equiv-J++44-verification-alpha1-cas-impair-reduites-zeta}
\begin{array}{l}
\epsilon_{1,3}\epsilon_{1,4}\left( 1-\omega +\eta\right)^{n_{1}}\sim\zeta^{n_{2}}\left(\left(\omega -1\right)^{2}-\eta^{2}\right)^{n_{3}}\textrm{ and}\\
\epsilon_{1,3}\epsilon_{1,4}\left( 1+\omega +\eta\right)^{n_{1}}\sim\left(\left(\omega +1\right)^{2}-\eta^{2}\right)^{n_{5}}.\\
\end{array}
\end{equation}
Taking the product of these two equivalences we obtain
\begin{equation}\label{equiv-J++44-verification-alpha1-cas-impair-final-zeta}
\left(\left(\eta +1\right)^{2}-\omega^{2}\right)^{n_{1}}\sim\zeta^{n_{2}}\left(\left(\omega -1\right)^{2}-\eta^{2}\right)^{n_{3}}\left(\left(\omega +1\right)^{2}-\eta^{2}\right)^{n_{5}}.
\end{equation}
Two elements equivalent (relatively to $\sim$) have the same sign. Since $\zeta$, $\left(\omega -1\right)^{2}-\eta^{2}$ and $\left(\omega +1\right)^{2}-\eta^{2}$ are positive, and $\left(\eta +1\right)^{2}-\omega^{2}$ is negative, 
the equi\-va\-lence \ref{equiv-J++44-verification-alpha1-cas-impair-final-zeta} implies the parity of $n_{1}$. As a consequence and because of Equivalences \ref{equiv-J++44-verification-alpha1-cas-impair-reduites-zeta} 
the coefficient $\epsilon_{1,3}\epsilon_{1,4}$ is positive. \\

\noindent The third step is to apply Proposition \ref{J++-mod-1-C}. Its assumptions are satisfied: $E=\frac{1-C}{2}$ is squarefree (its degree is $1$), $E$ and $D$ are coprime, $\delta$ and $E$ are coprime (since $\delta\in k^{\times}$). %
Using the congruence $-\delta D\equiv \zeta (E^{2}-D)\bmod E$ we deduce from Proposition \ref{J++-mod-1-C} that
\begin{equation}\label{Eqt-avt-cas1-mu34cst-verifC+zeta-Appli-144}
\textrm{either }\mu_{1,3}\mu_{3,4}\sim 1\bmod E\textrm{ or }\mu_{1,3}\mu_{3,4}\sim \zeta (E^{2}-D)\bmod E.
\end{equation}
\begin{description}
\item[Case $1$: if $\mu_{3,4}$ is a constant. ] Then $n_{2}=0$ and the equi\-va\-lence \ref{equiv-J++44-verification-alpha1-cas-impair-final-zeta} becomes
$$1\sim\left(\left(\omega -1\right)^{2}-\eta^{2}\right)^{n_{3}}\left(\left(\omega +1\right)^{2}-\eta^{2}\right)^{n_{5}}.$$
Thus Hypothesis \ref{hypothesisI} implies the triviality of $n_{3}$ and $n_{5}$. In particular Equivalences \ref{equiv-J++44-verification-alpha1-cas-impair-reduites-zeta} becomes $\epsilon_{1,3},\epsilon_{1,4}\sim 1$.

We apply Proposition \ref{etude_y2-A_en_p} with $\mathcal{P}$ the infinite place of $k(x)$ and \linebreak$A:=\delta^{2}D$ to study $\Xi_{\CbeC^{+}_{\zeta},4}(\beta )=[\mu_{1,4}\mu_{3,4}]$. 
As $\mu_{1,4}$ and $\mu_{3,4}$ are constants ($n_{1}=n_{2}=0$), we obtain an equivalence $\mu_{1,4}\mu_{3,4}=\epsilon_{1,4}\epsilon_{3,4}\sim 1$ i.e. the triviality of the classe $\alpha=([\mu_{1,3}\mu_{1,4}],[\mu_{1,3}\mu_{3,4}],[\mu_{1,4}\mu_{3,4}])$.


\item[Case $2$: if $\mu_{3,4}=\epsilon_{3,4}\left( x+b_{1}-1+\omega\right)$. ] We apply Proposition \ref{etude_y2-A_en_p} with $\mathcal{P}$ the infinite place of $k(x)$ and $A:=\delta^{2}D$. As $\mu_{1,4}$ is a constant and $D$ is a de degree $1$ polynomial with leading coefficient $4\left(\omega^{2}-\eta^{2}\right)$, we obtain an equivalence 
$\mu_{1,4}\mu_{3,4}\sim -(\omega^{2}-\eta^{2})(x+b_{1}-1+\omega)$ i.e. $\epsilon_{1,4}\epsilon_{3,4}\sim -(\omega^{2}-\eta^{2})$. 

Using the value of $\mu_{3,4}$ Assertion \ref{Eqt-avt-cas1-mu34cst-verifC+zeta-Appli-144} becomes 
$$
\epsilon_{1,3}\epsilon_{3,4}\sim -2\left(\omega^{2}-\eta^{2}-2\omega\right)\textrm{ or }
\epsilon_{1,3}\epsilon_{3,4}\sim -2\zeta\left(\omega^{2}-\eta^{2}+2\omega\right).
$$
We reformulate these equivalences thanks to Equivalences \ref{equiv-J++44-verification-alpha1-cas-impair-reduites-zeta}, and $\epsilon_{1,4}\epsilon_{3,4}\sim -(\omega^{2}-\eta^{2})$, and $\left(\epsilon_{1,3}\epsilon_{1,4}\right)\left(\epsilon_{1,3}\epsilon_{3,4}\right)\left(\epsilon_{1,4}\epsilon_{3,4}\right)\sim 1.$ We get
$$\begin{array}{l}
2\left(\omega^{2}-\eta^{2}\right)\left(\omega^{2}-\eta^{2}-2\omega\right)\left(\left(\omega +1\right)^{2}-\eta^{2}\right)^{n_{5}}\sim 1\textrm{ or}\\
2\left(\omega^{2}-\eta^{2}\right)\left(\omega^{2}-\eta^{2}+2\omega\right)\left(\left(\omega -1\right)^{2}-\eta^{2}\right)^{n_{3}}\sim 1.\\
\end{array}$$
Both equivalences are in contradiction with either hypothesis \ref{hypothesisG} or hypothesis \ref{hypothesisH}: Case $2$ never happens.
\fin%
\end{description}

\begin{sprop}\label{verification-J++1x} We use the notation \ref{Notations_Principales_21}. We also use the notation \ref{nota_genre2} and \ref{nota-equiv-mod-carre}. 
%
%
We assume that $\omega > |\eta |+1$, $\omega^{2}-\eta^{2}>2\omega $ and $b_{1}>1+\frac{\omega^{2}-\eta^{2}}{2}$. We also assume that the hypotheses \ref{hypothesis2}, \ref{hypothesis4}, \ref{hypothesisI}, \ref{hypothesisJ}, \ref{hypothesisK} and \ref{hypothesisL} are satisfied. 

Then, for every $\zeta\in k^{\times}$ greater than 0, the image of $\Pi_{\CbeC^{+}_{\zeta x}}$ is generated by the images of the $2$-torsion points of $\textrm{Jac}(\CbeC^{+}_{\zeta x})(k(x)).$
\end{sprop}
\dem%
We apply Propositions \ref{caracterisationimagePi}, \ref{etude_y2-A_en_p}, \ref{J++-mod-B-C} and \ref{J++-mod-1-C} with $\delta :=\zeta x$, $E(x):=\frac{1-C(x)}{2}$ and $D=\frac{(1+C(x))^{2}-4B(x)}{4}.$

Let $\beta$ be a $k(x)$-point of $\textrm{Jac}(\CbeC^{+}_{\zeta x})$. We use the notation \ref{nota-sect-44-gene1}, \ref{nota-sect-44-gene2} and \ref{nota-sect-44-gene3} (relative to $\Pi_{\CbeC^{+}_{\zeta x}}(\beta )$). 
As in the proof of Proposition \ref{verification-J++1} 
\begin{itemize}
        \item[* ] the assumptions made on $E$ and $D$ while stating the notation \ref{nota-sect-44-gene1}, \ref{nota-sect-44-gene2} and \ref{nota-sect-44-gene3} 
are satisfied;
        \item[* ] by definition of $\mu_{1,3}$, $\mu_{1,4}$, $\mu_{2,3}$ and $\mu_{3,4}$ we have 
$$\begin{array}{l}
\Xi_{\CbeC^{+}_{\zeta x}}(\beta )=([\mu_{1,3}\mu_{1,4}],[\mu_{2,3}],[\mu_{1,3}\mu_{2,3}\mu_{3,4}],[\mu_{1,4}\mu_{3,4}])\textrm{ and}\\
\Pi_{\CbeC^{+}_{\zeta x}}(\beta ) =\left(\left[\mu_{1,3}\mu_{1,4}\right] ,\left[\mu_{1,3}\mu_{3,4}\right] ,\left[\mu_{1,4}\mu_{3,4}\right]\right) .\\
\end{array}$$
\end{itemize}
We can eventually substitute $x\mu_{1,3}$, $x\mu_{1,4}$ and $x\mu_{3,4}$ to $\mu_{1,3}$, $\mu_{1,4}$ and $\mu_{3,4}$. Thus, without loss of generality, we can assume that $\mu_{1,4}$ and $x$ are coprime.

The polynomials $1-2E=C$ and $2E=1-C$ are irreducible, 
the polynomial $(1-E)^{2}-D=B$ is equal to $\left( x+b_{1}+\eta\right)\left( x+b_{1}-\eta\right)$ and %
the polynomial $E^{2}-D=B-C$ is equal to $\left( x+b_{1}-1+\omega\right)\left( x+b_{1}-1-\omega\right) .$ Moreover
$$\begin{array}{l}
\Xi_{\CbeC^{+}_{\zeta x}}(<y-\zeta xE,0>)=\left(\left[\zeta x\right] ,\left[ 2E\left( E^{2}-D\right)\right] ,\left[ 2\zeta xE\right] ,\left[ E^{2}-D\right]\right) ,\\
\Xi_{\CbeC^{+}_{\zeta x}}(<y^{2}-\zeta^{2}x^{2}E^{2},0>)=\left(\left[ 1-2E\right] ,\left[ -\zeta x\left( E^{2}-D\right)\right] ,\left[ -\zeta x\left( 1-2E\right)\left( E^{2}-D\right)\right],\left[ 1\right]\right)\\
\Xi_{\CbeC^{+}_{\zeta x}}(<y^{2}-\zeta^{2}x^{2}D,0>)=\left( \left[\left( 1-E\right)^{2}-D\right] ,\left[ E^{2}-D\right] ,\left[ E^{2}-D\right] ,\left[\left( 1-E\right)^{2}-D\right]\right)\\
\end{array}$$
are in the image of the $2$-torsion by $\Xi_{\CbeC^{+}_{\zeta x}}$. Thus we can also assume without loss of generality that
\begin{itemize}
\item[* ] $\mu_{1,3}$ is a constant $\epsilon_{1,3}\in k^{\times}$,
\item[* ] $\mu_{1,4}=\epsilon_{1,4}\left( x+b_{1}+\eta\right)^{n_{1}}$ with $\epsilon_{1,4}\in k^{\times}$ and $n_{1}\in\{ 0,1\}$, and
\item[* ] $\mu_{3,4}=\epsilon_{3,4}x^{n_{6}}\left( x+b_{1}-1+\omega\right)^{n_{2}}\left( x+b_{1}-1-\omega\right)^{n_{4}}$ with $\epsilon_{3,4}\in k^{\times}$ and $n_{2},n_{4},n_{6}\in\{ 0,1\}$.
\end{itemize} 
Since $E^{2}(0)-D(0)=(b_{1}-1)^{2}-\omega^{2}$ is greater than $0$ the polynomials $\delta =\zeta x$ and $E^{2}-D$ are coprime. Hence, as in the proof of Proposition \ref{verification-J++1}, we can apply Proposition \ref{J++-mod-B-C} : there are two integers $n_{3},$ $n_{5}\in\{ 0,1\}$ such that 
\begin{equation}\label{equiv-J++44-verification-alpha1-cas-impair-final-valeurE13E14}
\begin{array}{l}
\epsilon_{1,3}\epsilon_{1,4}\left( 1-\omega +\eta\right)^{n_{1}}\sim\zeta^{n_{2}}\left( 1-\omega-b_{1}\right)^{n_{2}}\left(\left(\omega -1\right)^{2}-\eta^{2}\right)^{n_{3}}\textrm{ and}\\
\epsilon_{1,3}\epsilon_{1,4}\left( 1+\omega +\eta\right)^{n_{1}}\sim\zeta^{n_{4}}\left( 1+\omega-b_{1}\right)^{n_{4}}\left(\left(\omega +1\right)^{2}-\eta^{2}\right)^{n_{5}}.\\
\end{array}
\end{equation}
Taking the product of these two equivalences we obtain
\begin{equation}\label{equiv-J++44-verification-alpha1-cas-impair-final}
\begin{array}{rcl}
\left(\left(\eta +1\right)^{2}-\omega^{2}\right)^{n_{1}}&\sim &\zeta^{n_{2}+n_{4}}\left( 1-\omega -b_{1}\right)^{n_{2}}\left( 1+\omega-b_{1}\right)^{n_{4}}\\
&&\times\left(\left(\omega -1\right)^{2}-\eta^{2}\right)^{n_{3}}\left(\left(\omega +1\right)^{2}-\eta^{2}\right)^{n_{5}}.\\
\end{array}
\end{equation}
Two elements equivalent (relatively to $\sim$) have the same sign. Moreover
\begin{itemize}
\item[* ] $\zeta$, $\left(\omega -1\right)^{2}-\eta^{2}$ and $\left(\omega +1\right)^{2}-\eta^{2}$ are positives ;
\item[* ] $\left(\eta +1\right)^{2}-\omega^{2}$, $1-\omega -b_{1}=-\frac{\omega^{2}-\eta^{2}+2\omega}{2}-\left( b_{1}-1-\frac{\omega^{2}-\eta^{2}}{2}\right)$ and \linebreak$1+\omega -b_{1}=-\frac{\omega^{2}-\eta^{2}-2\omega}{2}-\left( b_{1}-1-\frac{\omega^{2}-\eta^{2}}{2}\right)$ are negatives.
\end{itemize}
Thus the equi\-va\-lence \ref{equiv-J++44-verification-alpha1-cas-impair-final} implies the parity of $n_{1}+n_{2}+n_{4}$.

In the next steps we use Proposition \ref{J++-mod-1-C}. Its assumption are satisfied: $E=\frac{1-C}{2}$ is squarefree (its degree is $1$), $E$ and $D$ are coprime, and $\delta =\zeta x$ and $E$ are coprime (because $E(0)=1-b_{1}-\frac{\omega^{2}-\eta^{2}}{2}<0$).

\begin{description}
        \item[Case $(1)$: if $\mu_{1,4}$ is constant. ] Then $n_{1}=0$, and $n_{2}$ and $n_{4}$ have the same parity. Hence Equivalence \ref{equiv-J++44-verification-alpha1-cas-impair-final} becomes
$$1\sim\left( (b_{1}-1)^{2}-\omega^{2}\right)^{n_{2}}\left(\left(\omega -1\right)^{2}-\eta^{2}\right)^{n_{3}}\left(\left(\omega +1\right)^{2}-\eta^{2}\right)^{n_{5}}.$$
This equivalence and hypothesis \ref{hypothesisI} implies the triviality of $n_{2}$, $n_{4}$, $n_{3}$ and $n_{5}$. In particular the polynomial $\mu_{3,4}$ is equal to $\epsilon_{3,4}x^{n_{6}}$.
        \begin{description}
                \item[Subcase $(1a)$: if $\mu_{3,4}$ is constant. ] Since $n_{1}=n_{2}=n_{3}=0$, Equivalence \ref{equiv-J++44-verification-alpha1-cas-impair-final-valeurE13E14} becomes $\mu_{1,3}\mu_{1,4}=\epsilon_{1,3}\epsilon_{1,4}\sim1$.

As in the proof of Proposition \ref{verification-J++1} we apply Proposition \ref{etude_y2-A_en_p} with $\mathcal{P}$ the infinite place of $k(x)$ and $A:=\delta^{2}D$. 
Since $\mu_{1,4}$ and $\mu_{3,4}$ are constants ($n_{2}=n_{3}=0$), we obtain an equivalence $\mu_{1,4}\mu_{3,4}=\epsilon_{1,4}\epsilon_{3,4}\sim 1$ i.e. the triviality of the class $\Xi_{\CbeC^{+}_{\zeta x}}(\beta ) =([\mu_{1,3}\mu_{1,4}],[\mu_{1,3}\mu_{3,4}],[\mu_{1,4}\mu_{3,4}])$.
                \item[Subcase $(1b)$: if $\mu_{3,4}=\epsilon_{3,4}x$. ] We apply Proposition \ref{etude_y2-A_en_p} with $\mathcal{P}$ the infinite place of $k(x)$ and $A:=\delta^{2}D$. Since the degree of $\mu_{1,4}\mu_{3,4}$ is $1$, we obtain an equivalence $\epsilon_{1,4}\epsilon_{3,4}\sim -(\omega^{2}-\eta^{2})$.

From Proposition \ref{J++-mod-1-C} we have either $\mu_{1,3}\mu_{3,4}\sim 1\bmod E$ or $\mu_{1,3}\mu_{3,4}\sim \zeta x(E^{2}-D)\bmod E$. These equivalences can be reformulate in the form
$$-(\omega^{2}-\eta^{2})x\sim 1\bmod E\textrm{ or }-(\omega^{2}-\eta^{2})(E^{2}-D)\sim \zeta\bmod E$$
because $\epsilon_{1,4}\epsilon_{3,4}\sim -(\omega^{2}-\eta^{2})$ and $\epsilon_{1,3}\epsilon_{1,4}\sim 1$ (see Equivalence \ref{equiv-J++44-verification-alpha1-cas-impair-final-valeurE13E14} and use equalities $n_{1}=n_{2}=n_{3}=0$).

Using the equality $-2E=C-1=2(x+b_{1}-1)+\omega^{2}-\eta^{2}$, the first equivalence becomes $2\left(\omega^{2}-\eta^{2}\right)\left(2b_{1}-2+\omega^{2}-\eta^{2}\right)\sim 1.$ Hence this equivalence contradicts hypothesis \ref{hypothesisJ}.

The remainder of the Euclidean division of $E^{2}-D=B-C=(x+b_{1}-1)^{2}-\omega^{2}$ by $-2E=C-1=2(x+b_{1}-1)+\omega^{2}-\eta^{2}$ is $\frac{(\omega^{2}-\eta^{2})^{2}-4\omega^{2}}{4}$. This remainder and $\omega^{2}-\eta^{2}$ being positive and two equivalent elements (relatively to $\sim$) having the same sign, the equivalence $-(\omega^{2}-\eta^{2})(E^{2}-D)\sim \zeta\bmod E$ contradicts the positivity of $\zeta$.
        \end{description}
        \item[Case $(2)$: if $\mu_{1,4}=\epsilon_{1,4}(x+b_{1}+\eta ).$ ] Then $n_{2}$ and $n_{4}$ have different parities.
        \begin{description}
                \item[Subcase $(2a)$: if $n_{2}$ is odd. ] Then $n_{4}$ is trivial. We apply Proposition \ref{J++-mod-1-C}. As $\mu_{1,3}\mu_{3,4}=\epsilon_{1,3}\epsilon_{3,4}x^{n_{6}}(x+b_{1}-1+\omega )$ and $E^{2}-D=B-C=\left( x+b_{1}-1+\omega\right)\left( x+b_{1}-1-\omega\right)$ we have either
\begin{equation}\label{equiv-J++44-verification-alpha4-cas-delta-avluation-inconnue-E2-D-valuationimpaire-n2impair}
\begin{array}{l}
\epsilon_{1,3}\epsilon_{3,4}x^{n_{6}}(x+b_{1}-1+\omega )\sim 1\bmod E\textrm{ or }\\
\epsilon_{1,3}\epsilon_{3,4}x^{1-n_{6}}(x+b_{1}-1-\omega )\sim \zeta\bmod E.\\
\end{array}
\end{equation}
Since the degree of $\mu_{1,4}\mu_{3,4}$ is $2+n_{6}$, Proposition \ref{etude_y2-A_en_p} (applied with $\mathcal{P}$ the infinite place of $k(x)$ and $A:=\delta^{2}D$) gives us an equivalence $\epsilon_{1,4}\epsilon_{3,4}\sim (-(\omega^{2}-\eta^{2}))^{n_{6}}$. In particular, the signs of $\epsilon_{1,4}\epsilon_{3,4}$ and $(-1)^{n_{6}}$ are equal. Hence $\epsilon_{1,3}\epsilon_{1,4}$ being positive (see the second of the two equivalences \ref{equiv-J++44-verification-alpha1-cas-impair-final-valeurE13E14}), the signs of $\epsilon_{1,3}\epsilon_{3,4}$ and $(-1)^{n_{6}}$ are equal. 

The remainders $-(b_{1}-1+\frac{\omega^{2}-\eta^{2}}{2})$ and $-\frac{\omega^{2}-\eta^{2}-2\omega}{2}$ of the Euclidean divisions of $x$ and $x+b_{1}-1+\omega$ by $-2E=2(x+b_{1}-1)+\omega^{2}-\eta^{2}$ are negative. Thus a sign argument proves that only one of the two equivalences \ref{equiv-J++44-verification-alpha4-cas-delta-avluation-inconnue-E2-D-valuationimpaire-n2impair}: the second one i.e. the equivalence
$$(-2)^{n_{6}}\epsilon_{1,3}\epsilon_{3,4}\left( 2b_{1}-2+\omega^{2}-\eta^{2}\right)^{1-n_{6}}\left(\omega^{2}-\eta^{2}+2\omega\right)\sim\zeta .$$

From this equivalence and the equivalences %
$\epsilon_{1,4}\epsilon_{3,4}\sim\left( -(\omega^{2}-\eta^{2})\right)^{n_{6}}$ and 
$\epsilon_{1,3}\epsilon_{1,4}\sim\zeta\left( b_{1}-1+\omega\right)\left(\omega -1-\eta\right)^{1-n_{3}}\left(\omega -1+\eta\right)^{n_{3}}$ (see\linebreak Equivalences \ref{equiv-J++44-verification-alpha1-cas-impair-final-valeurE13E14}) 
we deduce an equivalence in contradiction with hypothesis \ref{hypothesisK}:
$$\begin{array}{rcl}1&\sim &2^{n_{6}}\left(\omega^{2}-\eta^{2}+2\omega\right)\left( b_{1}-1+\omega\right)\left(\omega^{2}-\eta^{2}\right)^{n_{6}}\times\\
&&\left( 2b_{1}-2+\omega^{2}-\eta^{2}\right)^{1-n_{6}}\left(\omega -1-\eta\right)^{1-n_{3}}\left(\omega -1+\eta\right)^{n_{3}},\\
\end{array}$$
                \item[Subcase $(2b)$: if $n_{2}$ is even. ] The proof is analogous to the proof of Subcase $(2a)$; we get a contradiction with Hypothesis \ref{hypothesisL}.
\fin%
        \end{description}
\end{description}

\end{subsection}

\begin{subsection}{The image of $\Pi_{\widehat{\CbeC}^{+}_{\delta}}.$}
\begin{snota}\label{nota-sect-45-gene1} Let $k$ be a characteristic $0$ field. Let $B,$ $C,$ $\delta\in k[x]$ be three nontrivial polynomials such that $1-C,$ $B-C$ and $(1+C)^{2}-4B$ are nontrivial. We consider the hyperelliptic curve $\CbeH$ defined on $k(x)$ by the affine equation 
$$\CbeH :z^{2}=(y+\delta (1+C))(y^{2}-4\delta^{2}B)(y^{2}-4\delta^{2}C).$$
\end{snota}
\rmk%
As the polynomials $B,$ $C,$ $\delta$, $1-C$, $B-C$ and $(1+C)^{2}-4B$ are nontrivial, the polynomial $(y+\delta (1+C))(y^{2}-4\delta^{2}B)(y^{2}-4\delta^{2}C)$ is squarefree.

\begin{snota}\label{nota-sect-45-gene2} We use the notation \ref{nota-sect-45-gene1}. Consider the three polynomials $f_{1}(y):=y+\delta (1+C)$, $f_{2}(y):=y^{2}-4\delta^{2}B$ and $f_{3}(y):=y^{2}-4\delta^{2}C$. We assume that $f_{2}(y)$ and $f_{3}(y)$ are irreducible (i.e. that $B$ and $C$ are not squares in $k(x)$). 

Denote by $K_{i}$ the field $k(x)[T]/(f_{i}(T))$ and by $y_{i}$ the class of $T$ in $K_{i}$. Let $\pi_{\CbeH}:\textrm{Jac}(\CbeH )(k(x))\longrightarrow\displaystyle\prod_{i=1}^{4}K_{i}^{\times}/K_{i}^{\times 2}$ be the morphism defined by proposition \ref{def-Cas-Sch+im-in-ker} and $\pi_{\CbeH ,i}:\textrm{Jac}(\CbeH )(k(x))\longrightarrow\displaystyle K_{i}^{\times}/K_{i}^{\times 2}$ be its $i$-th component. In this section we study the morphism $\Xi_{\CbeH}:\textrm{Jac}(\CbeH )(k(x))\longrightarrow\displaystyle\prod_{i=1}^{4}k(x)^{\times}/k(x)^{\times 2}$ with $i$-th component the morphism $\Xi_{\CbeH,i}:=N_{K_{i}/k(x)}\circ\pi_{\CbeH ,i}$.
\end{snota}

\begin{snota}\label{nota-sect-45-gene3} We use the notation \ref{nota-sect-45-gene1} and \ref{nota-sect-45-gene2}. We assume that %
$B$ and $C$ are coprime, %
$B$ and $1-C$ are coprime, %
$(1+C)^{2}-4B$ and $C$ are coprime, and %
$B-1$ and $C-1$ are coprime.

Let $\alpha$ be an element of the image $\textrm{Im}(\Xi_{\CbeH})$. From proposition \ref{caracterisationimagePi}, we know that $\alpha =([\mu_{1,2}\mu_{1,3}],[\mu_{1,2}\mu_{2,3}],[\mu_{1,3}\mu_{2,3}])$ with
\begin{itemize}
\item[* ] $\mu_{1,2}\in k[x]$ a squarefree divisor of $\delta ((1+C)^{2}-4B),$
\item[* ] $\mu_{1,3}\in k[x]$ a squarefree divisor of $\delta (1-C),$ and
\item[* ] $\mu_{2,3}\in k[x]$ a squarefree divisor of $\delta (B-C)$.
\end{itemize}
\end{snota}

\begin{sprop}\label{etude_mu_Jchap}%
%
We use the notation \ref{nota-sect-45-gene1}, \ref{nota-sect-45-gene2} and \ref{nota-sect-45-gene3}. Let $p$ be a prime factor of $B-C$ such that\begin{itemize}
\item[* ] $v_{p}(\delta )=v_{p}(1-C)=v_{p}(C)=0$, and 
\item[* ] there is no $\lambda\in k$ such that $C\equiv\lambda^{2}\bmod p$.
\end{itemize}
Then for every $\alpha$ in the image of $\Xi_{\CbeH}$ there exist three polynomials $\alpha_{1},$ $\alpha_{2}$, $\alpha_{3}\in k[x]$ such that $\alpha =\left(\left[\alpha_{1}\right] ,\left[\alpha_{2}\right] ,\left[\alpha_{3}\right]\right)$, and $v_{p}(\alpha_{1})=0$ and
$$\alpha_{1}\sim 1\bmod p.$$
\end{sprop}
\dem%
Let $<u,v>\in\textrm{Jac}(\CbeH )(k(x))$ be a $k(x)$-point of $\textrm{Jac}(\CbeH )$. Without loss of generality, we can assume that $u$ and $f_{1}f_{2}f_{3}$ are coprime (eventually we need to add a $2$-torsion point to $<u,v>$).

Write $u(y)=\frac{u_{2}}{\lambda}y^{2}+\frac{u_{1}}{\lambda}y+\frac{u_{0}}{\lambda}$ with $u_{0},$ $u_{1},$ $u_{2}$, $\lambda\in k[x]$ four coprime polynomials. The polynomial $u$ being monic, $\lambda$ is the leading coefficient of $\widetilde{u}(y):=u_{2}y^{2}+u_{1}y+u_{0}.$ In particular $u_{0},$ $u_{1}$ and $u_{2}$ are coprime.

As in the proof for Proposition \ref{J++-mod-B-C} we look at the polynomial 
\begin{equation}\label{defn-psi_u(T)-4-5-11}
\Psi_{\widetilde{u}}(T):=\textrm{res}_{y}\left((-1)^{\deg (\widetilde{u})}\widetilde{u}(y),y^{2}-T\right) =\left( u_{0}+Tu_{2}\right)^{2}-Tu_{1}^{2}.
\end{equation}

\begin{description}
        \item[Case 1: if $p$ does not divide $\Psi_{\widetilde{u}}(4\delta^{2}C)$. ] Then $\Psi_{\widetilde{u}}\left(\delta^{2}C\right)$ has an inverse modulo $p$. As in the proof for Proposition \ref{J++-mod-B-C} we use Taylor's formula at $4\delta^{2}C$ to compute $\Psi_{\widetilde{u}}(4\delta^{2}B)$: since $p$ divides $4\delta^{2}B-4\delta^{2}C$, Taylor's formula implies the equivalence $\Psi_{\widetilde{u}}\left( \delta^{2}B\right)\sim\Psi_{\widetilde{u}}\left(\delta^{2}C\right)\bmod p$. This is enough to conclude: following Proposition \ref{def-Cas-Sch+im-in-ker} we have 
$$\begin{array}{rcl}
\Xi_{\CbeH ,1}(<u,v>) &=&\Xi_{\CbeH ,2}(<u,v>)\Xi_{\CbeH ,3}(<u,v>)\\
&=&\left[\Psi_{\widetilde{u}}\left(4\delta^{2}B\right)\Psi_{\widetilde{u}}\left(4\delta^{2}C\right)\right] .\\
\end{array}$$
        \item[Case 2: if $p$ divides $\Psi_{\widetilde{u}}(4\delta^{2}C)$. ] Then the congruence 
\begin{equation}\label{congruence-base-cas-2-ds-4-5-11}
\left( u_{0}+4\delta^{2}Cu_{2}\right)^{2}\equiv 4\delta^{2}Cu_{1}^{2}\bmod p
\end{equation} 
holds (see Definition \ref{defn-psi_u(T)-4-5-11}).
        \begin{description}
                \item[Subcase $(a)$: if $p$ and $u_{1}$ are coprime. ] Then $C$ is a square modulo $p$. This is in contradiction with our choice for $p$.
                \item[Subcase $(b)$: if $p$ divide $u_{1}$. ] Then Equation \ref{congruence-base-cas-2-ds-4-5-11} gives a congruence %
\begin{equation}\label{equation-etude_mu_Jchap45}
u_{0}\equiv -4\delta^{2}Cu_{2}\bmod p. %
\end{equation}
The polynomial $u_{0}$, $u_{1}$ and $u_{2}$ being corpime, $p$ does not divide $u_{2}$. In particular the degree of $\widetilde{u}$ is $2$ and its 
leading coefficient $\lambda$ is equal to $u_{2}$. Thus from the divisibility of $u_{1}$ by $p$ and from Equation \ref{equation-etude_mu_Jchap45} we get a congruence 
$$(-1)^{\deg (\widetilde{u})}\lambda\widetilde{u}(\delta (1+C))\equiv u_{2}^{2}\delta^{2}(1-C)^{2}\bmod p.$$
We conlude by noticing 
that $\Xi_{\CbeH ,1}(<u,v>)$ is equal to the class $\left[ (-1)^{\deg (\widetilde{u})}\lambda\widetilde{u}(\delta (1+C))\right] .$
\fin
\end{description}
\end{description}

\begin{sprop}\label{verification-Cchap+zeta} We use the notation of Theorem \ref{theo-final-pas-somme-3-carres-46}. %
For every \linebreak$\delta\in k (x)^{\times}$ we consider the hyperelliptic curve $\widehat{\CbeC}^{+}_{\delta}$ on $k(x)$ given by the affine equation:
$$\widehat{\CbeC}^{+}_{\delta}:z^{2}=(y+\delta (1+C))(y^{2}-4\delta^{2}B)(y^{2}-4\delta^{2}C).$$
We use the notation \ref{nota_genre2} and \ref{nota-equiv-mod-carre}. We assume that $\left(\omega^{2}-\eta^{2}\right)^{2}-4\eta^{2}$, \linebreak$B(0)-C(0)=\left( b_{1}-1\right)^{2}-\omega^{2}$ and $\eta$ are different from $0$. 


We also assume that the hypotheses \ref{hypothesis2}, \ref{hypothesis3}, \ref{hypothesis4}, \ref{hypothesis5}, \ref{hypothesis6}, \ref{hypothesisA}, \ref{hypothesisI} (at least in the two cases $n_{1}=n_{3}=0$ and $n_{1}=n_{2}=0$), \ref{hypothesisM} and \ref{hypothesisN} are satisfied.


\noindent Then for every $\zeta\in k$ greater than $0$
\begin{itemize}
\item[* ] $\textrm{Im}(\Pi_{\widehat{\CbeC}^{+}_{\zeta}})$ is the image of the $2$-torsion of $\textrm{Jac}(\widehat{\CbeC}^{+}_{\zeta})(k(x))$ by $\Pi_{\widehat{\CbeC}^{+}_{\zeta}}$ and
\item[* ] $\textrm{Im}(\Pi_{\widehat{\CbeC}^{+}_{\zeta x}})$ is the image of the $2$-torsion of $\textrm{Jac}(\widehat{\CbeC}^{+}_{\zeta x})(k(x))$ by $\Pi_{\widehat{\CbeC}^{+}_{\zeta x}}$.
\end{itemize}
\end{sprop}
\dem%
Let $\zeta$ be a positive element of $k^{\times}$. Let $\delta$ be either $\zeta$ or $\zeta x$. Let $\beta$ be a $k(x)$-point of $\textrm{Jac}(\widehat{\CbeC}^{+}_{\zeta})$. We use notation \ref{nota-sect-45-gene1}, \ref{nota-sect-45-gene2} and \ref{nota-sect-45-gene3} (relative to $\Pi_{\widehat{\CbeC}^{+}_{\delta}}(\beta )$). The assumptions made on $B$, $C$ and $\delta$ while stating those notation are satisfied: $\delta$ is nontrivial and
\begin{enumerate}
        \item $C$ is not a square in $k(x)$ (its degree is $1$);
        \item $B=(x+b_{1})^{2}-\eta^{2}$ is not a square in $k(x)$ (in fact it does not even have a double root since $\eta\neq 0$);
        \item the remainder of the Euclidean division of $4B$ by $C-1$ is the element $(\omega^{2}-\eta^{2}-2)^{2}-4\eta^{2}$; hence
$B$ and $C-1$ are coprime (see Hypothethis \ref{hypothesis2}) and %
$B-1$ and $C-1$ are coprime (see Hypothethis \ref{hypothesis3});
        \item the remainder of the Euclidean division of $4B$ by $C$ is the element $(\omega^{2}-\eta^{2}-1)^{2}-4\eta^{2}$; hence
the polynomials $B$ and $C$ are coprime (see Hypothethis \ref{hypothesis4}), and the polynomials %
$(1+C)^{2}-4B$ and $C$ are coprime (notice that $\pgcd ((1+C)^{2}-4B,C)=\pgcd (1-4B,C)$ and use Hypothethis \ref{hypothesis5}); in particular $B$ and $C$ are nontrivial and distincts;
\end{enumerate}
Proposition \ref{caracterisationimagePi} applies and asserts the existence of
\begin{itemize}
\item[* ] $\mu_{1,2}\in k[x]$ a squarefree divisor of $\delta ((1+C)^{2}-4B),$
\item[* ] $\mu_{1,3}\in k[x]$ a squarefree divisor of $\delta (1-C),$ and
\item[* ] $\mu_{2,3}\in k[x]$ a squarefree divisor of $\delta (B-C)$
\end{itemize}
such that $\Xi_{\widehat{\CbeC}^{+}_{\delta}}(\beta )=([\mu_{1,2}\mu_{1,3}],[\mu_{1,2}\mu_{2,3}],[\mu_{1,3}\mu_{2,3}])$.

The class $([(1+C)^{2}-4B],[(1+C)^{2}-4B],[1])$ is the image under $\Xi_{\widehat{\CbeC}^{+}_{\delta}}$ of the $2$-torsion point $<y+\delta (1+C),0>$. Moreover the degree of $(1-C)^{2}-4B$ is $1$. Thus we can assume without loss of generality that $\mu_{1,2}$ is a divisor of $\delta$ (eventually we add $<y+\delta (1+C),0>$ to $\beta$).

The remainder $(1+\omega )^{2}-\eta^{2}$ of the Euclidean division of $B=(x+b_{1})^{2}-\eta^{2}$ by $p:=x+b_{1}-1-\omega$ is not a square in $k$ (see hypothesis \ref{hypothesisI}). An application of Proposition \ref{etude_y2-A_en_p-cas-pair} with $A:=\delta^{2}B$ and $\mathcal{P}$ the place with uniformizer $p$ gives 
the coprimality of $\mu_{1,3}\mu_{2,3}$ and $p$. Since $p$ divides $B-C$, since $\mu_{1,3}$ divides $\delta (1-C)$, and since $B-C$ and $\delta (1-C)$ are coprime we get the coprimality of $p$ and $\mu_{2,3}$. In the same way, we can deduce from 
hypothesis \ref{hypothesisI} the coprimality of $\mu_{2,3}$ and $x+b_{1}-1+\omega$. 
This means that $\mu_{2,3}$ is a divisor of $\delta$.

Assume temporarily that $\delta$ is equal to $\zeta x$. From Hypothesis \ref{hypothesisM} we know that $B(0)=b_{1}^{2}-\eta^{2}$ is not a square in $k$. 
Applying Proposition \ref{etude_y2-A_en_p-cas-pair} with $A:=\delta^{2}B$ and $\mathcal{P}$ the place with uniformizer $x$ we show that 
$v_{x}(\mu_{1,2}\mu_{2,3})$ is even (notice that following the definition of $\mu_{1,2}$ and $\mu_{2,3}$, there is a polynomial $u\in k(x)[y]$ such that $\mu_{1,2}\mu_{2,3}\sim N_{K_{2}/k(x)}(u(y_{2}))$). In the same way hypothesis \ref{hypothesisN} allows us to apply Proposition \ref{etude_y2-A_en_p-cas-pair} with $A:=\delta^{2}C$ and $\mathcal{P}$ the place with uniformizer $x$. Doing so we prove that $v_{p}(\mu_{1,3}\mu_{2,3})$ is even. 

We do not assume $\delta$ is $\zeta x$ anymore. The valuations $v_{x}(\mu_{1,2})$, $v_{x}(\mu_{1,3})$ and $v_{x}(\mu_{2,3})$ have the same parity. Hence after eventually replacing $\mu_{i,j}$ by $x^{-1}\mu_{i,j}$ we can assume without loss of generality
\begin{itemize}
        \item[* ] the existence of $\epsilon_{1,2},\epsilon_{2,3}\in  k^{\times}$ such that $\mu_{1,2}=\epsilon_{1,2}$ and $\mu_{2,3}=\epsilon_{2,3}$,
        \item[* ] and the existence of $\epsilon\in  k^{\times}$ such that $\mu_{1,3}=\epsilon$ or $\mu_{1,3}=\epsilon (1-C)$ (the degree of $1-C$ is $1$).
\end{itemize}
Now we study $\Pi_{\widehat{\CbeC}^{+}_{\delta},3}(\beta )=[\mu_{1,3}\mu_{2,3}]$ by applying Proposition \ref{etude_y2-A_en_p} with $A:=4\delta^{2}C$ and $\mathcal{P}$ the infinite place of $k(x)$. Since $1-C$ and $-C$ have the same leading coefficient this proposition asserts that $\epsilon_{2,3}\epsilon\sim 1$.

The polynomial $\mu_{1,2}\mu_{2,3}$ is a constant. For every prime factor $p$ of $B$ the equivalence $\mu_{1,2}\mu_{2,3}\sim 1\bmod p$ holds (apply Proposition \ref{etude_y2-A_en_p} with $A:=4\delta^{2}B$ and $\mathcal{P}$ the place with uniformizer $p$; notice that $v_{p}(B)=1$ since $B$ is squarefree). In particular we have $\mu_{1,2}\sim\mu_{2,3}$. Since $\epsilon_{2,3}\sim\epsilon$, this means that that one of the equivalences $\mu_{1,2}\mu_{1,3}\sim 1$ or $\mu_{1,2}\mu_{1,3}\sim 1-C$ holds.

Assume that $\mu_{1,2}\mu_{1,3}\sim 1-C$. We know the coprimality of $B-C$ and $C$, and the coprimality of $B-C$ and $C-1$. Moreover $B-C$ and $x$ are coprime (by assumption $B(0)-C(0)=(b_{1}-1)^{2}-\omega^{2}$ is nontrivial). The polynomial $B-C$ factorizes as $B-C=(x+b_{1}-1-\omega )(x+b_{1}-1+\omega )$. The remainder of the Euclidean division of $C$ by $x+b_{1}-1-\omega$ is $(\omega +1)^{2}-\eta^{2}$. This is not a square in $k$ (see Hypothesis \ref{hypothesisI}). In the same way we prove that $C$  is not a square modulo $x+b_{1}-1+\omega$. Hence Proposition \ref{etude_mu_Jchap} applies: $\mu_{1,2}\mu_{1,3}\sim 1\bmod p$ for every prime divisor $p$ of $B-C$. In particular $\eta^{2}-\omega^{2}-2\omega$ and $\eta^{2}-\omega^{2}+2\omega$ are squares in $k$ (notice that $\mu_{1,2}\mu_{1,3}\sim 1-C$). This is in contradiction with hypothesis \ref{hypothesisA}. As a consequence $\mu_{1,2}\mu_{1,3}\sim 1$ and thus $\Xi_{\widehat{\CbeC}^{+}_{\delta}}(\beta )=([\mu_{1,2}\mu_{1,3}],[\mu_{1,2}\mu_{2,3}],[\mu_{1,3}\mu_{2,3}])$ is trivial.
\fin%

\end{subsection}


\end{section}

\begin{section}{A proof of Theorem \ref{theo-final-pas-somme-3-carres-46}}

Assumptions \ref{Asummption22} implies the positivity of $P(x^{2},y^{2})$. Let us show that $P(x^{2},y^{2})$ is not a sum of three squares in $\R (x,y)$.

We use the notation of Theorem \ref{theo-critere-nullite-rang}. Applying Proposition \ref{prop-verif-C-zeta41final}, Proposition \ref{prop-verif-C-zetax41final}, Proposition \ref{prop-verif-Cchap-zeta42final}, Proposition \ref{prop-verif-Cchap-zetax42final}, Proposition \ref{verification-J++1}, Proposition \ref{verification-J++1x} and Proposition \ref{verification-Cchap+zeta}, we know that for every positive element $\zeta\in k^{\times}$ the images of 
$$
\gamma_{\CbeC^{-}_{\zeta}},\textrm{ } \gamma_{\CbeC^{-}_{\zeta x}},\textrm{ } \gamma_{\widehat{\CbeC}^{-}_{\zeta}},\textrm{ } \gamma_{\widehat{\CbeC}^{-}_{\zeta x}},\textrm{ } \Pi_{\CbeC^{+}_{\zeta}},\textrm{ } \Pi_{\CbeC^{+}_{\zeta x}},\textrm{ } \Pi_{\widehat{\CbeC}^{+}_{\zeta}}\textrm{ and }\Pi_{\widehat{\CbeC}^{+}_{\zeta x}}
$$ 
are respectively the images of the $k(x)$-rational torsion subgroups of 
$$
\CbeC^{-}_{\zeta},\textrm{ } \CbeC^{-}_{\zeta x},\textrm{ } \widehat{\CbeC}^{-}_{\zeta},\textrm{ } \widehat{\CbeC}^{-}_{\zeta x},\textrm{ } \textrm{Jac}(\CbeC^{+}_{\zeta}),\textrm{ } \textrm{Jac}(\CbeC^{+}_{\zeta x}),\textrm{ } \textrm{Jac}(\widehat{\CbeC}^{+}_{\zeta})\textrm{ and }\textrm{Jac}(\widehat{\CbeC}^{+}_{\zeta x}).
$$ 
Moreover the assumptions of Theorem \ref{theo-critere-nullite-rang} (i.e. of Corollary \ref{stheo-3-3-final-desc-finie}) are satisfied:
\begin{itemize}
\item[* ] From Assumptions \ref{Asummption22} we get that $\omega$, $\omega^{2}-\eta^{2}$, $\omega^{2}-\eta^{2}+2\omega$ and $\omega^{2}-\eta^{2}-2\omega$ are nontrivial;
\item[* ] From Assumptions \ref{Asummption23} we know $\eta$, $\rho$, $\omega^{2}-\eta^{2}-2+2\eta$, $\omega^{2}-\eta^{2}-2-2\eta$, $\omega^{2}-\eta^{2}-1+2\eta$ and  $\omega^{2}-\eta^{2}-1-2\eta$ are nontrivial;
\item[* ] From Assumptions \ref{Asummption24} we deduce the nontriviality of $2b_{1}+\omega^{2}-\eta^{2}-1$ (see hypothesis \ref{hypothesisN}), $2b_{1}+\omega^{2}-\eta^{2}-2$ (see hypothesis \ref{hypothesisB}), $b_{1}+\eta$, $b_{1}-\eta$ (see hypothesis \ref{hypothesisM}), $b_{1}-1+\omega$ and $b_{1}-1-\omega$ (see hypothesis \ref{hypothesisF}).
\end{itemize}
From this theorem we deduce that the $\R (x)$-Mordell-Weil rank of $\textrm{Jac}(\CbeC )$ is $0$. In particular, since $\textrm{Jac}(\CbeC )$ has no antineutral torsion point (see Theorem \ref{existence-pt-anti-exemple}, notice that its hypotheses have been checked while applying Corollary \ref{finitude-jac-CbeC-final} in the proof of Corollary \ref{stheo-3-3-final-desc-finie}), the jacobian $\textrm{Jac}(\CbeC )$ has no antineutral point. To conclude we use Proposition \ref{prop-finale-HM}: since $\textrm{Jac}(\CbeC )$ has no antineutral point, $P(x^{2},y^{2})$ is not a sum of three squares in $\R (x,y)$.
\fin%

\end{section}

\addcontentsline{toc}{section}{Biliography}

\nocite{Bochnak-Coste-Roy,Bost-Mestre,Cassels83,Igusa,Knapp,Lam,Mahe1992,Mahe1990,Oort-Ueno,Pfister1967,Corps-Locaux,Gpes-Alg,Hindry-Silverman,Stoll2001,Raynaud}
\bibliography{biblio}

\end{document}